\documentclass[mathpazo]{cicp}
\usepackage{caption}
\usepackage{amsmath, amsthm, amsfonts, amssymb}
\usepackage{bm}
\usepackage{graphicx}
\usepackage{float}
\usepackage{xcolor}
\usepackage{algorithmicx,algorithm}
\usepackage{algpseudocode}
\usepackage{booktabs}
\usepackage{hyperref}
\usepackage{overpic}
\usepackage{subcaption}

\DeclareMathOperator*{\argmin}{arg\,min}

\begin{document}
\title{Integral regularization PINNs for evolution equations}

\author[Author]{Xiaodong Feng\affil{1}, Haojiong Shangguan\affil{1}\comma\corrauth, Tao Tang\affil{2,3}, Xiaoliang Wan\affil{4}}
\address{\affil{1}\ Faculty of Science and Technology, Beijing Normal-Hong Kong Baptist University, Zhuhai 519087, China. \\
    \affil{2}\ School of Mathematics and Statistics, Guangzhou Nanfang College, Guangzhou 510970, China. \\
    \affil{3}\ Zhuhai SimArk Technology Co., LTD, Zhuhai 519087, China.\\
    \affil{4}\ Department of Mathematics and Center for Computation and Technology, Louisiana State University, Baton Rouge 70803, USA.}
\email{{\tt xiaodongfeng@uic.edu.cn} (X. Feng), {\tt shangguanhaojiong@uic.edu.cn} (H. Shangguan), {\tt ttang@nfu.edu.cn} (T. Tang), {\tt xlwan@lsu.edu} (X. Wan)}

\begin{abstract}
    Evolution equations, including both ordinary differential equations (ODEs) and partial differential equations (PDEs), play a pivotal role in modeling dynamic systems. However, achieving accurate long-time integration for these equations remains a significant challenge. While physics-informed neural networks (PINNs) provide a mesh-free framework for solving PDEs, they often suffer from temporal error accumulation, which limits their effectiveness in capturing long-time behaviors. To alleviate this issue, we propose integral regularization PINNs (IR-PINNs), a novel approach that enhances temporal accuracy by incorporating an integral-based residual term into the loss function. This method divides the entire time interval into smaller sub-intervals and enforces constraints over these sub-intervals, thereby improving the resolution and correlation of temporal dynamics. Furthermore, IR-PINNs leverage adaptive sampling to dynamically refine the distribution of collocation points based on the evolving solution, ensuring higher accuracy in regions with sharp gradients or rapid variations. Numerical experiments on benchmark problems demonstrate that IR-PINNs outperform original PINNs and other state-of-the-art methods in capturing long-time behaviors, offering a robust and accurate solution for evolution equations.
\end{abstract}


\maketitle

\section{Introduction}
Evolution equations are fundamental in modeling a wide range of physical, biological, and engineering phenomena, spanning from fluid dynamics to material science \cite{courant2008methods}. Among these, evolution equations, characterized by their dependence on both temporal and spatial variables, play a crucial role in describing dynamic systems. However, achieving accurate and efficient solutions to these equations remains a significant challenge, particularly for problems requiring long-time integration or involving complex and high-dimensional domains.

In recent years, physics-informed neural networks (PINNs) \cite{raissi2019physics} have emerged as a promising alternative for solving evolution equations. By embedding the governing equations into the loss function and leveraging the expressive power of neural networks, PINNs can approximate solutions without relying on predefined grids or explicit discretization schemes. Despite their versatility, challenges remain in applying PINNs to evolution equations \cite{krishnapriyan2021characterizing, wang2021understanding, wang2022and, daw2022rethinking, rathore2024challenges}. One of the most pressing issues is the accumulation of temporal errors during long-time integration, prompting significant efforts to address this limitation. Several training strategies have been proposed to improve temporal accuracy, including sequential learning \cite{wight2020solving, mattey2022novel}, causal training \cite{wang2024respecting, penwarden2023unified, jung2024ceens} and operator learning\cite{wang2021learning, wang2023long, lippe2023pde, wu2025coast}. Additionally, hybrid strategies have been developed to combine classical numerical methods with deep learning techniques. These approaches either adapt neural networks to augment classical PDE solvers \cite{bruno2022fc, dresdner2022learning} or incorporate classical numerical methods to enhance the performance of PINNs \cite{chiu2022can, gu2022deep, feng2024hybrid}.

For evolution equations, the solution at any given time is inherently dependent on its state at previous times, reflecting strong temporal correlation. However, original PINNs treat temporal collocation points in isolation, failing to explicitly account for these correlations. This limitation often leads to challenges in capturing long-time dynamics and results in temporal error accumulation. Inspired by the integral form of evolution equations, we introduce a regularization term into the training process, proposing a novel framework termed integral regularization PINNs (IR-PINNs). Our main contributions can be summarized as:
\begin{itemize}
    \item We propose integral regularization PINNs (IR-PINNs) for evolution equations by dividing the entire time interval into smaller sub-intervals, reformulating the evolution equation into an integral form, and incorporating an integral-based residual term into the loss function. This approach enhances temporal accuracy by enforcing constraints over specific temporal subintervals, thereby improving the resolution and correlation of temporal dynamics.
    \item We extend IR-PINNs with an adaptive sampling strategy, which dynamically refines the distribution of spatial collocation points, ensuring higher accuracy in regions with sharp gradients or rapid variations.
    \item We conduct numerical experiments on benchmark problems, including both linear and nonlinear evolution equations, to demonstrate the effectiveness of IR-PINNs in capturing long-time behaviors and outperforming original PINNs and other state-of-the-art methods.
\end{itemize}

This paper is structured as follows. In Section \ref{sec:PINNs}, we provide a brief overview of PINNs, followed by a critical analysis of their limitations in solving evolution equations through a simple case study. In Section \ref{sec:methodology}, we introduce the proposed IR-PINNs framework and further enhance it by incorporating an adaptive sampling strategy to improve computational efficiency. In Section \ref{sec:numerical_experiments}, we demonstrate comprehensive numerical experiments with detailed discussions, highlighting the advantages of the proposed method. Finally, Section \ref{sec:conclusion} summarizes the paper and discusses potential future research directions.

\section{Physics-informed neural networks}\label{sec:PINNs}
We begin with a brief overview of physics-informed neural networks (PINNs) \cite{raissi2019physics}. On an open domain $\Omega\subset\mathbb{R}^d$, consider a general time-dependent PDE
\begin{equation}
    u_t(t, \bm{x}) + \mathcal{N}[u](t, \bm{x})=f(t, \bm{x}), \quad t\in[0,T], \ \bm{x}\in\Omega,
    \label{General_time_dependent_PDE}
\end{equation}
subject to the initial and boundary conditions
\begin{equation}
    \begin{aligned}
         & u(0, \bm{x}) = g(\bm{x}),\quad \bm{x}\in \Omega,                                     \\
         & \mathcal{B}[u](t, \bm{x}) = b(t, \bm{x}),\quad t\in[0,T], \ \bm{x}\in\partial\Omega,
    \end{aligned}
    \label{General_time_dependent_PDE_condition}
\end{equation}
where $\mathcal{N}[\cdot]$ is a linear or nonlinear differential operator, and $\mathcal{B}[\cdot]$ is a boundary operator corresponding to Dirichlet, Neumann, Robin or periodic boundary conditions.

Our goal is to seek a neural network $u_N(t, \bm{x}; \theta)$, where $\theta$ denotes trainable parameters, to approximate the solution $u(t, \bm{x})$ by minimizing the following composite loss function
\begin{equation}
    \small
    \begin{aligned}
        \mathcal{L}(\theta)
         & = \mathcal{L}_{r}(\theta) + \lambda_1\mathcal{L}_{ic}(\theta) + \lambda_2\mathcal{L}_{bc}(\theta)                                                                                                                                                    \\
         & = \Vert r(t, \bm{x}; \theta)\Vert_{L^2([0, T] \times \Omega)}^2 + \lambda_1\Vert u_N(0, \bm{x}; \theta)-g(\bm{x})\Vert_{L^2(\Omega)}^2 + \lambda_2\Vert \mathcal{B}[u_N](t, \bm{x}; \theta)-b(t, \bm{x})\Vert_{L^2([0, T] \times \partial\Omega)}^2,
    \end{aligned}
    \label{Continuous_loss_function}
\end{equation}
where $\lambda_1$ and $\lambda_2$ are hyper-parameters to balance the interplay between different loss terms during the training process, $\Vert \cdot \Vert_{L^2}$ represents $L^2$ norm and
\begin{equation}
    r(t, \bm{x}; \theta) = (u_N)_t(t, \bm{x}; \theta) + \mathcal{N}[u_N](t, \bm{x}; \theta) - f(t, \bm{x}).
    \label{Residual_differential}
\end{equation}
In practice, we often choose three sets of uniformly distributed collocation points  $S_r = \{(t_{r}^{(i)}, \bm{x}_{r}^{(i)})\}_{i=1}^{N_r}$on $[0,T]\times \Omega$, $S_{ic}=\{\bm{x}_{ic}^{(i)}\}_{i=1}^{N_{ic}}$ on $\Omega$ and $S_{bc} = \{(t_{bc}^{(i)}, \bm{x}_{bc}^{(i)})\}_{i=1}^{N_{bc}}$ on $[0,T]\times\partial \Omega$ and the loss function is discretized numerically as
\begin{equation}
    \begin{aligned}
        \widehat{\mathcal{L}}(\theta)
         & = \widehat{\mathcal{L}}_r(\theta) + \hat{\lambda}_1\widehat{\mathcal{L}}_{ic}(\theta) + \hat{\lambda}_2\widehat{\mathcal{L}}_{bc}(\theta)                                                                  \\
         & = \Vert r(t, \bm{x}; \theta)\Vert_{S_r}^2 + \hat{\lambda}_1\Vert u_N(0, \bm{x}; \theta)-g(\bm{x})\Vert_{S_{ic}}^2 + \hat{\lambda}_2\Vert \mathcal{B}[u_N](t, \bm{x}; \theta)-b(t, \bm{x})\Vert_{S_{bc}}^2,
    \end{aligned}
    \label{Discretized_loss_function}
\end{equation}
where
\begin{equation}
    \begin{aligned}
         & \Vert r(t, \bm{x}; \theta)\Vert_{S_r}^2 = \frac{1}{N_r}\sum_{i=1}^{N_r}\left|(u_N)_t(t_{r}^{(i)}, \bm{x}_{r}^{(i)}; \theta) + \mathcal{N}[u_N](t_{r}^{(i)}, \bm{x}_{r}^{(i)}; \theta) - f(t_{r}^{(i)}, \bm{x}_{r}^{(i)})\right|^2, \\
         & \Vert u_N(0, \bm{x}; \theta)-g(\bm{x})\Vert_{S_{ic}}^2 = \frac{1}{N_{ic}}\sum_{i=1}^{N_{ic}}\left|u_N(0, \bm{x}_{ic}^{(i)}; \theta) - g(\bm{x}_{ic}^{(i)})\right|^2,                                                               \\
         & \Vert \mathcal{B}[u_N](t, \bm{x}; \theta)-b(t, \bm{x})\Vert_{S_{bc}}^2 = \frac{1}{N_{bc}}\sum_{i=1}^{N_{bc}}\left|\mathcal{B}[u_N](t_{bc}^{(i)}, \bm{x}_{bc}^{(i)};\theta) - b(t_{bc}^{(i)}, \bm{x}_{bc}^{(i)})\right|^2,
    \end{aligned}
    \label{Different_discretized_loss_terms}
\end{equation}
and $\hat{\lambda}_1 = \frac{\lambda_1}{T}, \hat{\lambda}_2 = \frac{\lambda_2\vert \partial\Omega\vert}{\vert \Omega\vert}$ such that $\widehat{\mathcal{L}}(\theta)$ is a Monte Carlo approximation of $\mathcal{L}(\theta)$ up to a constant scaling factor $T\vert \Omega\vert$.

To provide intuition for the key ideas behind our method, we first analyze a simple case by considering the following ordinary differential equation (ODE):

\begin{equation}
    \begin{aligned}
         & u'(t) = \lambda u(t), \quad t\in [0, T], \\
         & u(0) = u_0,
    \end{aligned}
\end{equation}
where the exact solution is $u(t) = u_0e^{\lambda t}$. Conventional PINNs are known to struggle with large values of $\lambda T$ since the ODE system will become stiff. Here we set $\lambda = 1$, $T = 5$ and $u_0 = 1$. In the framework of PINNs, we represent the latent variable $u(t)$ by a fully-connected neural network $u_N(t; \theta)$ with Tanh activation function, $3$ hidden layers, and $64$ neurons per hidden layer, and exactly impose the initial condition by
\begin{equation}
    \tilde{u}_N(t; \theta) = tu_N(t; \theta) + u_0.
\end{equation}
Then the discretized loss function simplifies to
\begin{equation}
    \widehat{\mathcal{L}}(\theta) = \widehat{\mathcal{L}}_r(\theta) = \Vert r(t; \theta)\Vert_{S_r}^2 = \frac{1}{N_r}\sum_{i=1}^{N_r}\left|\tilde{u}_N'(t_{r}^{(i)}; \theta) - \lambda\tilde{u}_N(t_{r}^{(i)}; \theta)\right|^2,
\end{equation}
where $S_r = \{t_r^{(i)}\}_{i=1}^{N_r}$ are uniform mesh points on $[0, T]$ and $N_r = 40$. The Adam algorithm is used to optimize the loss function with an initial learning rate $10^{-3}$, which decays every $1,000$ epochs at a rate of $0.9$. The maximum number of training epochs is $80,000$.
\begin{figure}[H]
    \centering
    \includegraphics[width=0.33\linewidth]{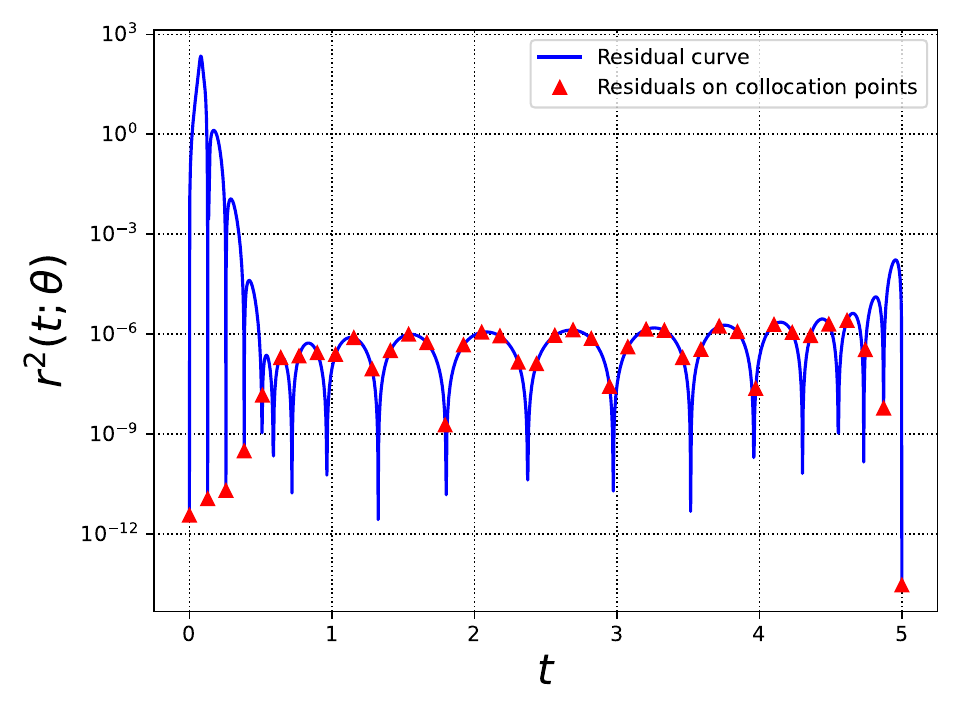}
    \includegraphics[width=0.33\linewidth]{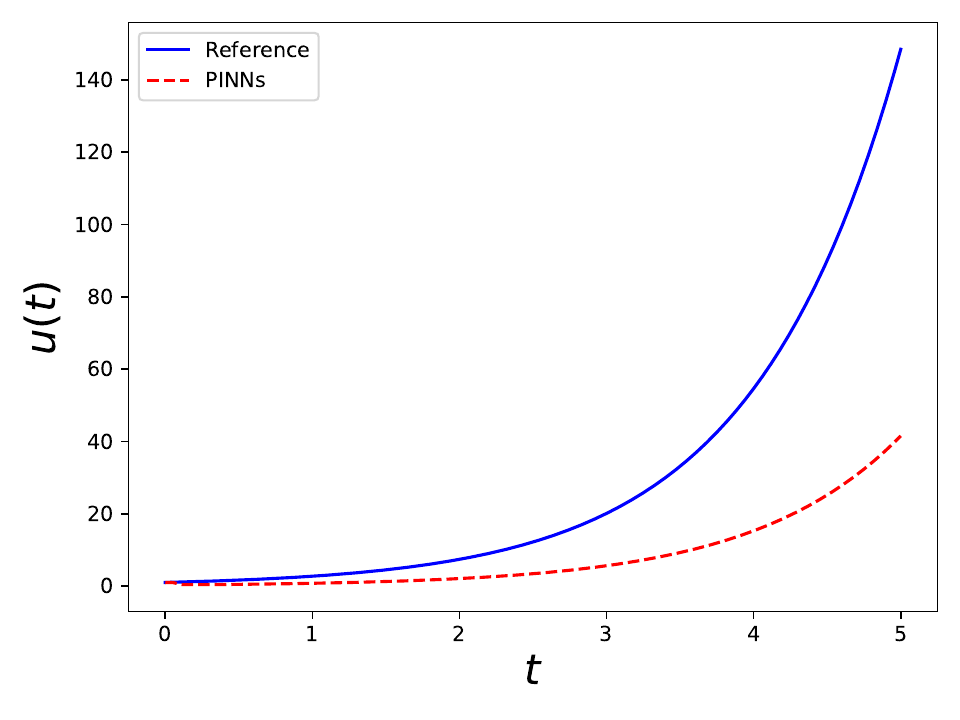}
    \caption{\textit{Simple ODE:} Left: Residual curve. Right: Reference solution versus numerical solution. Relative $L_2$ error: $7.2033e-01$.}
    \label{fig:Simple_ODE_PINNs_Uniform}
\end{figure}
As shown in Figure \ref{fig:Simple_ODE_PINNs_Uniform}, although the residuals on the training set are sufficiently small, conventional PINNs fail to capture the exponential growth of the solution. This suggests that the model tends to overfit when the temporal collocation points are insufficiently dense.

To investigate the failure mode of PINNs in this scenario, we analyze the role of the automatic differentiation technique. While automatic differentiation can accurately compute the derivative at individual temporal collocation points, it does not account for the information in the small neighborhood around each collocation point. This may lead to relatively large residuals in regions outside the training points. From the residual curve, it is evident that when $t < 0.5$, the overall residual is large. Given the propagation property of the evolution equation, PINNs cannot converge to the correct solution, even if the residual later becomes small. This analysis suggests that when solving long-time integration problems, it is crucial to consider the correlations between collocation points to reduce the residuals across the entire temporal interval,  a principle that traditional numerical schemes often emphasize. In the next section, we propose a simple and effective strategy based on integral form of the evolution equation to strengthen the temporal correlation and mitigate the overfitting problem.

\section{Methodology}\label{sec:methodology}
\subsection{Integral regularization PINNs (IR-PINNs)}
We propose a novel method for solving evolution equations, termed integral regularization PINNs (IR-PINNs), which enhances the standard loss function with a regularization term designed to inherently capture temporal correlation in the solution.

More precisely, we consider the integral form of the general time-dependent PDE \eqref{General_time_dependent_PDE} over a temporal interval $[a, b]$ ($a < b$):
\begin{equation}
    u(b, \bm{x}) - u(a, \bm{x}) + \int_{a}^{b}\mathcal{N}[u](t, \bm{x})\mathrm{d} t = \int_{a}^{b}f(t, \bm{x})\mathrm{d} t, \quad \forall a, b\in[0,T], \ \bm{x}\in\Omega.
    \label{Integral_equation}
\end{equation}
Then we can define an integral-based residual function corresponding to the above equation \eqref{Integral_equation} as follows:
\begin{equation}
    r_{\rm int}(a, b, \bm{x}) = \frac{u(b, \bm{x}) - u(a, \bm{x}) + \int_{a}^{b}\left(\mathcal{N}[u](t, \bm{x}) - f(t, \bm{x})\right)\mathrm{d} t}{b-a}.
    \label{Residual_integral}
\end{equation}
Given a numerical quadrature scheme defined by a set of nodes $Q = \{s_i\}_{i=1}^{M}$ and their associated weights $\{w_i\}_{i=1}^{M}$, the integral term can be approximated as follows
\begin{equation}
    \begin{aligned}
        I_{Q}[\mathcal{N}[u](t, \bm{x}) - f(t, \bm{x})] = \sum_{i=1}^{M}w_i\left(\mathcal{N}[u](s_i, \bm{x}) - f(s_i, \bm{x})\right).
    \end{aligned}
\end{equation}

We divide the temporal domain $[0, T]$ into $N$ equidistant subintervals, denoted by $0 = t_0 < t_1 < t_2 < \cdots < t_N = T$ and two distinct options for the regularization term are provided:
\begin{itemize}
    \item (IR-PINNs1) $\mathcal{L}_{\rm int}^{(\rm R1)}(\theta) = \frac{1}{N}\sum_{k=0}^{N-1}\Vert r_{\rm int}(t_k, t_{k+1}, \bm{x}; \theta)\Vert_{L^2(\Omega)}^2.$
    \item (IR-PINNs2) $\mathcal{L}_{\rm int}^{(\rm R2)}(\theta) = \frac{1}{N}\sum_{k=0}^{N-1}\Vert r_{\rm int}(t_0, t_{k+1}, \bm{x}; \theta)\Vert_{L^2(\Omega)}^2.$
\end{itemize}

Given numerical quadrature points $Q_r^{(k)} = \{t_r^{(k,i)}\}_{i=1}^M$ for each subinterval $[t_{k}, t_{k+1}]$, the total numerical quadrature points set is $Q = \cup_{k=0}^{N-1}Q_r^{(k)}$ and the new residual loss function can be written as
\begin{equation}
    \mathcal{L}_{\rm couple}(\theta) = \mathcal{L}_{r}(\theta) + \alpha\mathcal{L}_{\rm int}(\theta),
    \label{New_residual_loss}
\end{equation}
where $\mathcal{L}_{r}(\theta)$ is modified to
\begin{equation}
    \mathcal{L}_r(\theta) = \frac{1}{NM}\sum_{k=0}^{N-1}\sum_{i=1}^{M}\Vert r(t_{r}^{(k, i)}, \bm{x}; \theta)\Vert_{L^2(\Omega)}^2,
\end{equation}
$\mathcal{L}_{\rm int}(\theta)$ stands for $\mathcal{L}_{\rm int}^{(\rm R1)}(\theta)$ or $\mathcal{L}_{\rm int}^{(\rm R2)}(\theta)$ and $\alpha$ is a hyper-parameter to balance the two residual loss terms.

In practice, we generate uniformly distributed points $\{\bm{x}_r^{(j)}\}_{j=1}^{N_{\bm{x}}}$ in spatial direction and obtain two sets of collocation points $S_{r} = \cup_{k=0}^{N-1}\cup_{i=1}^{M}S_{r}^{(k, i)} = \cup_{k=0}^{N-1}\cup_{i=1}^{M}\{(t_{r}^{(k, i)}, \bm{x}_{r}^{(j)})\}_{j=1}^{N_{\bm{x}}}$ and $S_{\rm int} = \cup_{k=0}^{N-1}S_{\rm int}^{(k)} = \cup_{k=0}^{N-1}\{\bm{x}_{r}^{(j)}\}_{j=1}^{N_{\bm{x}}}$, the discretized residual loss function is
\begin{equation}
    \widehat{\mathcal{L}}_{\rm couple}(\theta) = \widehat{\mathcal{L}}_{r}(\theta) + \hat{\alpha}\widehat{\mathcal{L}}_{\rm int}(\theta)
\end{equation}
where
\begin{gather}
    \widehat{\mathcal{L}}_{r}(\theta) = \frac{1}{NM}\sum_{k=0}^{N-1}\sum_{i=1}^{M}\Vert r(t_{r}^{(k, i)}, x; \theta)\Vert_{S_{r}^{(k, i)}}^2, \\
    \widehat{\mathcal{L}}_{\rm int}^{(\rm R1)}(\theta) = \frac{1}{N}\sum_{k=0}^{N-1}\Vert r_{\rm int}(t_k, t_{k+1}, \bm{x}; \theta)\Vert_{S_{\rm int}^{(k)}}^2, \quad \widehat{\mathcal{L}}_{\rm int}^{(\rm R2)}(\theta) = \frac{1}{N}\sum_{k=0}^{N-1}\Vert r_{\rm int}(t_0, t_{k+1}, \bm{x}; \theta)\Vert_{S_{\rm int}^{(k)}}^2.
    \label{Discretized_residual_loss}
\end{gather}
More precisely,
\begin{equation}
    \Vert r(t_{r}^{(k, i)}, x; \theta)\Vert_{S_{r}^{(k, i)}}^2 = \frac{1}{N_{\bm{x}}}\sum_{j=1}^{N_{\bm{x}}}\left|(u_N)_t(t_{r}^{(k, i)}, \bm{x}_{r}^{(j)}; \theta) + \mathcal{N}[u_N](t_{r}^{(k, i)}, \bm{x}_{r}^{(j)}; \theta) - f(t_{r}^{(k, i)}, \bm{x}_{r}^{(j)})\right|^2,
    \label{Differential_residual_loss}
\end{equation}
for all $0 \le k \le N-1, 1 \le i \le M$ and
\begin{equation}
    \begin{aligned}
          & \Vert r_{\rm int}(t_m, t_n, \bm{x}; \theta)\Vert_{S_{\rm int}^{(k)}}^2                                                                                                                                                                               \\
        = & \frac{1}{N_{\bm{x}}}\sum_{j=1}^{N_{\bm{x}}}\left|\frac{u_N(t_n, \bm{x}_r^{(j)}; \theta) - u_N(t_m, \bm{x}_r^{(j)}; \theta) + \sum_{k=m}^{n-1}I_{Q_{r}^{(k)}}[\mathcal{N}[u_N](t, \bm{x}_r^{(j)}; \theta) - f(t, \bm{x}_r^{(j)})]}{t_n-t_m}\right|^2,
    \end{aligned}
    \label{Integral_residual_loss}
\end{equation}
for all $0 \le m < n \le N$. Here $\hat{\alpha} = \alpha$ such that $\widehat{\mathcal{L}}_{\rm couple}(\theta)$ is a Monte Carlo approximation of $\mathcal{L}_{\rm couple}(\theta)$ up to a constant scaling factor $\vert \Omega\vert$. Finally we obtain an estimator $\theta^{\star}$ by solving the following optimization problem
\begin{equation}
    \theta^{\star} = \argmin_{\theta}\widehat{\mathcal{L}}(\theta) = \argmin_{\theta}\left(\widehat{\mathcal{L}}_{\rm couple}(\theta) + \hat{\lambda}_1\widehat{\mathcal{L}}_{ic}(\theta) + \hat{\lambda}_2\widehat{\mathcal{L}}_{bc}(\theta)\right).
    \label{Optimization_problem}
\end{equation}

We now proceed to implement IR-PINNs to solve the simple ODE problem described in the preceding section. The network architecture, number of collocation points, training parameters, maximum number of iterations, and optimizer remain unchanged. The principle enhancement involves the incorporation of an integral-based regularization term into the loss function. Here we set $\alpha = 1$ and employ the trapezoidal rule to approximate the integral in the regularization term.
\begin{figure}[H]
    \centering
    \includegraphics[width=0.33\linewidth]{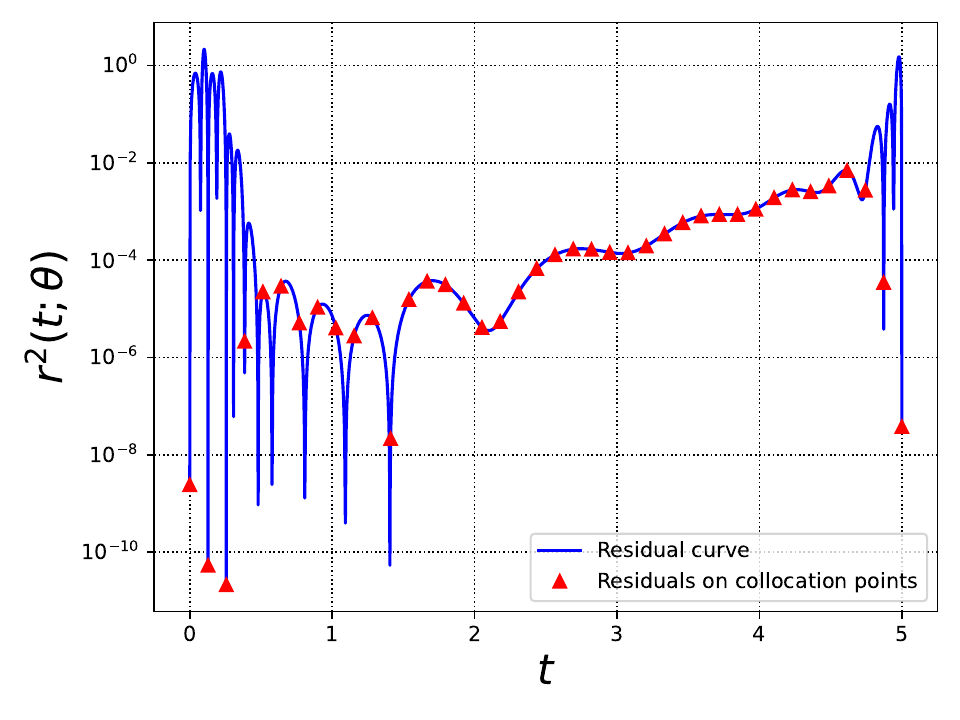}
    \includegraphics[width=0.33\linewidth]{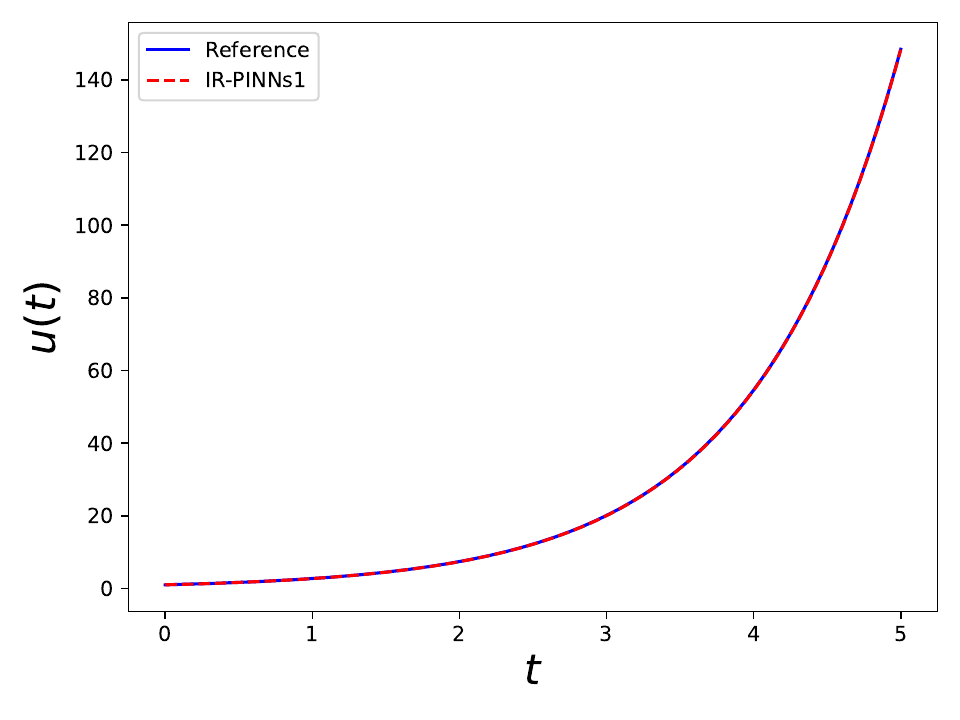}
    \caption{\textit{Simple ODE:} Left: Residual curve. Right: Reference solution versus numerical solution. Relative $L_2$ error: $1.3156e-03$.}
    \label{fig:Simple_ODE_IR-PINNs1}
\end{figure}
\begin{figure}[H]
    \centering
    \includegraphics[width=0.33\linewidth]{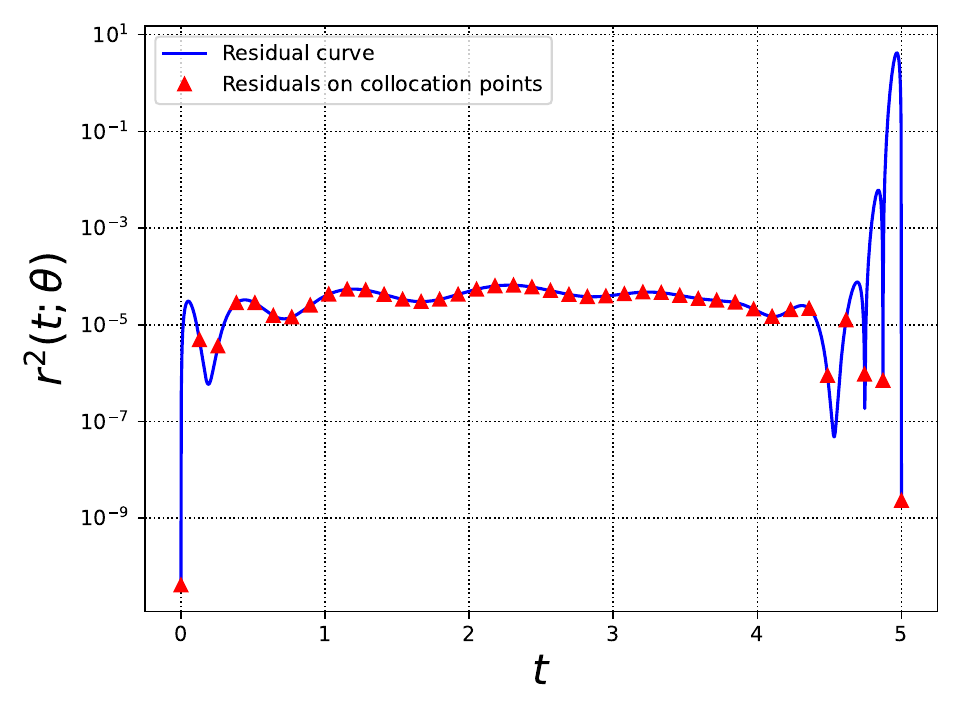}
    \includegraphics[width=0.33\linewidth]{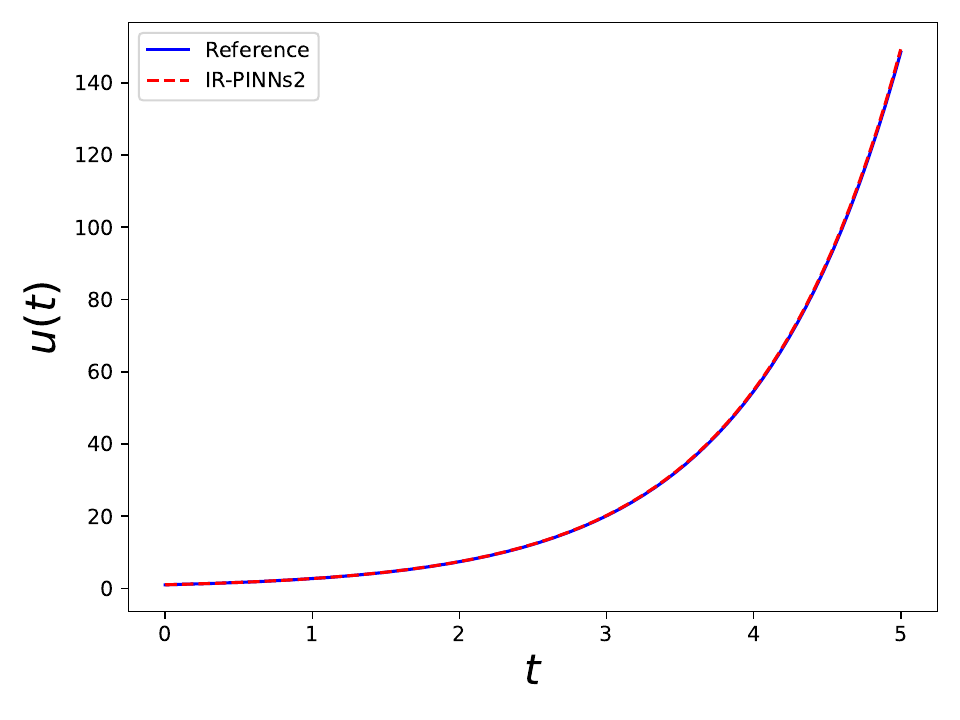}
    \caption{\textit{Simple ODE:} Left: Residual curve. Right: Reference solution versus numerical solution. Relative $L_2$ error: $4.9012e-03$.}
    \label{fig:Simple_ODE_IR-PINNs2}
\end{figure}
The results of IR-PINNs1 and IR-PINNs2 are summarized in Figure \ref{fig:Simple_ODE_IR-PINNs1} and Figure \ref{fig:Simple_ODE_IR-PINNs2}, respectively. One can see that the predicted solutions achieve an excellent agreement with the reference solution, yielding a reduction in relative $L_2$ error of two orders of magnitude compared to PINNs. This improvement can be primarily attributed to the incorporation of the regularization term, which inherits the benefits of traditional methods by strengthening the correlations among collocation points, effectively reducing residuals, particularly in the earlier time periods. This ensures proper temporal propagation of the governing equation's information during the training process. It is noteworthy that, due to the exponential growth nature of the solution, the residual range of the numerical solution gradually increases over time, which is entirely expected.

\paragraph{The choice of number of subintervals $N$:}
The selection of $N$  significantly influences both the numerical precision and the efficacy of the regularization process. With a fixed number of quadrature points for each subinterval, a larger $N$ is essential for accurate numerical integration, particularly when the temporal interval is long or the solution exhibits complex behaviors. Additionally, $N$ determines the number of integral-based residual constraints, and increasing it appropriately enhances the impact of regularization term on the training process, thereby improving the ability of model to capture temporal correlation.

\paragraph{The choice of hyperparameter $\alpha$:}
The the hyperparameter $\alpha$ determines the balance between the two residual loss terms $\mathcal{L}_{r}(\theta)$ and $\mathcal{L}_{\rm int}(\theta)$ and plays a critical role in the training process of the neural network. If $\alpha$ is too small, the method degenerates to the original PINNs, where the regularization term has almost no effect on the training process. On the other hand, if $\alpha$ is too large, the regularization term dominates the training process. However, due to limitations in integral accuracy and optimization difficulty, this can lead to unstable training and hinder the achievement of desired accuracy. Therefore, selecting an appropriate $\alpha$ is crucial. Here we propose a simple strategy for selecting $\alpha$. From the mean value theorem, there exists some $\xi \in [a, b]$ such that
\begin{equation}
    \begin{aligned}
        r_{\rm int}(a, b, \bm{x})
         & = \frac{u(b, \bm{x}) - u(a, \bm{x}) + \int_{a}^{b}(\mathcal{N}[u](t, \bm{x}) - f(t, \bm{x}))\mathrm{d} t}{b - a} \\
         & = \frac{\int_{a}^{b}(u_t(t, \bm{x}) + \mathcal{N}[u](t, \bm{x}) - f(t, \bm{x}))\mathrm{d} t}{b - a}              \\
         & = u_t(\xi, \bm{x}) + \mathcal{N}[u](\xi, \bm{x}) - f(\xi, \bm{x})                                                \\
         & = r(\xi, \bm{x}).
    \end{aligned}
\end{equation}
This implies that the two residuals $r$ and $r_{\rm int}$ are of similar order of magnitude, suggesting $\alpha = 1$ as a natural choice. Certainly, a more effective approach would be to adaptively adjust the weight $\alpha$, as discussed in \cite{wang2021understanding, wang2022and}. However, this is not the primary focus of the current work. Therefore, for all numerical experiments conducted in this study, we consistently set $\alpha = 1$ as a fixed value.

\paragraph{Computational cost of regularization term:}
While introducing a new regularization term, it is essential to consider its computational cost, particularly in the context of numerical quadrature. From equations \eqref{Differential_residual_loss} and \eqref{Integral_residual_loss}, we observe that the set of residual collocation points $S_{r}$ is used for discretizing the two residual loss terms. Before each update of the loss function, we can compute and store the following values:
\begin{equation}
    \mathcal{N}[u_N](t_{r}^{(i)}, \bm{x}_{r}^{(j)}; \theta) - f(t_{r}^{(i)}, \bm{x}_{r}^{(j)}), \quad \forall i = 1, \cdots, N_t, \ j = 1, \cdots, N_{\bm{x}}.
    \label{Repeated_values}
\end{equation}
These precomputed values can then be reused to evaluate the two residual loss terms, minimizing redundant computations. As a result, the additional computational cost is negligible and remains within acceptable limits.

\begin{figure}[H]
    \centering
    \includegraphics[width=1.0\linewidth]{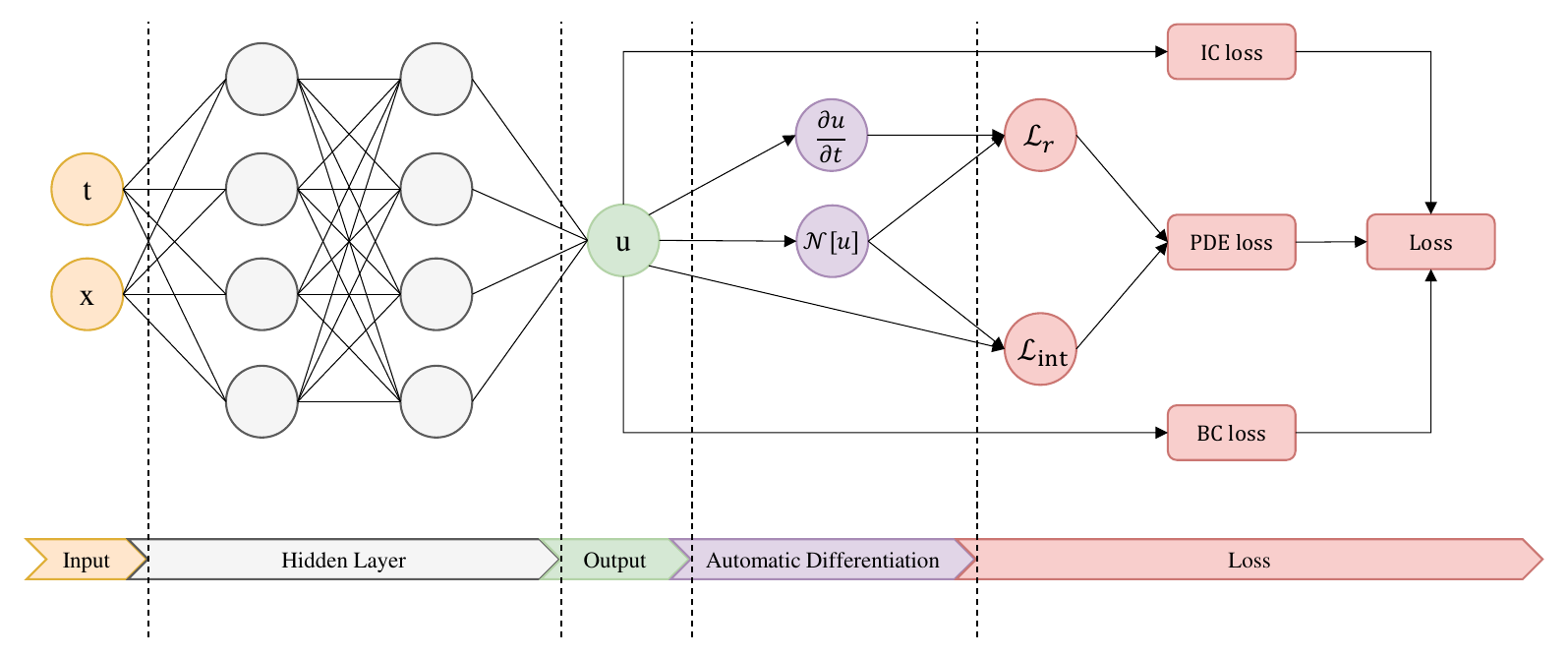}
    \caption{Schematic diagram of integral regularization PINNs.}
    \label{fig:Schematic_diagram}
\end{figure}

The schematic diagram of the proposed approach is shown in Figure \ref{fig:Schematic_diagram}, and the corresponding algorithm is summarized in Algorithm \ref{alg:1}.

\begin{algorithm}
    \caption{Solving evolution equations via integral regularization PINNs}
    \label{alg:1}
    \begin{algorithmic}
        \State \textbf{Input:} Terminal time $T$, number of subintervals $N$, numerical quadrature points set $Q$, number of training epochs $N_e$, initial learning rate $l_r^{(0)}$, decay rate $\eta$, step size $n_s$.
        \State Generate the training datasets $S_r$, $S_{\rm int}$, $S_{ic}$ and $S_{bc}$.
        \State Set $l_r$ to $l_r^{(0)}$.
        \For{$e = 1,\cdots,N_e$}
        \If{$e \% n_s==0$}
        \State $l_r=\eta*l_r$.
        \EndIf
        \State Compute the shared values in \eqref{Repeated_values} using the datasets $S_r$ and $S_{\rm int}$.
        \State Compute the residual loss function $\widehat{\mathcal{L}}_{r, \rm new}(\theta)$ using the datasets $S_r$, $S_{\rm int}$ and shared values.
        \State Compute the initial and boundary loss functions $\widehat{\mathcal{L}}_{ic}(\theta)$ and $\widehat{\mathcal{L}}_{bc}(\theta)$ using the dataset $S_{ic}$ and $S_{bc}$.
        \State Compute the total loss function ${\widehat{\mathcal{L}}}(\theta)$ in equation \eqref{Optimization_problem}.
        \State Update $\theta$ by stochastic gradient descent method.
        \EndFor
        \State \textbf{Output:} The predicted solutions $u_N(t,\bm{x};\theta^*)$.
    \end{algorithmic}
\end{algorithm}

\subsection{Adaptive sampling strategy for IR-PINNs}
In this section, we present an extension of the integral regularization PINNs (IR-PINNs) by incorporating the adaptive sampling strategy proposed in \cite{feng2024hybrid}. This extension further enhances the performance of IR-PINNs when addressing more complex and challenging scenarios, such as low-regularity problems or unbounded domain problems.

The basic idea is to construct a continuous joint probability density model $p(t, \bm{x}; \theta_f)$ to approximate the distribution induced by the residual $r^2(t, \bm{x}; \theta)$ in \eqref{Residual_differential} as
\begin{equation}
    p(t, \bm{x}; \theta_f) = p_{\mathrm{poly}}(t;\theta_{f,1}) p_{\mathrm{B-KRnet}}(\bm{x}|t;\theta_{f,2}),
\end{equation}
where $p_{\mathrm{poly}}(t;\theta_{f,1})$ is the bounded polynomial layer \cite{muller2019neural}, $p_{\mathrm{B-KRnet}}(\bm{x}|t;\theta_{f,2})$ is the bounded KR-net \cite{zeng2023bounded} and $\theta_f = \{\theta_{f,1}, \theta_{f,2}\}$. To seek the optimal parameter $\theta_f$, we can minimize the following objective
\begin{equation}
    \begin{aligned}
          & D_{\mathrm{KL}}(\hat{r}(t,\bm{x};\theta) || p(t,\bm{x}; \theta_f))                                    \\
        = & D_{\mathrm{KL}}\big(\hat{r}(t,\bm{x};\theta)|| p_{\mathrm{poly}}
        (t;\theta_{f,1})p_{\mathrm{B-KRnet}}(\bm{x}|t;\theta_{f,2})\big)                                          \\
        = & \iint \hat{r}(t,\bm{x};\theta) \log \left(\hat{r}(t,\bm{x};\theta)\right) \mathrm{d}\bm{x}\mathrm{d}t
        - \iint \hat{r}(t,\bm{x};\theta) \log \left( p_{\mathrm{poly}}(t;\theta_{f,1})p_{\mathrm{B-KRnet}}(\bm{x}|t;\theta_{f,2})\right) \mathrm{d}\bm{x}\mathrm{d}t,
    \end{aligned}
    \label{Density_kl}
\end{equation}
where $D_{\mathrm{KL}}$ indicates the Kullback-Leibler (KL) divergence and $\hat{r}(t,\bm{x};\theta)$ is the distribution induced by $r^2(t,\bm{x};\theta)$ and the well-trained parameters $\theta_f^* = \{\theta_{f,1}^*, \theta_{f,2}^*\}$ is obtained. The first term on the right-hand side in \eqref{Density_kl} represents the negative differential entropy of
$\hat{r}(t,\bm{x}; \theta)$, which does not affect the optimization process with respect to
$\theta_f$. So minimizing the KL divergence is equivalent to minimizing the cross
entropy between $\hat{r}(t,\bm{x}; \theta)$ and $p(t,\bm{x};\theta_f)$
\begin{equation}
    \begin{aligned}
        H(\hat {r}(t,\bm{x}; \theta), p(t,\bm{x};\theta_f)) &
        = - \iint \hat{r}(t, \bm{x};\theta) \log \left( p_{\mathrm{poly}}(t;\theta_{f,1})
        p_{\mathrm{B-KRnet}}(\bm{x}|t;\theta_{f,2})\right) \mathrm{d}\bm{x}\mathrm{d}t,
    \end{aligned}
\end{equation}
Since the samples from $\hat{r}(t,\bm{x};\theta)$ are not available, we approximate the cross entropy up to a constant scaling factor $T\vert\Omega\vert$ using the importance sampling technique:
\begin{equation}
    \begin{aligned}
                & H(\hat{r}(t,\bm{x};\theta),p(t,\bm{x};\theta_f)) \\
        \approx & -\frac{1}{N_r}\sum_{i=1}^{N_r}
        \frac{\hat{r}(t_i,\bm{x}_i;\theta)}{p_{\mathrm{poly}}(t_i;\hat{\theta}_{f,1})
            p_{\mathrm{B-KRnet}}(\bm{x}_i|t_i;\hat{\theta}_{f,2})}
        \left(\log p_{\mathrm{poly}}(t_i;\theta_{f,1}) + \log p_{\mathrm{B-KRnet}}
        (\bm{x}_i|t_i;{\theta}_{f,2})\right),
    \end{aligned}
    \label{cross_entropy}
\end{equation}
where
\begin{equation}
    t_i\sim p_{\mathrm{poly}}(\cdot;\hat{\theta}_{f,1}),\quad  \bm{x}_i\sim p_{\mathrm{B-KRnet}}(\cdot|t_i;\hat{\theta}_{f,2}).
    \label{sample_from_joint_distribution}
\end{equation}
More details about the adaptive sampling approach can be found in \cite{feng2024hybrid}.

\begin{figure}[H]
    \centering
    \includegraphics[width=1.0\linewidth]{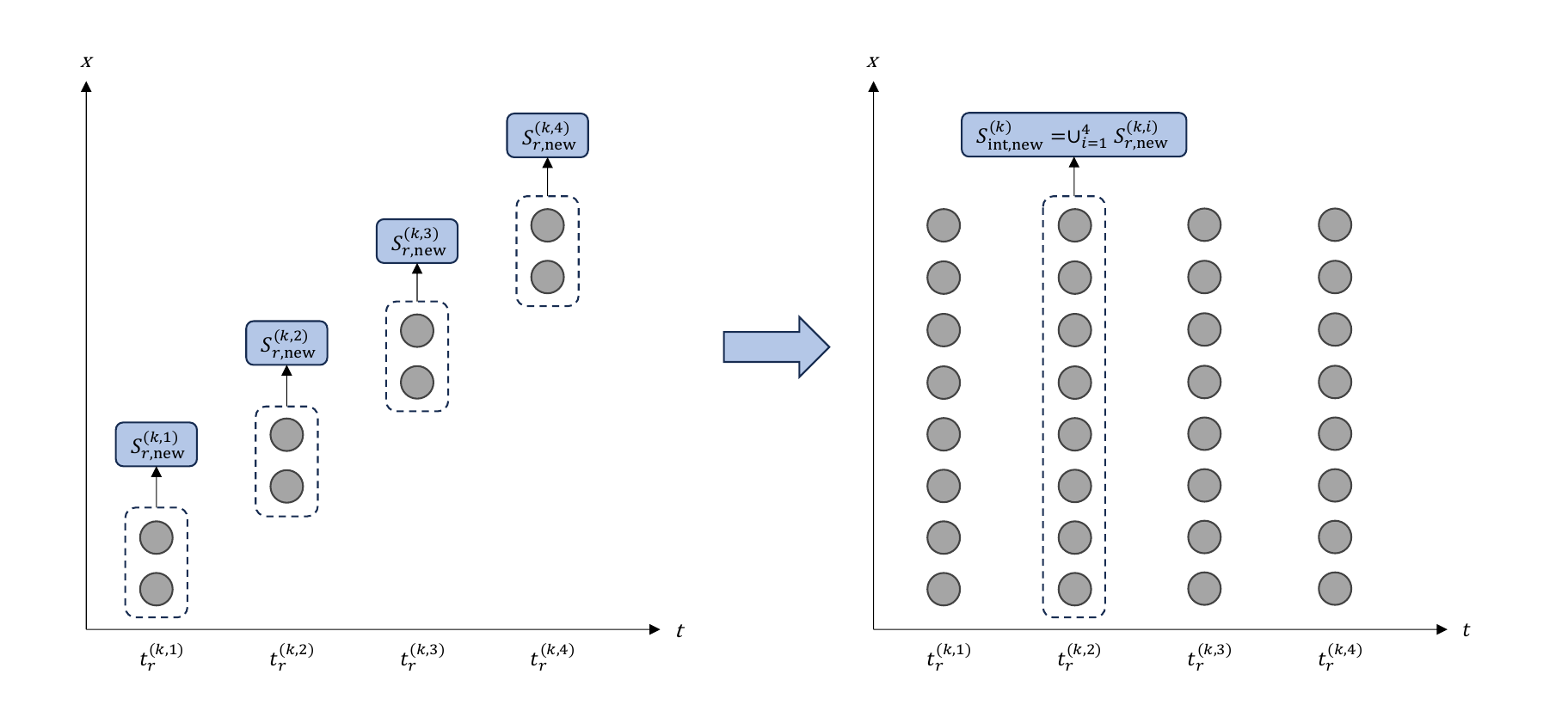}
    \caption{Schematic diagram of adaptive sampling.}
    \label{fig:Schematic_diagram_adap}
\end{figure}

Considering that the residuals $r^2(t, \bm{x};\theta)$ are computed only at the discrete temporal points set $Q = \cup_{k=0}^{N-1}Q_r^{(k)}$, we employ a discrete distribution $p_{\rm dis}(t; \theta_{f,1}^*)$ based on $p_{\rm poly}(t; \theta_{f,1})$. Subsequently, we generate $N_{\rm new}$ samples $\{t^{(j)}\}_{j=1}^{N_{\rm new}}$ from the discrete distribution $p_{\mathrm{dis}}(t; \theta_{f,1})$, along with their corresponding spatial points $\bm{x}^{(j)}\sim p_{\mathrm{B-KRnet}} (\bm{x}|t^{(j)};\theta_{f,2})$ for each $t^{(j)}$. Finally, we reorder the newly generated data into the original spatial domain
\begin{align}
     & \bm{x}^{(j)} \in S_{r, \rm new}^{(k, i)}, \quad \mathrm{if} \ t^{(j)} = t_r^{(k, i)}, \forall j=1,\cdots,N_{\rm new}, \ k=0,\cdots,N-1, \ i=1,\cdots,M, \\
     & \bm{x}^{(j)} \in S_{\rm int, \rm new}^{(k)}, \quad \mathrm{if} \ t^{(j)} \in T^{(k)}, \forall j=1,\cdots,N_{\rm new}, \ k=0,\cdots,N-1,
\end{align}
where $T^{(k)} = Q_r^{(k)}$ for IR-PINNs1 and $T^{(k)} = \cup_{p=0}^{k}Q_r^{(p)}$ for IR-PINNs2. Now we refine the training set for the residual loss function $\widehat{\mathcal{L}}_{r, \rm new}(\theta)$ as follows:
\begin{align}
     & S_{r} = \cup_{k=0}^{N-1}\cup_{i=1}^{M}(S_{r}^{(k, i)}\cup S_{r, \rm new}^{(k, i)}), \label{NewSr} \\
     & S_{\rm int} = \cup_{k=0}^{N-1}S_{\rm int, \rm new}^{(k)}.\label{NewSint}
\end{align}
Then we continue training with the updated training datasets $S_{r}$ and $S_{\rm int}$ until the training is terminated and the final solution is obtained. A schematic diagram representing adaptive sampling is shown in Figure \ref{fig:Schematic_diagram_adap}. The detailed algorithm is summarized in Algorithm \ref{alg:2}.

\begin{algorithm}
    \caption{Adaptive sampling strategy for IR-PINNs}
    \label{alg:2}
    \begin{algorithmic}
        \State \textbf{Input:} Number of adaptive iteration $N_{\mathrm{adaptive}}$, number of adaptive training epochs $N_a$, number of newly added points $N_{\mathrm{new}}$, initial training datasets $S_{r}^{(0)}$, $S_{\rm int}^{(0)}$, $S_{ic}$ and $S_{bc}$, initial probability density model $p(t, \bm{x};\theta_{f}^{*, (0)})$.
        \State Solve evolution equation via Algorithm \ref{alg:1} to obtain $u(t,\bm{x};\theta^{*,(0)})$.
        \For{$k=0,\cdots,N_{\mathrm{adaptive}}-1$}
        \State // Train probability density model.
        \For {$j=1,\cdots,N_a$}
        \State Generate samples from $p(t, \bm{x};\theta_{f}^{*, (k)})$.
        \State Update $p(t, \bm{x};\theta_{f}^{*, (k+1)})$ by descending the stochastic gradient of \eqref{cross_entropy}.
        \EndFor
        \State // Refine training datasets.
        \State Generate new training samples $S_{r, \rm new}$ through $p(t, \bm{x};\theta_{f}^{*, (k+1)})$.
        \State Update training datasets $S_{r}^{(k+1)}$ and $S_{\rm int}^{(k+1)}$ based on \eqref{NewSr} and \eqref{NewSint}.
        \State Solve evolution equation via Algorithm \ref{alg:1} to obtain $u(t,\bm{x};\theta^{*,(k+1)})$.
        \EndFor
        \State \textbf{Output:} The predicted solution $u_N(t,\bm{x};\theta^*)$.
    \end{algorithmic}
\end{algorithm}

\section{Numerical experiments}\label{sec:numerical_experiments}
In this section, we conduct some numerical experiments to demonstrate the effectiveness of the proposed method. Specifically, we will approximate the Lorentz system, the Kuramoto-Sivashinsky equation, the Boussinesq-Burgers equations, and the nonlinear time-dependent Fokker-Planck equation. Throughout all benchmarks, we will employ the fully-connected neural network equipped with hyperbolic tangent activation functions (Tanh) and initialized using the Glorot normal scheme \cite{glorot2010understanding}. All neural networks are trained via full-batch gradient descent using the Adam optimizer with default settings \cite{kingma2014adam}. For the approximation of the integral in the regularization term, we employ Gaussian quadrature formulas. Considering the highly chaotic behavior of the solutions in the first two examples, we utilize 64-bit double-precision floating-point numbers in our numerical simulations. All experiments are implemented by JAX \cite{bradbury2018jax}.

In order to test the validity of the method, we use the following relative $L_2$ error:
\begin{equation}
    err_{L_2} = \frac{\sqrt{\sum_{i=1}^{N_\mathrm{test}}|u_N(t_i, \bm{x}_i;\theta) - u(t_i, \bm{x}_i)|^2}}{\sqrt{\sum_{i=1}^{N_\mathrm{test}}|u(t_i, \bm{x}_i)|^2}},
\end{equation}
where $N_\mathrm{test}$ represents the total number of test points chosen randomly in the domain, and $u_N(t_i, \bm{x}_i;\theta)$ and $u(t_i, \bm{x}_i)$ represent the numerical and the reference solution values, respectively.

\subsection{Lorentz system}
We start with the chaotic Lorentz system, an ODE system. It is well known that this system exhibits strong sensitivity to its initial conditions, which can trigger divergent trajectories in finite time if the numerical predictions sought are not sufficiently accurate. The system is described by the following ordinary differential equations
\begin{equation}
    \begin{aligned}
         & \frac{\mathrm{d} x}{\mathrm{d} t}=\sigma(y-x), \\
         & \frac{\mathrm{d} y}{\mathrm{d} t}=x(\rho-z)-y, \\
         & \frac{\mathrm{d} z}{\mathrm{d} t}=x y-\beta z.
    \end{aligned}
\end{equation}
Here, we consider a classical setting with $\sigma=3$, $\rho=28$ and $\beta=8/3$. The temporal domain is $[0, 20]$ and initial conditions are $[x(0), y(0), z(0)] = [1, 1, 1]$.

In this example, we employ the time-marching strategy and partition the temporal domain $[0, 20]$ into $40$ subdomains. The details of the time-marching strategy can be founded in Appendix \ref{sec:time_marching}. For each subdomain, we construct a $5$-layer fully-connected neural network with $512$ neurons per hidden layer. The initial learning rate is $10^{-3}$, with an exponential decay rate of $0.9$ applied every $2,000$ training epochs. The maximum number of training epochs is $150,000$. We divide the subdomain into $N = 64$ subintervals and in each subinterval, the number of Gaussian quadrature points is $M = 4$. Since Lorentz system is highly sensitive to the initial condition, we exactly impose the initial condition by
\begin{align}
    \tilde{x}_{N}(t; \theta) & = t x(t; \theta) + x(0), \\
    \tilde{y}_{N}(t; \theta) & = t y(t; \theta) + y(0), \\
    \tilde{z}_{N}(t; \theta) & = t z(t; \theta) + z(0).
\end{align}

\begin{figure}[H]
    \centering
    \begin{subfigure}[t]{0.3\linewidth}
        \centering
        \includegraphics[width=\linewidth]{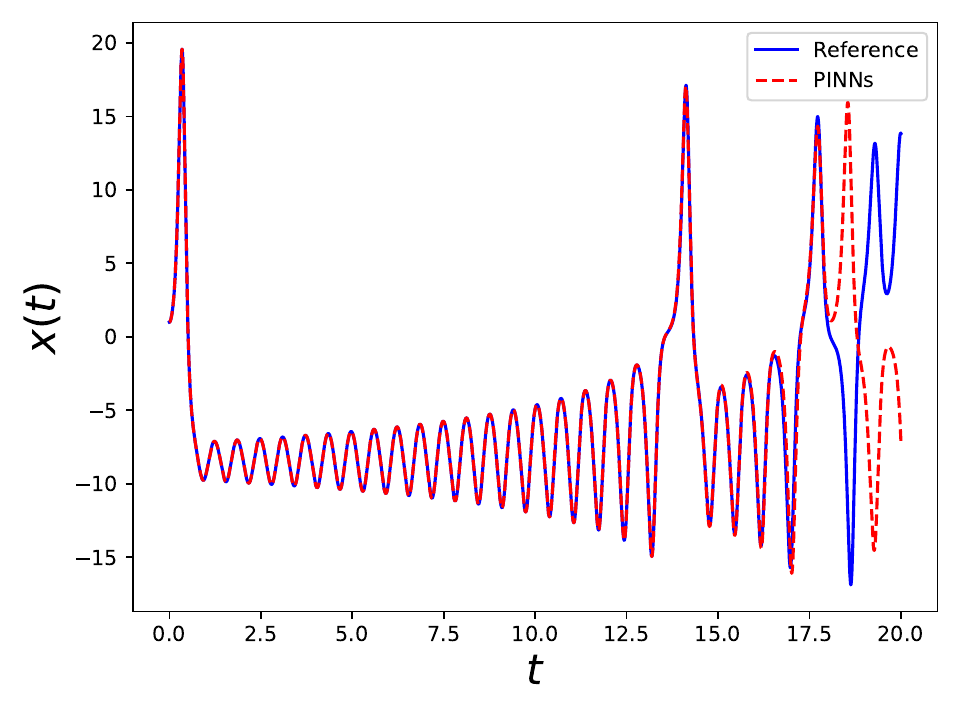}
        \includegraphics[width=\linewidth]{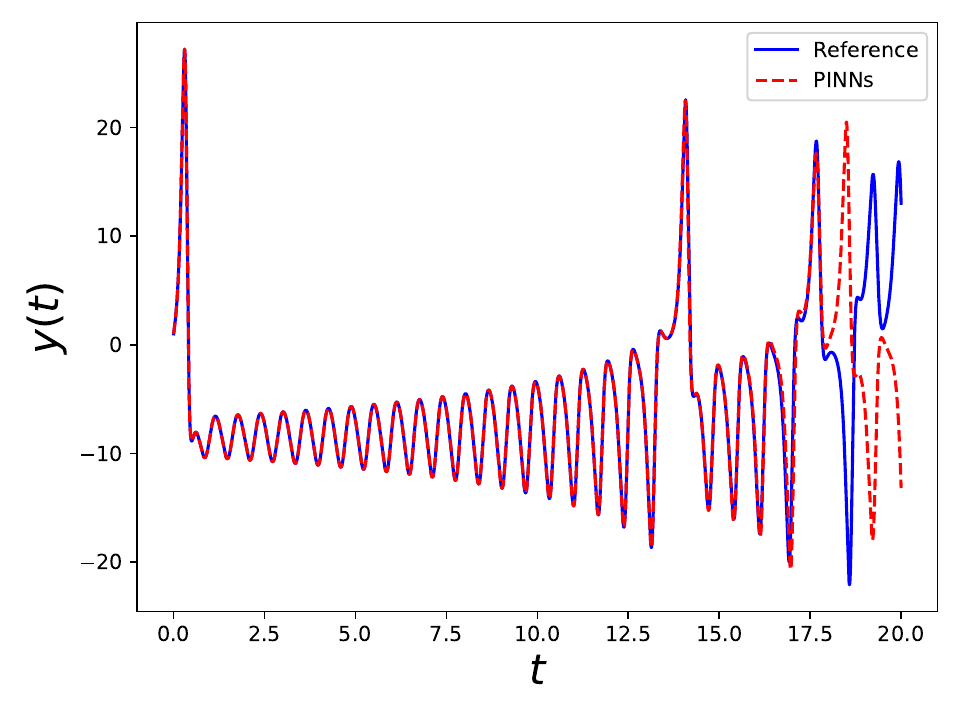}
        \includegraphics[width=\linewidth]{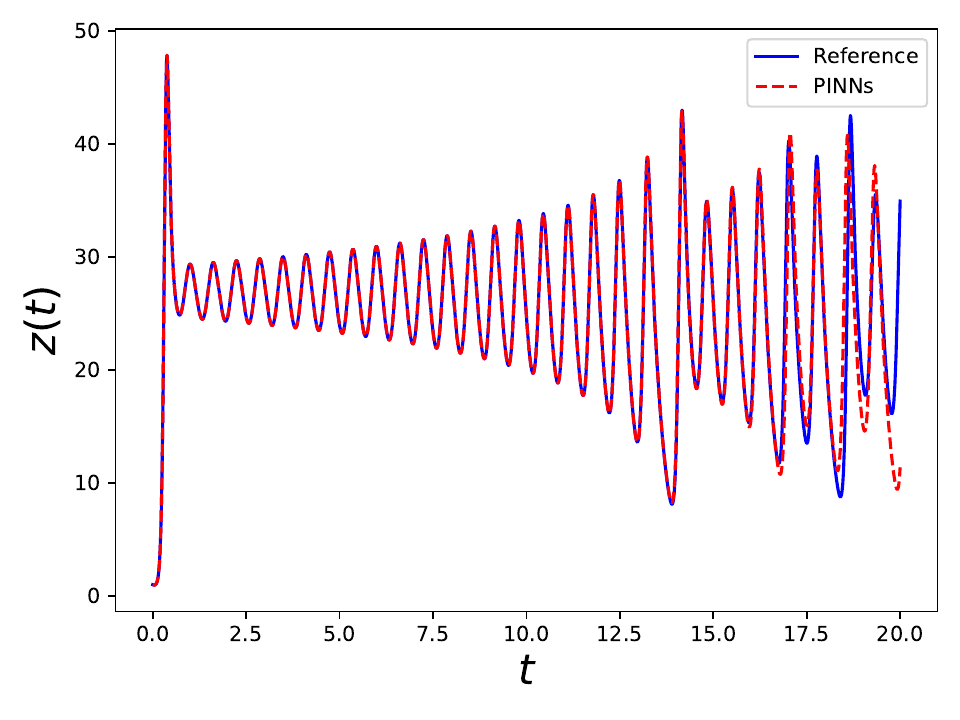}
        \caption{PINNs}
    \end{subfigure}
    \begin{subfigure}[t]{0.3\linewidth}
        \centering
        \includegraphics[width=\linewidth]{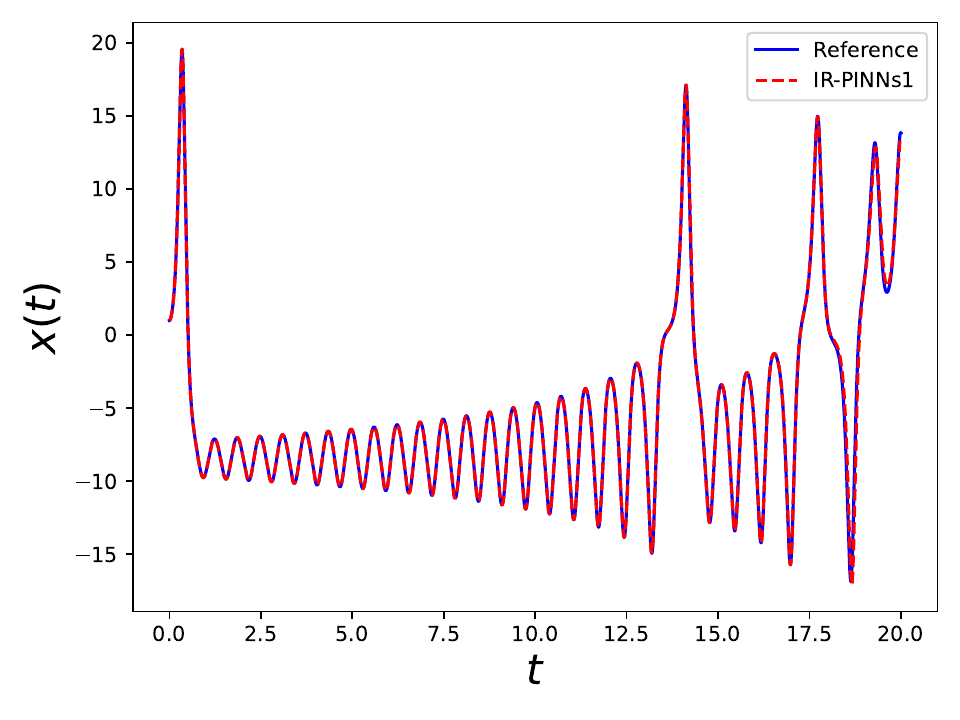}
        \includegraphics[width=\linewidth]{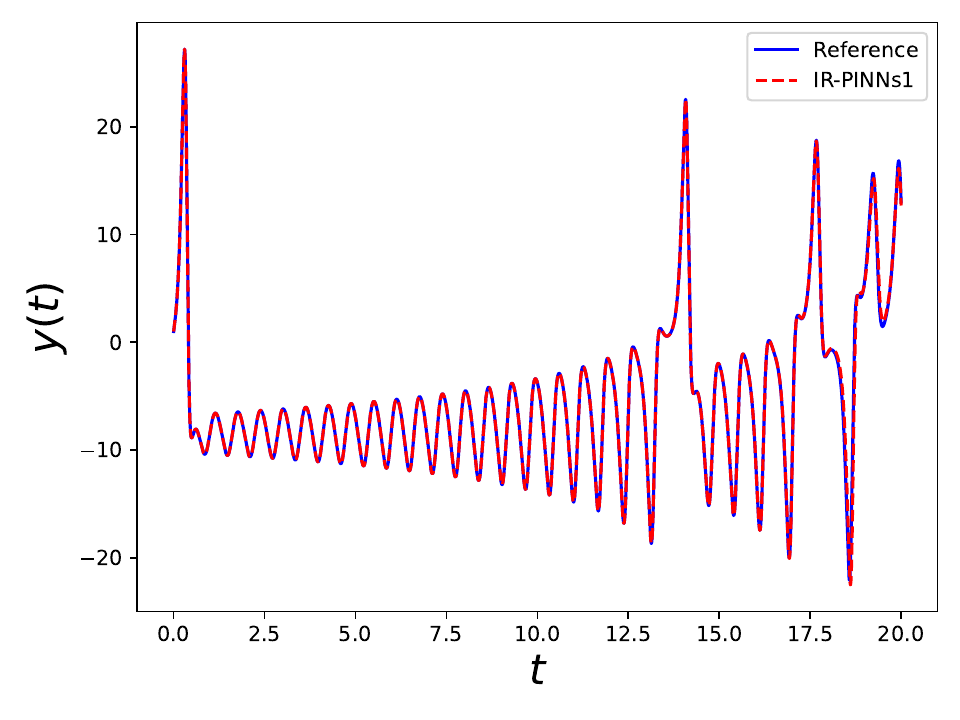}
        \includegraphics[width=\linewidth]{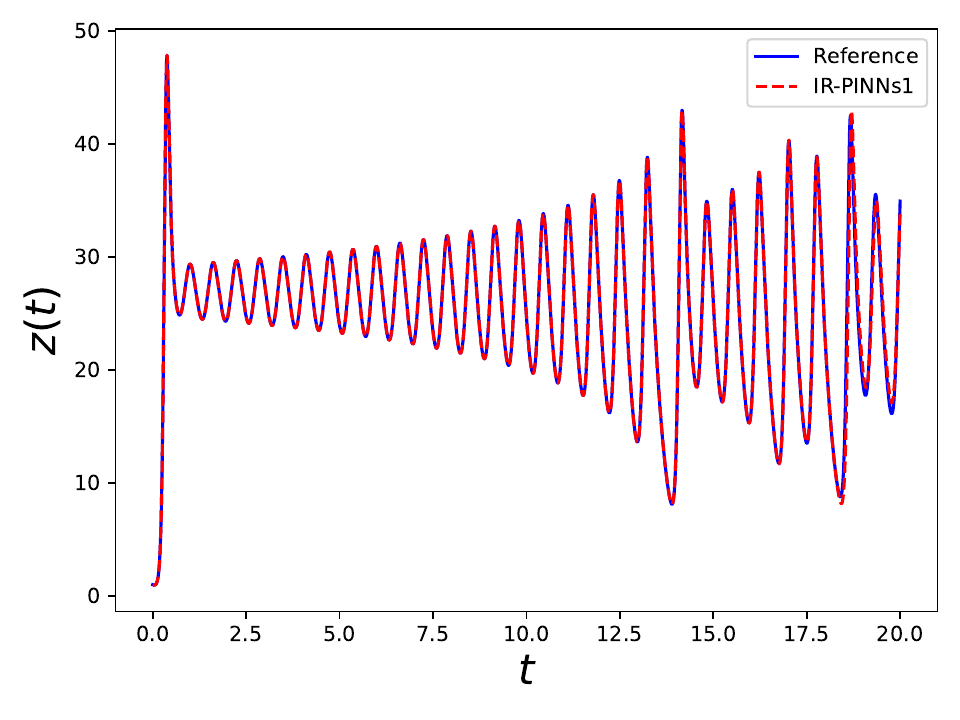}
        \caption{IR-PINNs1}
    \end{subfigure}
    \begin{subfigure}[t]{0.3\linewidth}
        \centering
        \includegraphics[width=\linewidth]{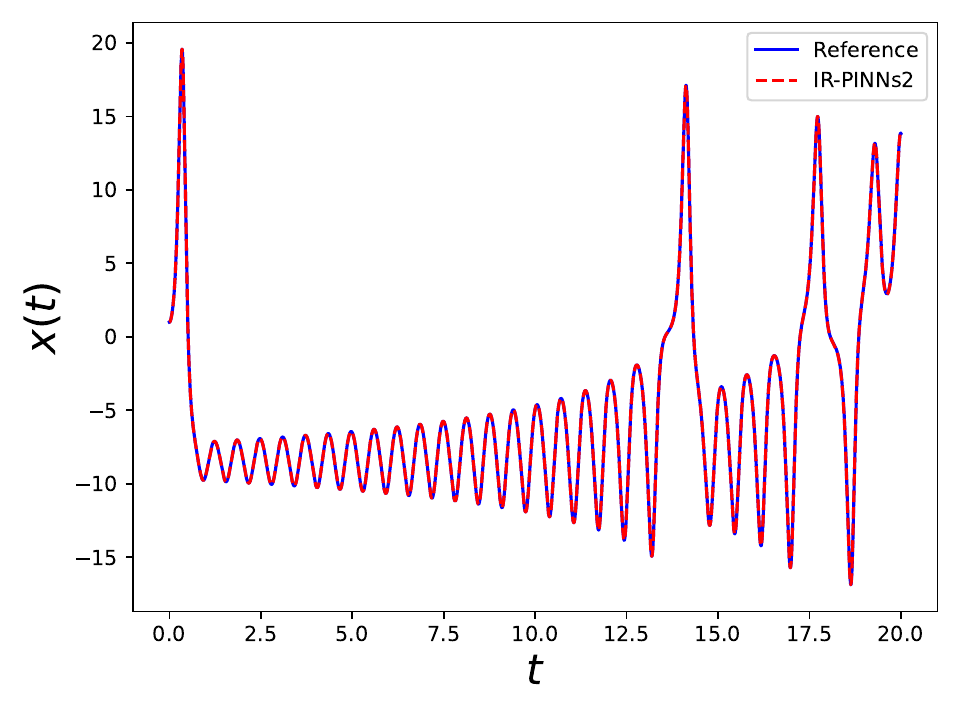}
        \includegraphics[width=\linewidth]{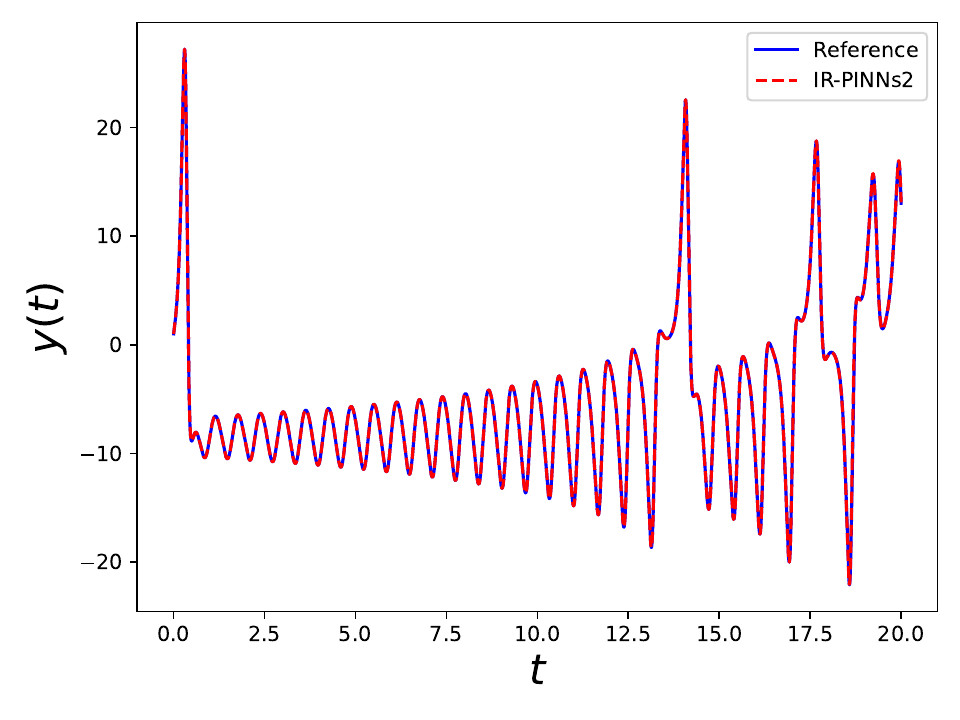}
        \includegraphics[width=\linewidth]{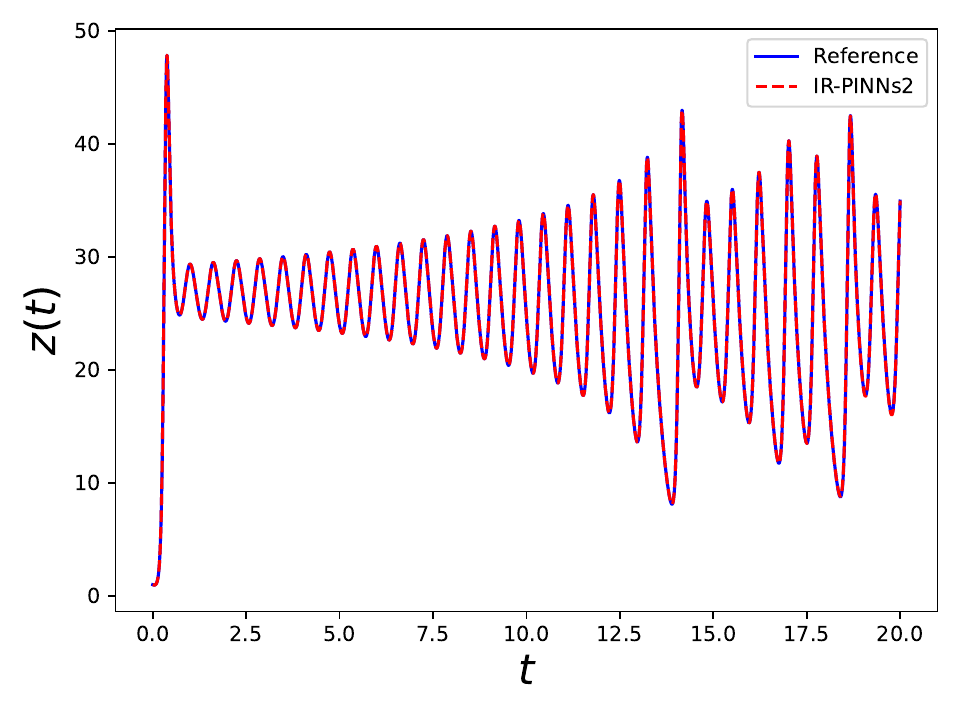}
        \caption{IR-PINNs2}
    \end{subfigure}

    \caption{\textit{Lorentz system:} Comparison between the reference and numerical solutions.}
    \label{fig:Lorentz_solution}
\end{figure}

\begin{figure}[H]
    \centering
    \includegraphics[width=0.3\linewidth]{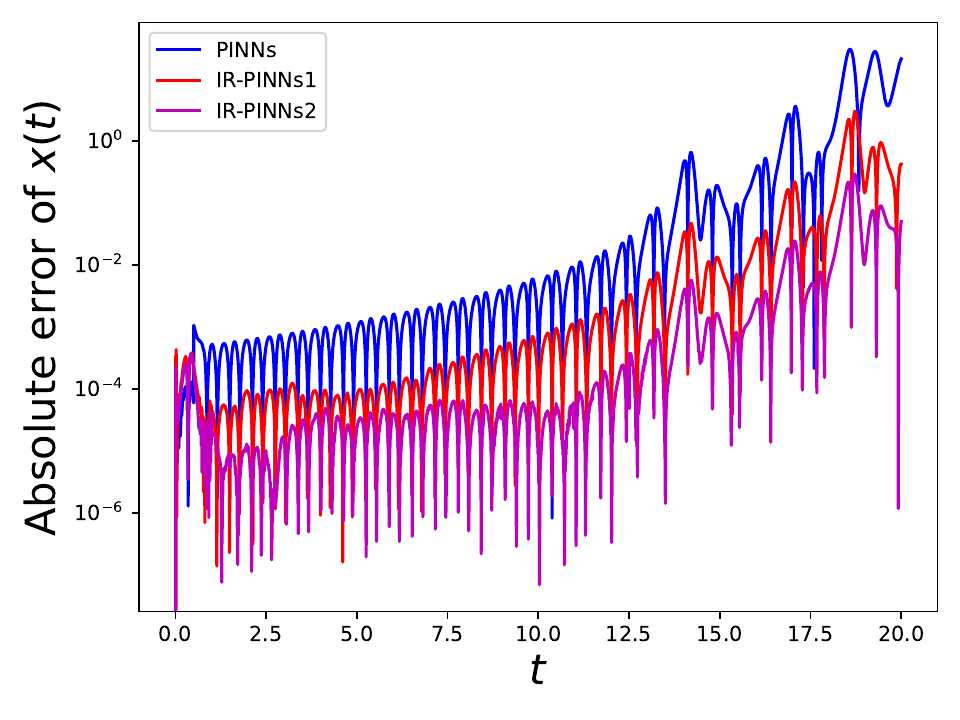}
    \includegraphics[width=0.3\linewidth]{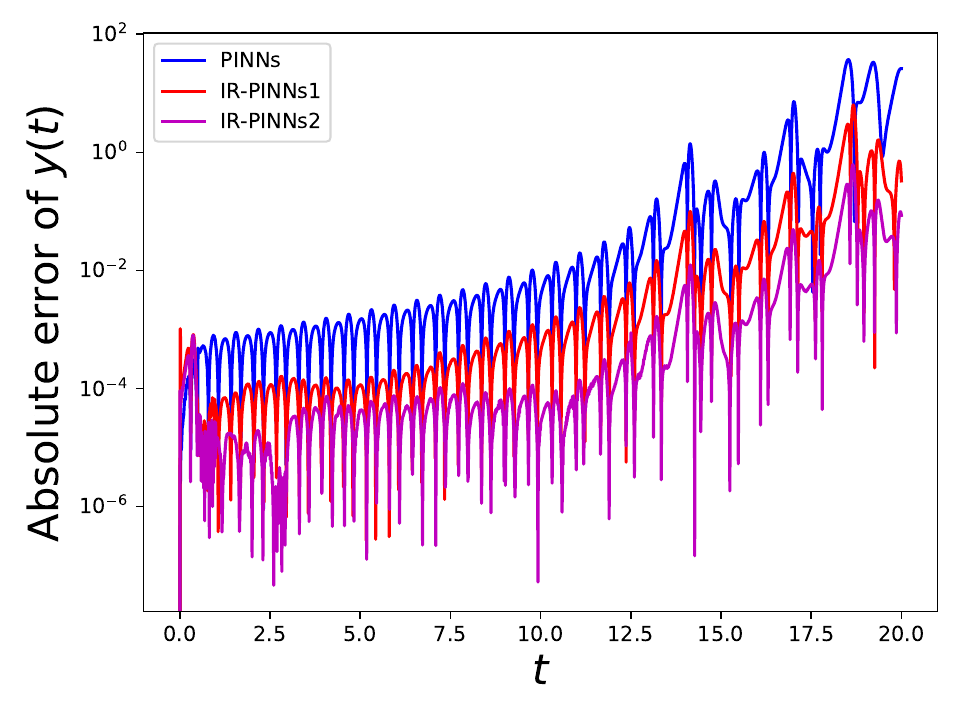}
    \includegraphics[width=0.3\linewidth]{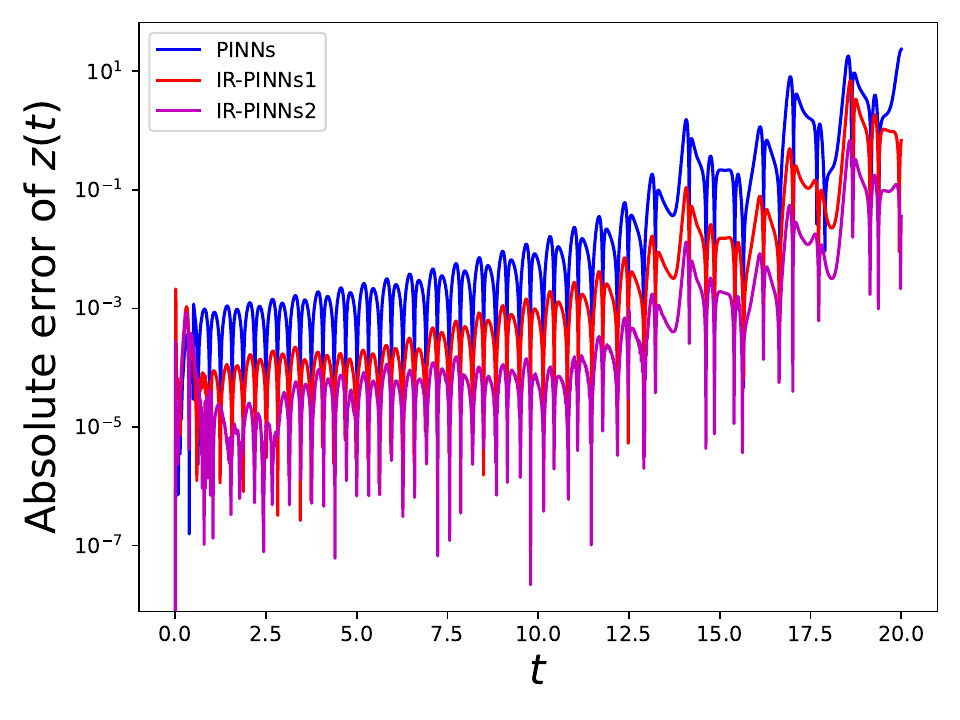}
    \caption{\textit{Lorentz system:} Absolute errors of $x(t)$, $y(t)$, and $z(t)$.}
    \label{fig:Lorentz_error}
\end{figure}

\begin{table}[H]
    \centering
    \begin{tabular}{cccc}
        \toprule
        Relative $L_2$ error & PINNs      & IR-PINNs1  & IR-PINNs2  \\
        \midrule
        $x$                  & 5.5445e-01 & 3.8330e-02 & 3.5086e-03 \\
        $y$                  & 6.0381e-01 & 5.5684e-02 & 5.1053e-03 \\
        $z$                  & 9.1354e-02 & 2.3450e-02 & 2.1651e-03 \\
        \midrule
        Running time (hours) & 1.557      & 1.961      & 1.931      \\
        \bottomrule
    \end{tabular}

    \caption{\textit{Lorentz system:} Relative $L_2$ errors and running time of different methods.}
    \label{table:lorentz}
\end{table}

For a fair comparison of our method with PINNs, we employ the same hyperparameter settings. Figure \ref{fig:Lorentz_solution} and Figure \ref{fig:Lorentz_error} present the predicted trajectories in comparison with the reference trajectories, which are obtained using \texttt{scipy.integrate.odeint} with default settings. In Table \ref{table:lorentz}, we provide the relative $L_2$ errors and running time, demonstrating that IR-PINNs1 and IR-PINNs2 achieve superior accuracy at an acceptable computational cost.

\begin{table}[H]
    \centering
    \begin{tabular}{ccc}
        \toprule
        Relative $L_2$ error & IR-PINNs2  & Causal PINNs \\
        \midrule
        $x$                  & 3.5086e-03 & 1.139e-02    \\
        $y$                  & 5.1053e-03 & 1.656e-02    \\
        $z$                  & 2.1651e-03 & 7.038e-03    \\
        \bottomrule
    \end{tabular}

    \caption{\textit{Lorentz system:} Relative $L_2$ errors between the proposed method and causal PINNs.}
    \label{table:lorentz_causal}
\end{table}

We also compare the performance and computational cost of the proposed method with causal PINNs in Table \ref{table:lorentz_causal}. In \cite{wang2024respecting}, causal PINNs generate the same number of equidistant points instead of Gaussian quadrature points in the temporal direction. We can see that the proposed method shows better performance compared to causal PINNs. Additionally, casual PINNs require a total of more than $6,000,000$ training epochs, indicating that the running time is at least greater than that of PINNs, and similar to or even more than that of the proposed method.

\subsection{Kuramoto-Sivashinsky equation}
\begin{figure}[H]
    \centering
    \includegraphics[width=0.3\linewidth]{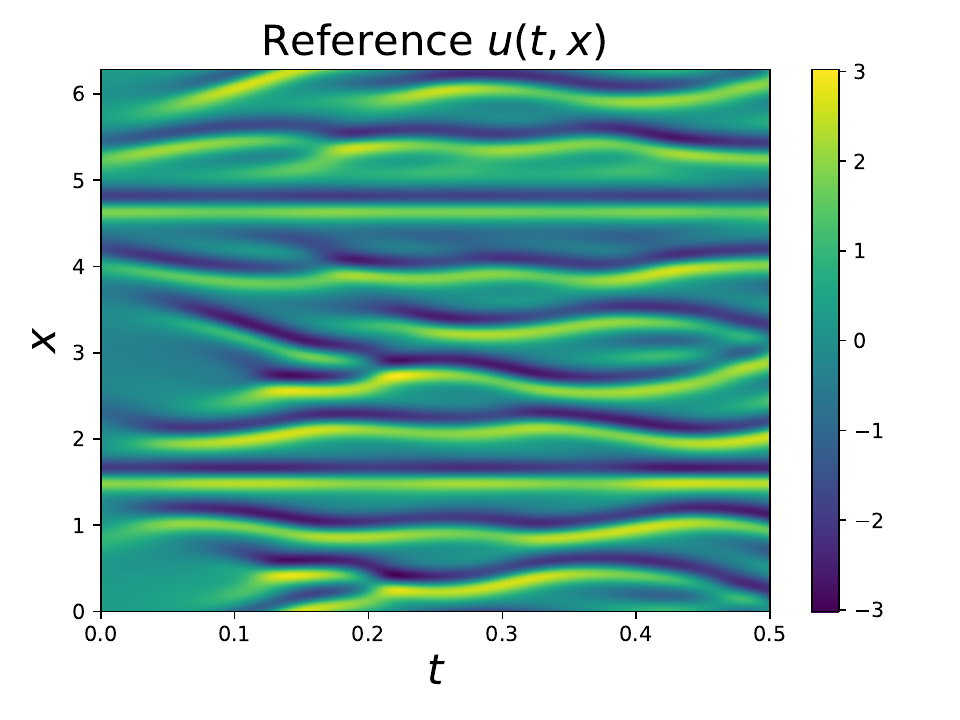}
    \caption{\textit{Kuramoto-Sivashinsky equation:} Reference solution.}
    \label{fig:KS_reference}
\end{figure}

The next example aims to demonstrate the effectiveness of the proposed method in tackling spatial-temporal chaotic systems. Consider one-dimensional Kuramoto-Sivashinsky equation, which has been independently derived in the context of reaction-diffusion systems and flame front propagation. The Kuramoto-Sivashinsky equation is a classic model describing spatial-temporal chaotic dynamical behavior. The equation is expressed as follows:
\begin{equation}
    u_{t} + \alpha uu_x + \beta u_{xx}+ \gamma u_{xxxx} = 0,
\end{equation}
subject to periodic boundary conditions and an initial condition
\begin{equation}
    u(0, x) = u_0(x).
\end{equation}
Here we set  $\alpha = 100/16$, $\beta=100/16^2, \gamma=100/16^4$, for a fixed spatial domain in $[0, 2 \pi]$. For validation, we use the Chebfun package \cite{driscoll2014chebfun} with a spectral Fourier discretization with $512$ modes and a fourth-order stiff time-stepping scheme (ETDRK4) with time-step size $10^{-5}$. Starting from the initial condition $u_0(x) = \cos(x)(1+\sin(x))$, we select the numerical solution at $t = 0.4$ as our initial condition in the chaotic regime. Our goal is to learn the associated solution from initial time $t_0 = 0.4$ to terminal time $T = 0.9$. The reference solution is shown in Figure \ref{fig:KS_reference}.

In this example, we employ the time-marching strategy and partition the temporal domain $[0, 0.5]$ into $5$ subdomains. For each subdomain, we construct a $5$-layer fully-connected neural network with $512$ neurons per hidden layer. The initial learning rate is $10^{-3}$, with an exponential decay rate of $0.9$ applied every $5,000$ training epochs. The maximum number of training epochs is $600,000$. We divide the subdomain into $N = 8$ subintervals and in each subinterval, the number of Gaussian quadrature points is $M = 4$. Spatially we choose $256$ equidistant points for enforcing the PDE residual and $512$ equidistant initial points. To further simplify the training objective, we also strictly impose the periodic boundary conditions by embedding the input coordinates into Fourier expansion (see Appendix \ref{sec:periodic_BC}). The initial loss weight coefficient is set to $\lambda_2 = 10,000$.

\begin{figure}[H]
    \begin{subfigure}[t]{0.3\linewidth}
        \centering
        \includegraphics[width=\linewidth]{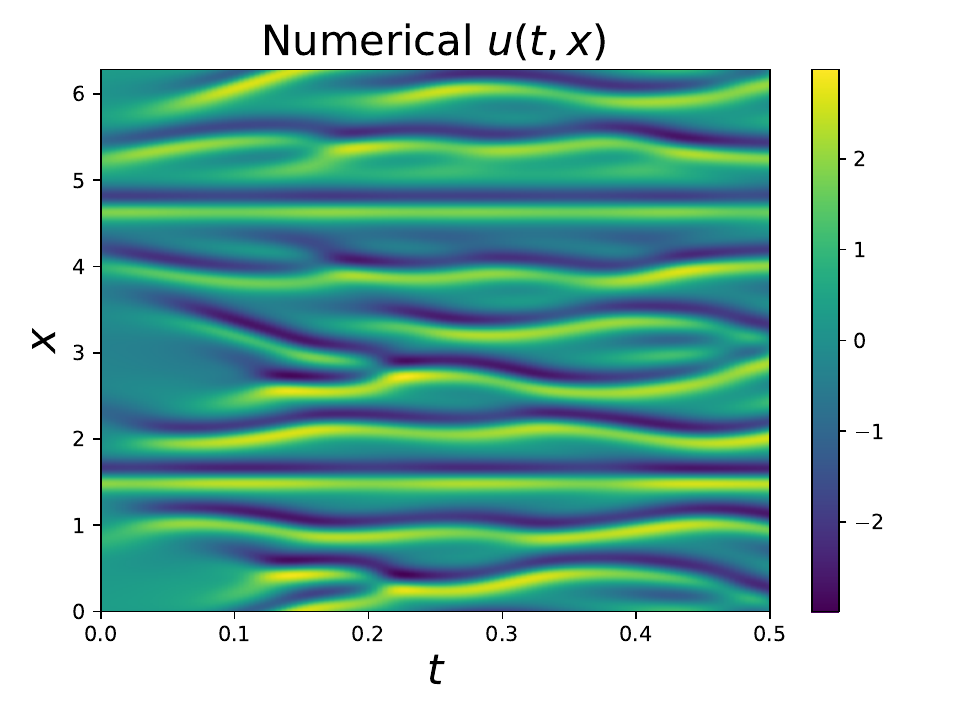}
        \includegraphics[width=\linewidth]{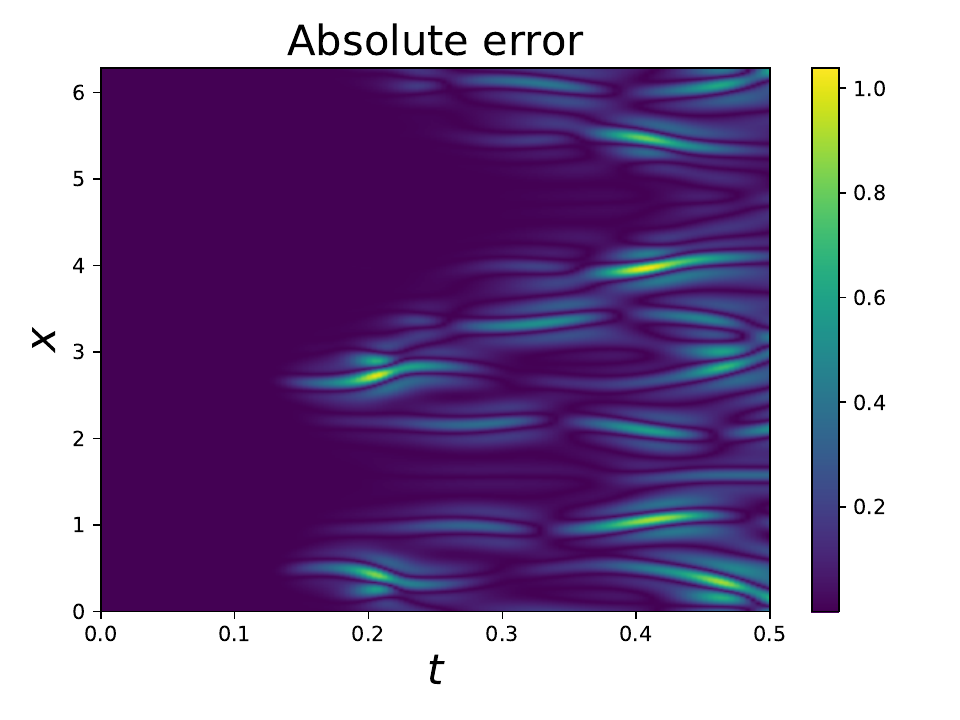}
        \caption{PINNs}
    \end{subfigure}
    \begin{subfigure}[t]{0.3\linewidth}
        \centering
        \includegraphics[width=\linewidth]{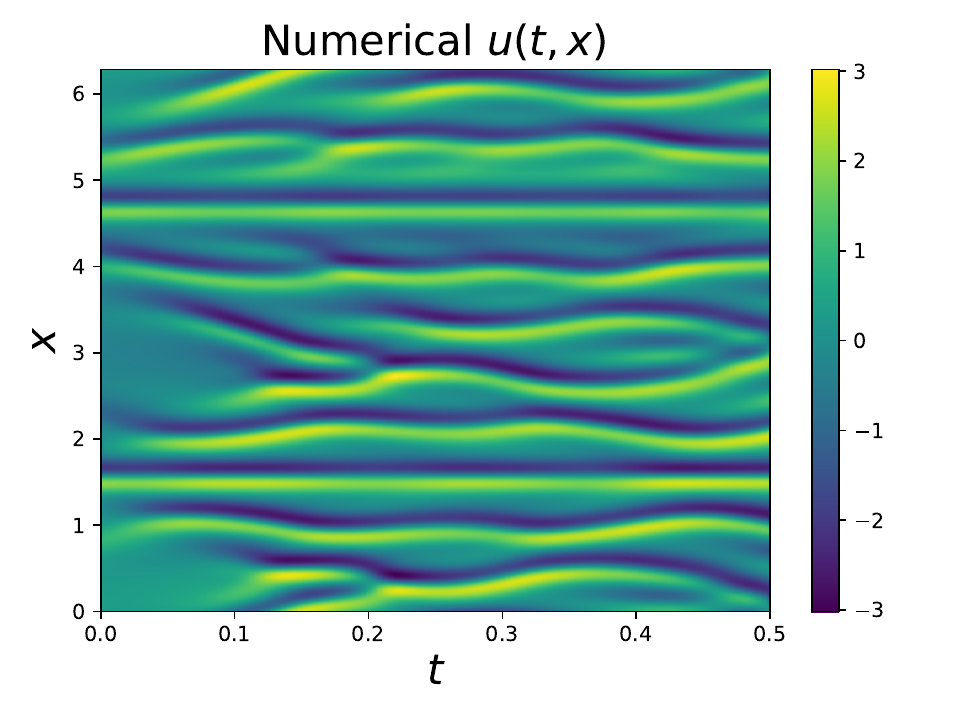}
        \includegraphics[width=\linewidth]{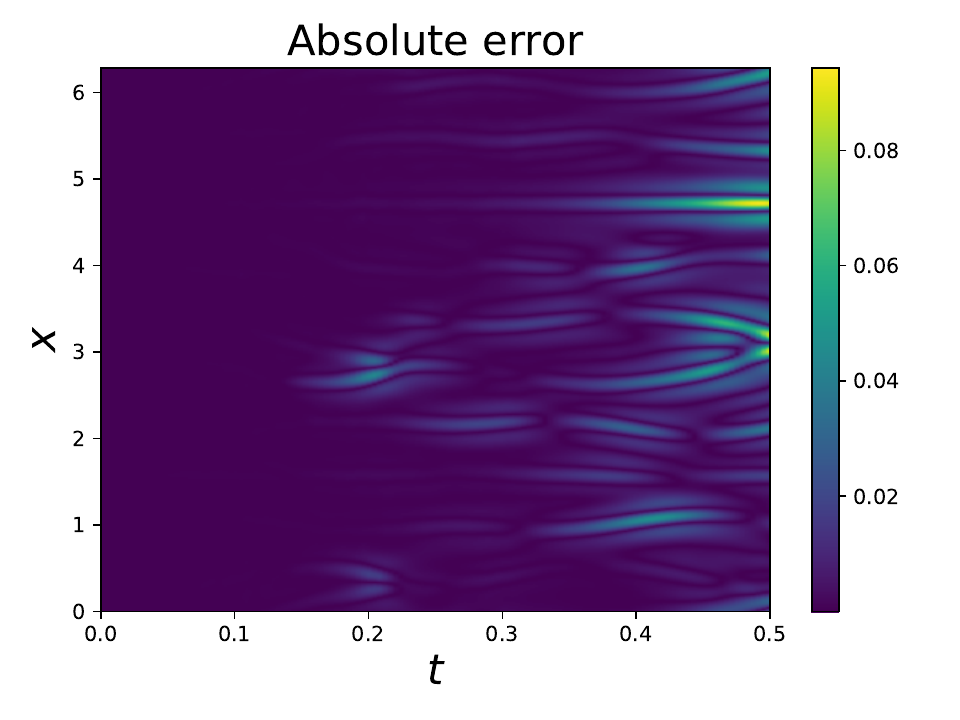}
        \caption{IR-PINNs1}
    \end{subfigure}
    \begin{subfigure}[t]{0.3\linewidth}
        \centering
        \includegraphics[width=\linewidth]{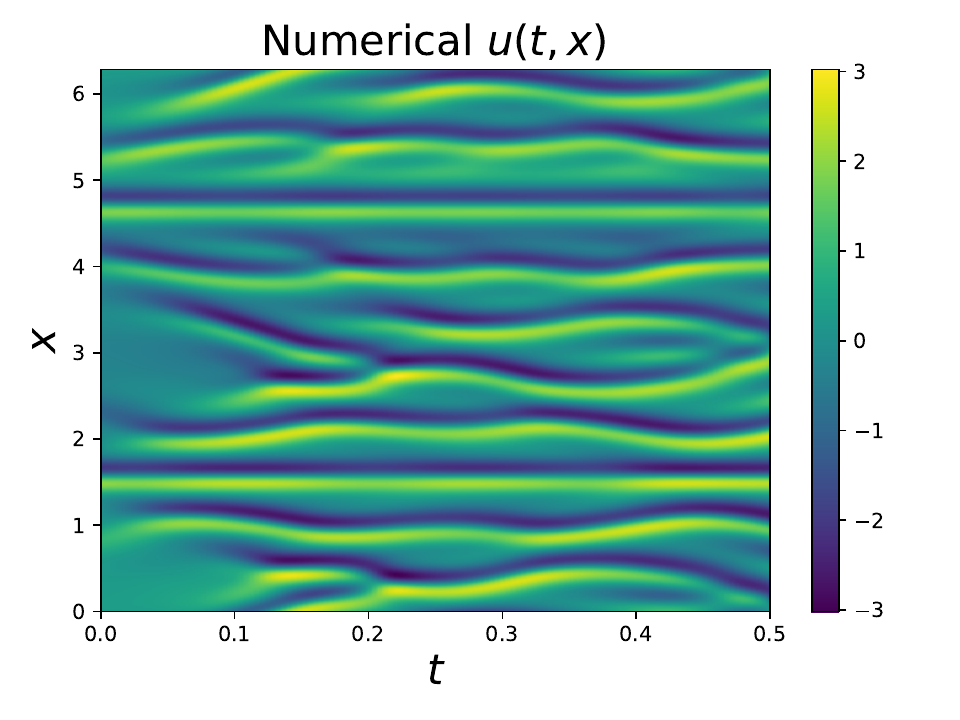}
        \includegraphics[width=\linewidth]{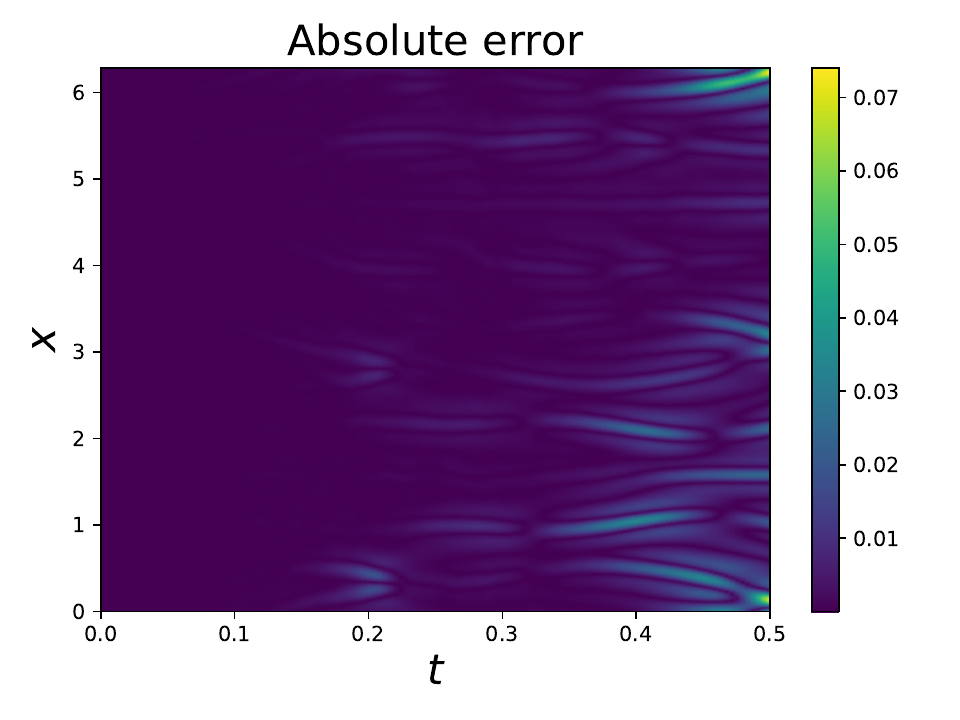}
        \caption{IR-PINNs2}
    \end{subfigure}
    \caption{\textit{Kuramoto-Sivashinsky equation:} Numerical solutions and absolute errors.}
    \label{fig:KS_solution}
\end{figure}

\begin{figure}[H]
    \centering
    \includegraphics[width=0.33\linewidth]{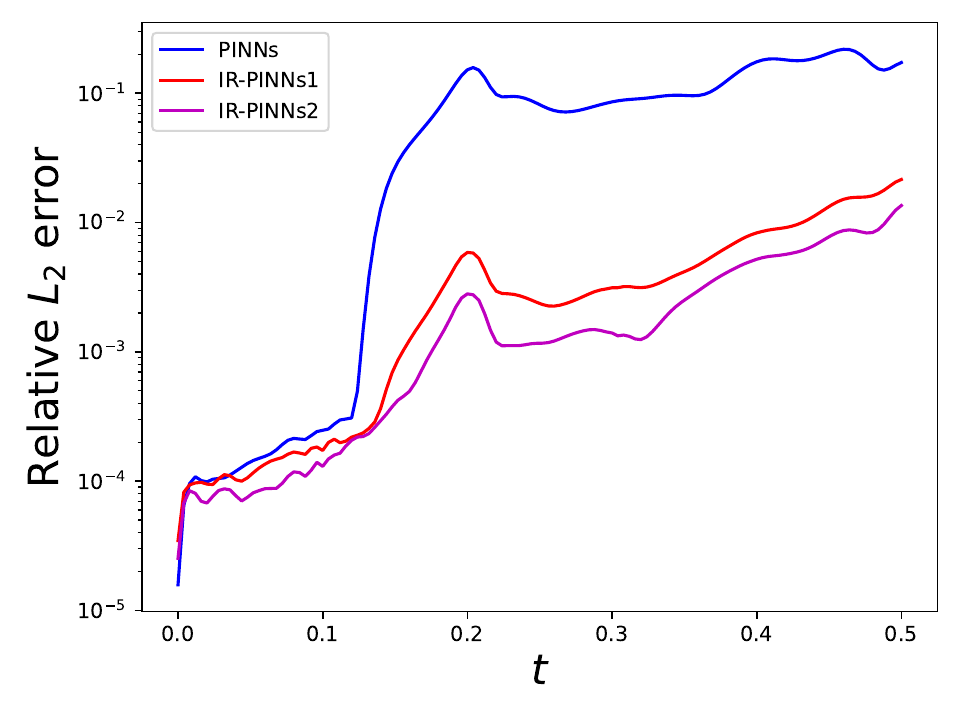}
    \caption{\textit{Kuramoto-Sivashinsky equation:} Relative $L_2$ errors over time.}
    \label{fig:KS_error}
\end{figure}

\begin{table}[H]
    \centering
    \begin{tabular}{cccc}
        \toprule
        Relative $L_2$ error & PINNs      & IR-PINNs1  & IR-PINNs2  \\
        \midrule
        $u$                  & 1.1278e-01 & 6.9356e-03 & 3.9945e-03 \\
        \midrule
        Running time (hours) & 46.91      & 49.36      & 48.23      \\
        \bottomrule
    \end{tabular}

    \caption{\textit{Kuramoto-Sivashinsky equation:} Relative $L_2$ errors and running time of different methods.}
    \label{table:KS}
\end{table}

The numerical results are summarized in Figure \ref{fig:KS_solution}, from which one can see that the numerical solutions of IR-PINNs1 and IR-PINNs2 exhibit good agreements with the reference solution. In Figure \ref{fig:KS_error}, We can clearly observe that for PINNs, the relative $L_2$ error experiences a sudden and significant increase within the temporal subdomain $[0.1, 0.2]$, resulting in a substantial deviation between the initial values and the true values in subsequent temporal subdomains, ultimately leading to the failure of the simulation. In contrast, both IR-PINNs1 and IR-PINNs2 effectively control the rise of the error. Table \ref{table:KS} provides the relative $L_2$ errors and running time, demonstrating that our proposed method is equally effective in describing the PDEs of complex chaotic systems.

\begin{table}[H]
    \centering
    \begin{tabular}{ccc}
        \toprule
        Relative $L_2$ error & IR-PINNs2  & Causal PINNs \\
        \midrule
        $u$                  & 3.9945e-03 & 2.2600e-02   \\
        \bottomrule
    \end{tabular}

    \caption{\textit{Kuramoto-Sivashinsky equation:} Relative $L_2$ errors between the proposed method and causal PINNs.}
    \label{table:KS_causal}
\end{table}

We also compare the performance and computational cost of the proposed method with causal PINNs in Table \ref{table:KS_causal}. Similar to the previous example, the proposed method demonstrates a more significant improvement in accuracy compared to causal PINNs. Additionally, casual PINNs require a total of more than $3,000,000$ training epochs, indicating that the running time is at least greater than that of PINNs, and similar to or even more than that of the proposed method.

\subsection{Boussinesq-Burgers equations}
In order to demonstrate the effectiveness of the adaptive sampling strategy, we consider the Boussinesq-Burgers equations, a coupled system of nonlinear partial differential equations that model the interaction of nonlinearity and dispersion in shallow water wave dynamics, with Dirichlet boundary conditions,
\begin{equation}
    \begin{aligned}
         & u_t + 2uu_x - \frac{1}{2}v_x = 0,                                              \\
         & v_t + 2(uv)_x - \frac{1}{2}u_{xxx} = 0, \quad t\in[t_0, t_1], x\in [x_0, x_1].
    \end{aligned}
\end{equation}

\begin{figure}[H]
    \centering
    \includegraphics[width=0.3\linewidth]{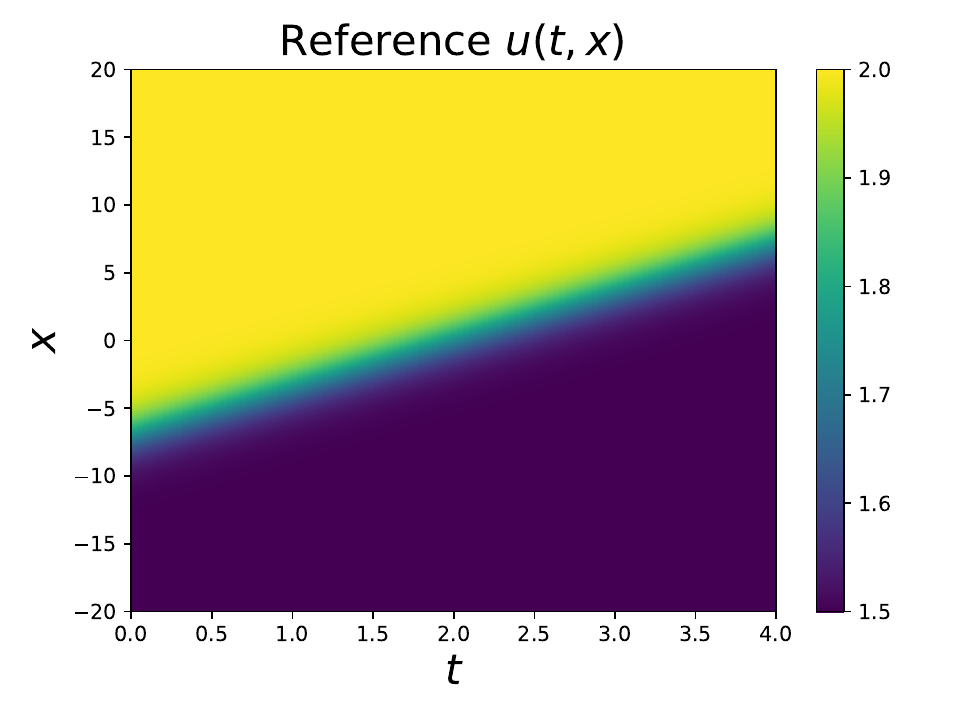}
    \includegraphics[width=0.3\linewidth]{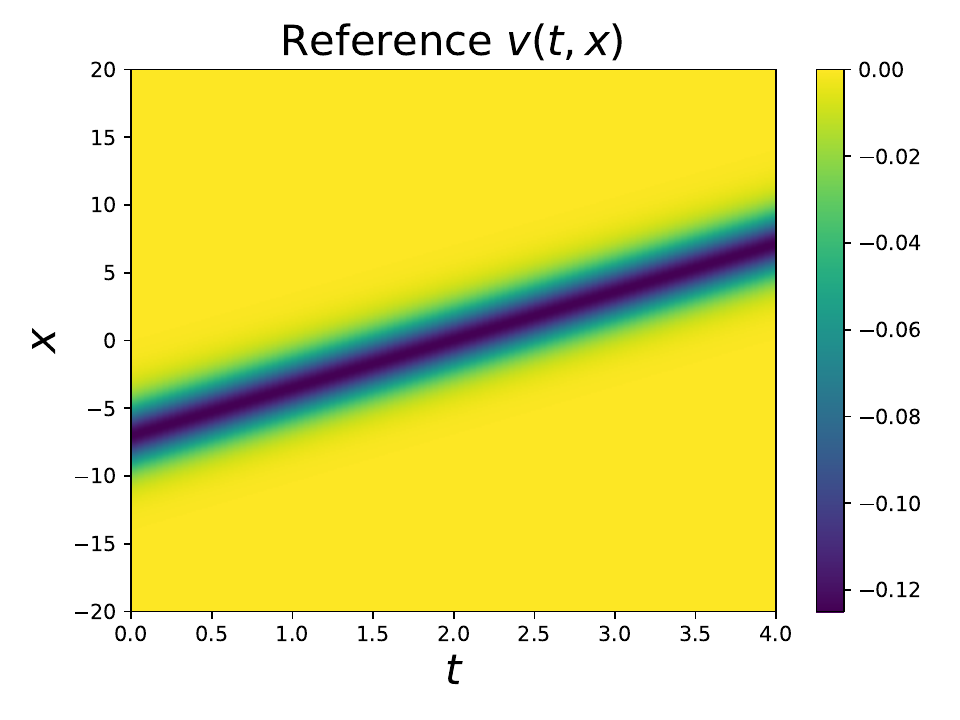}
    \caption{\textit{Boussinesq-Burgers equations:} Reference solutions of $u(t, x)$ and $v(t, x)$.}
    \label{fig:BBE_reference_solution}
\end{figure}

The one-soliton solution of the Boussinesq-Burgers equations is as follows
\begin{equation}
    \begin{aligned}
         & u(t, x) = -\frac{\exp(-x+\frac{7}{2}t-7)}{2(1+\exp(-x+\frac{7}{2}t-7))} + 2, \\
         & v(t, x) = -\frac{\exp(-x+\frac{7}{2}t-7)}{2(1+\exp(-x+\frac{7}{2}t-7))^2},
    \end{aligned}
\end{equation}
where the spatial-temporal domain is $[t_0, t_1] \times [x_0, x_1] = [0, 4] \times [-20, 20]$, as shown in Figure \ref{fig:BBE_reference_solution}. Initial and boundary conditions are obtained from the exact solution.

For the training process, we construct two $4$-layer fully-connected neural networks with $64$ neurons per hidden layer. The initial learning rate is $10^{-3}$, with an exponential decay rate of $0.9$ applied every $500$ training epochs. The maximum number of training epochs is $20,000$. We divide the subdomain into $N = 8$ subintervals and in each subinterval, the number of Gaussian quadrature points is $M = 4$. Before the first training, we select $8\times4\times16$ collocation points for enforcing the PDE residual, $256$ equidistant initial points and $256$ equidistant boundary points. The boundary and initial loss weight coefficients are both set to $\lambda_2 = \lambda_3 = 10$.

The number of adaptive iterations is $2$. For each adaptive iteration, we take $20$ bins for bounded polynomial layer and $8$ CDF coupling layers for bounded KR-net, with each layer parameterized by a $2$-layer fully-connected neural network with $64$ neurons. We randomly generate $8\times4\times100$ samples for evaluating PDE residual as our initial training set and the maximum number of adaptive training epochs for training $p(t, \bm{x}; \theta_f)$ is $3,000$. To refine the training set, we generate $8\times4\times8$ new samples from the density model $p(t, \bm{x}; \theta_f)$.

\begin{figure}[H]
    \centering
    \begin{subfigure}[t]{0.3\linewidth}
        \centering
        \includegraphics[width=\linewidth]{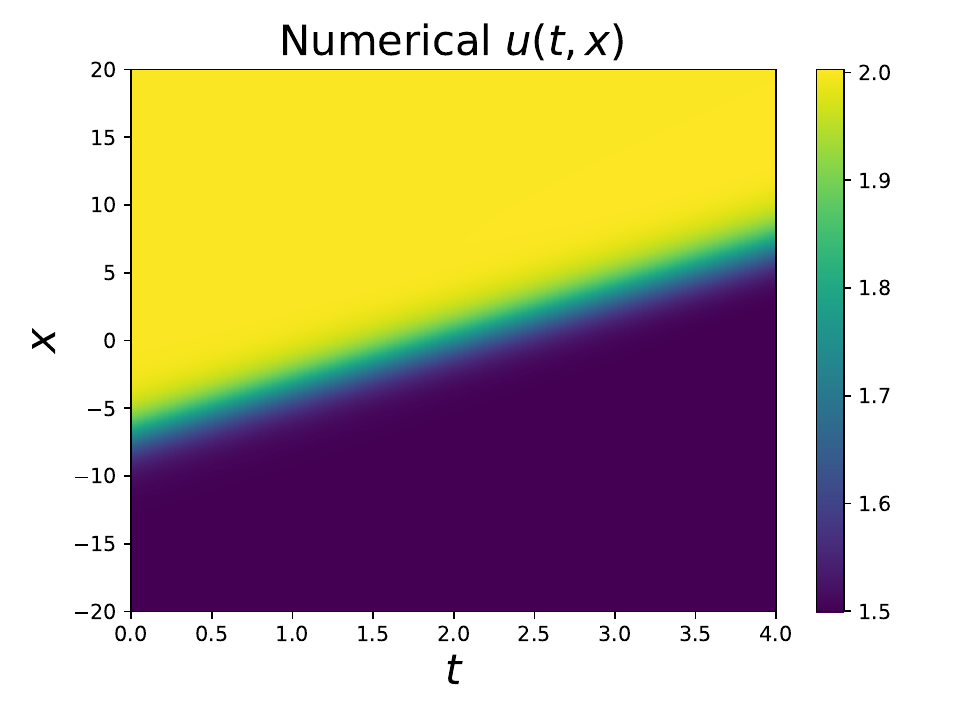}
        \includegraphics[width=\linewidth]{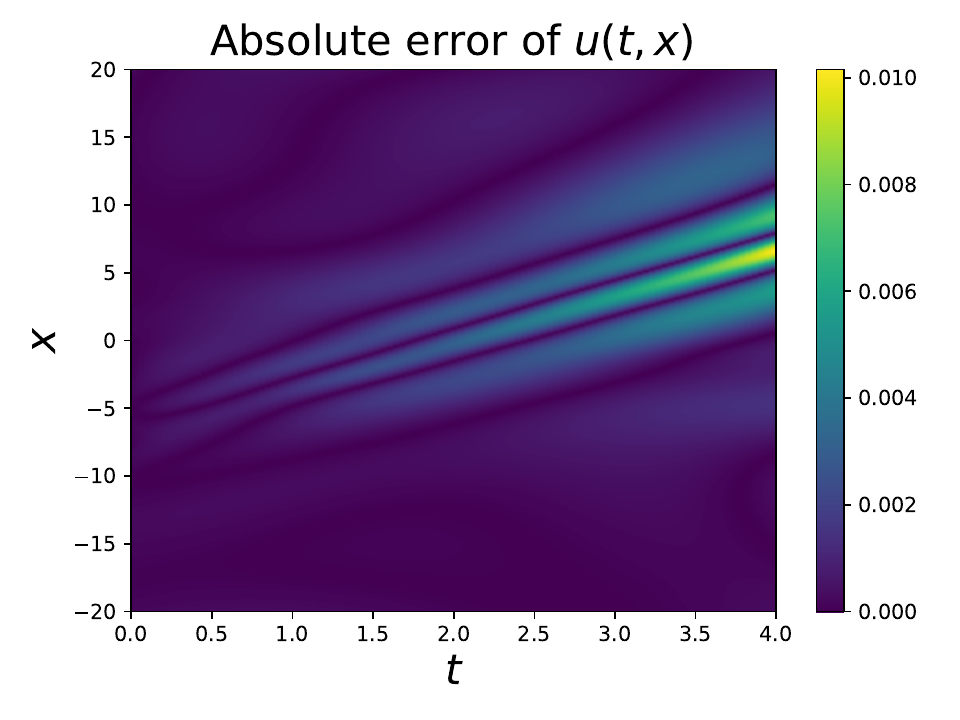}
        \caption{PINNs}
    \end{subfigure}
    \begin{subfigure}[t]{0.3\linewidth}
        \centering
        \includegraphics[width=\linewidth]{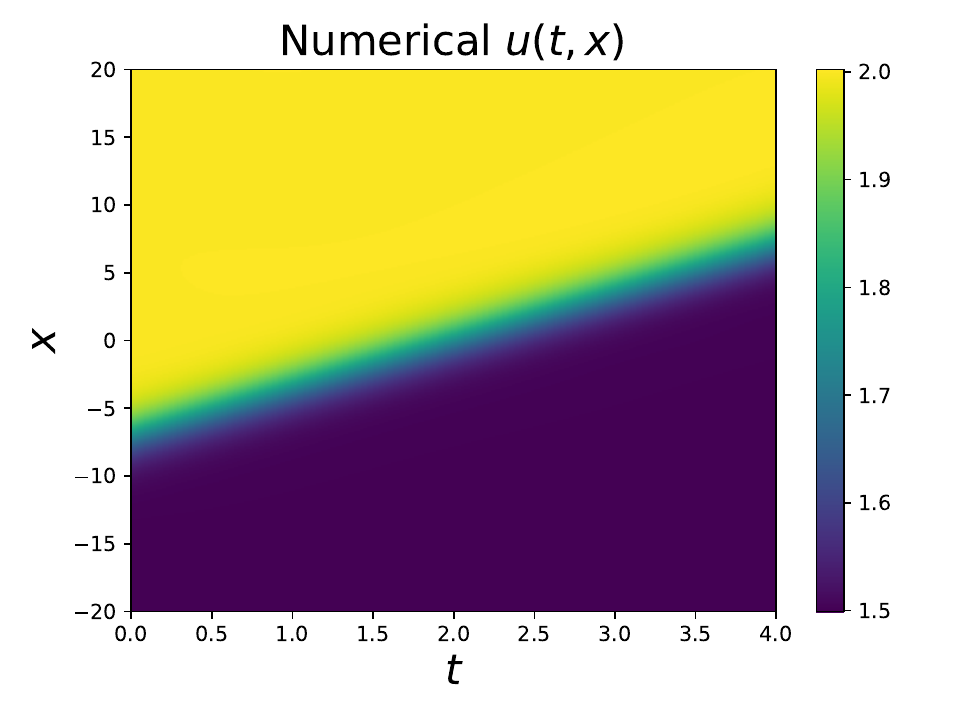}
        \includegraphics[width=\linewidth]{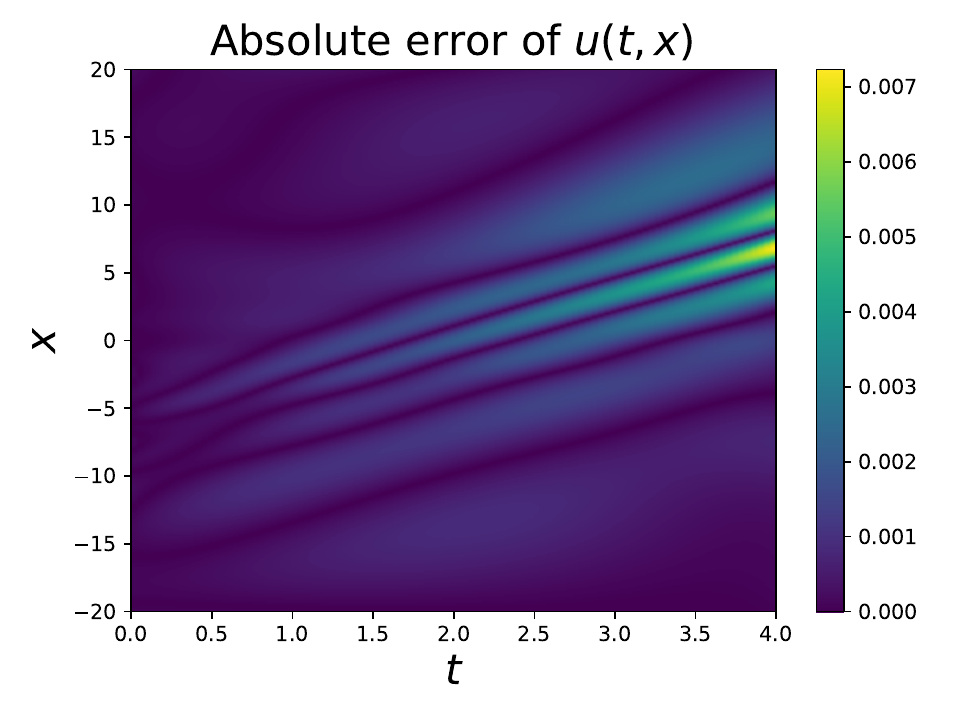}
        \caption{IR-PINNs1}
    \end{subfigure}
    \begin{subfigure}[t]{0.3\linewidth}
        \centering
        \includegraphics[width=\linewidth]{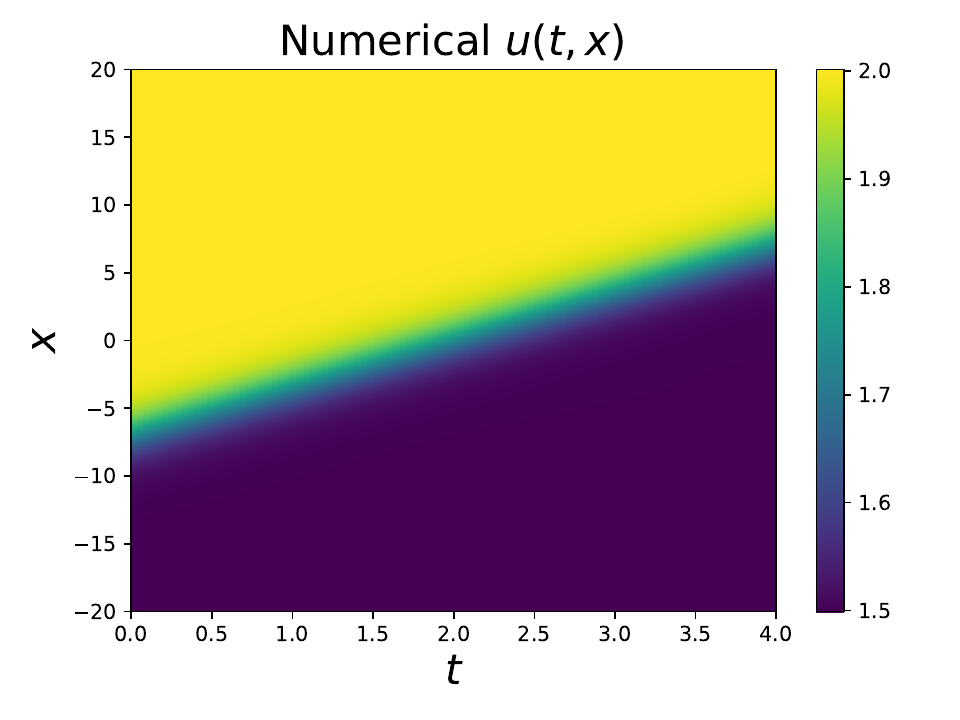}
        \includegraphics[width=\linewidth]{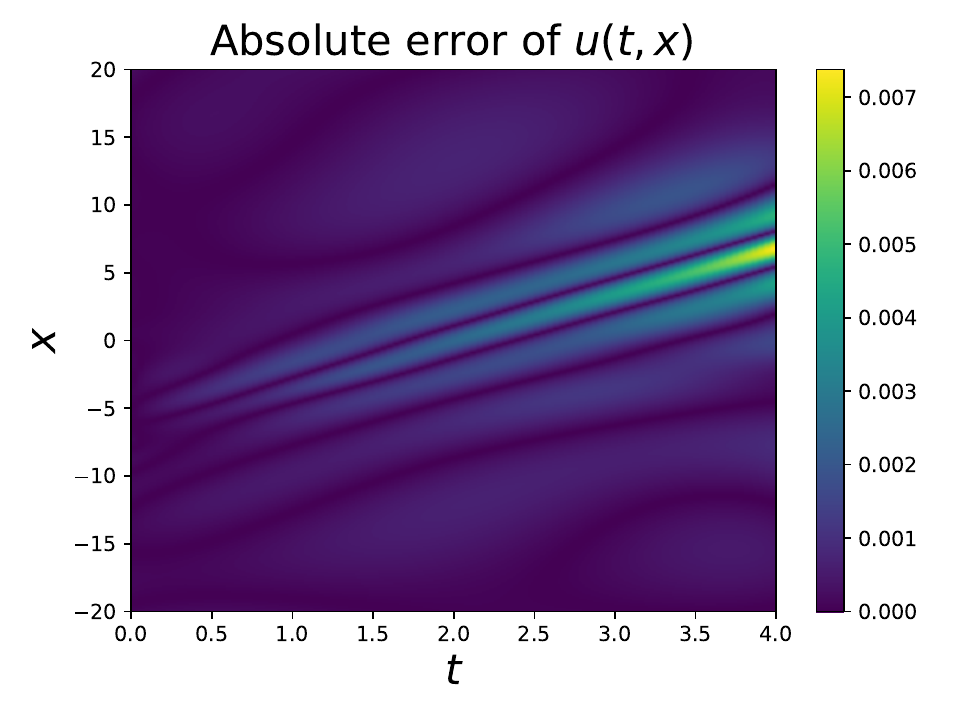}
        \caption{IR-PINNs2}
    \end{subfigure}
    \caption{\textit{Boussinesq-Burgers equations:} Numerical solutions and absolute errors of $u(t, x)$.}
    \label{fig:BBE_solution_u}
\end{figure}

\begin{figure}[H]
    \centering
    \begin{subfigure}[t]{0.3\linewidth}
        \centering
        \includegraphics[width=\linewidth]{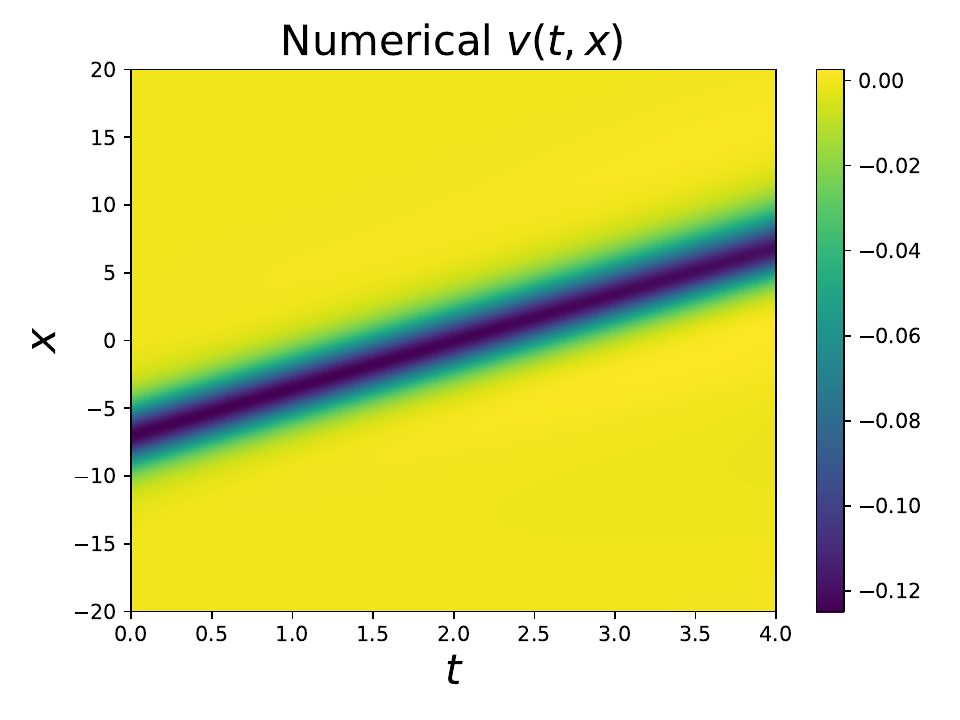}
        \includegraphics[width=\linewidth]{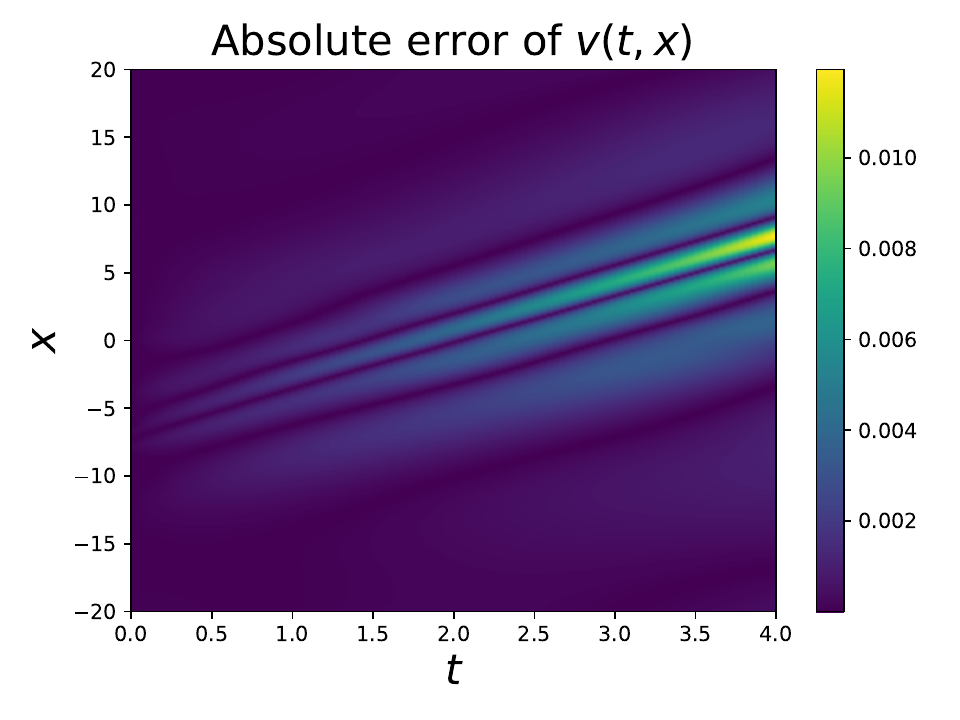}
        \caption{PINNs}
    \end{subfigure}
    \begin{subfigure}[t]{0.3\linewidth}
        \centering
        \includegraphics[width=\linewidth]{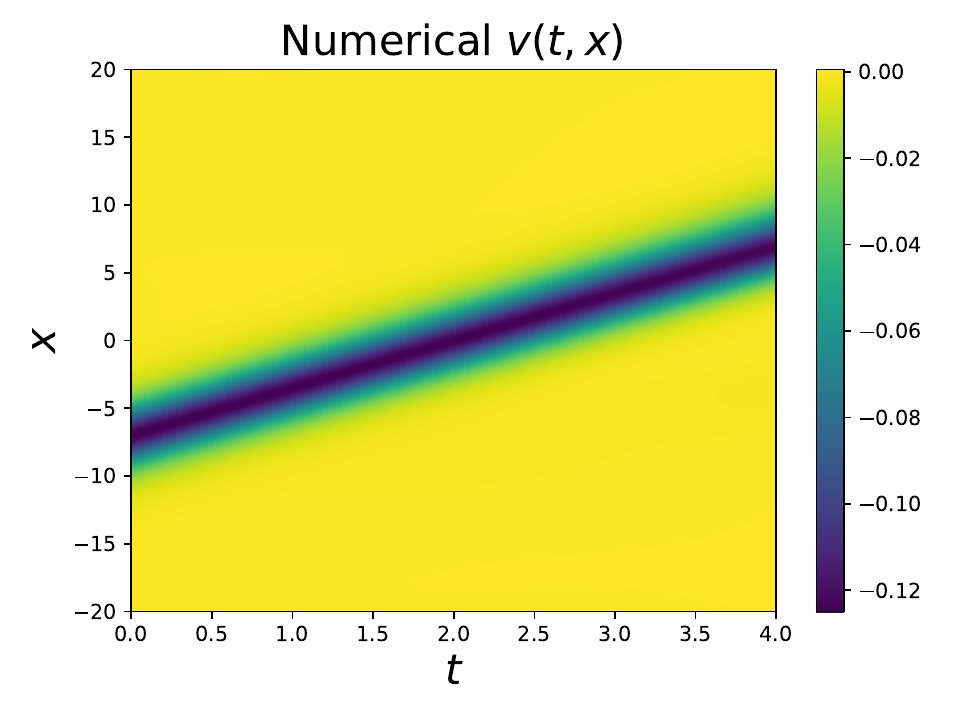}
        \includegraphics[width=\linewidth]{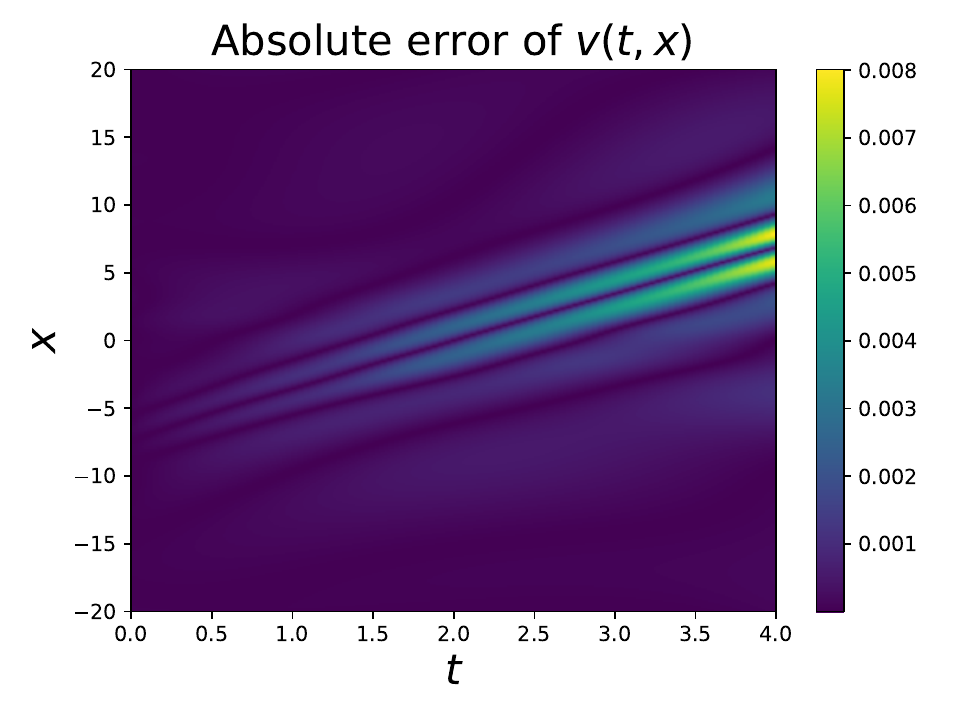}
        \caption{IR-PINNs1}
    \end{subfigure}
    \begin{subfigure}[t]{0.3\linewidth}
        \centering
        \includegraphics[width=\linewidth]{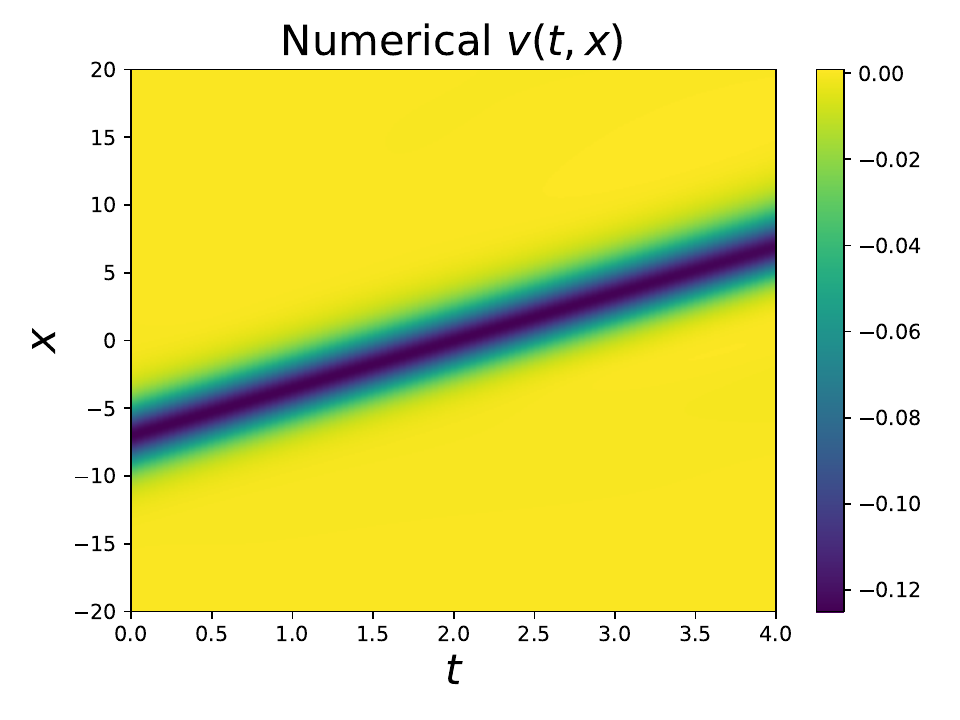}
        \includegraphics[width=\linewidth]{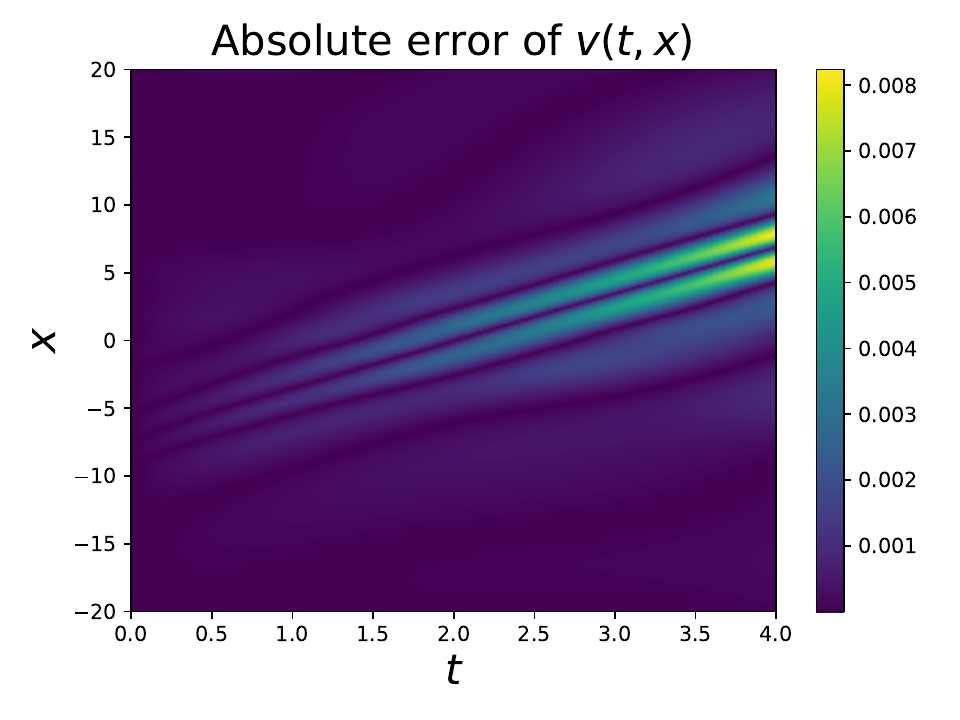}
        \caption{IR-PINNs2}
    \end{subfigure}
    \caption{\textit{Boussinesq-Burgers equations:} Numerical solutions and absolute errors of $v(t, x)$.}
    \label{fig:BBE_solution_v}
\end{figure}

\begin{figure}[H]
    \centering
    \includegraphics[width=0.3\linewidth]{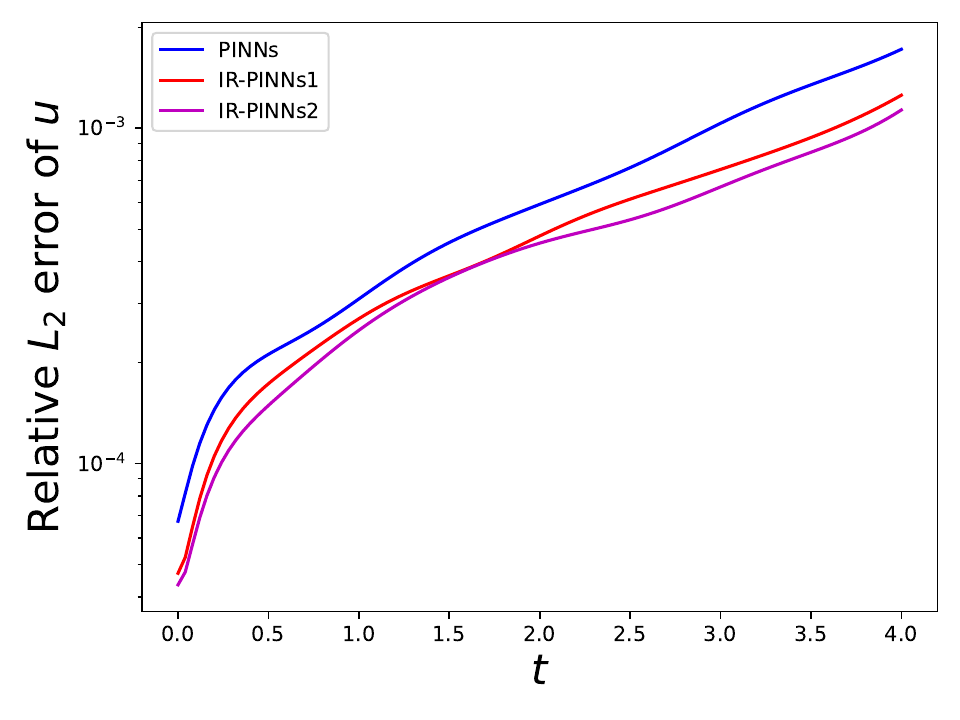}
    \includegraphics[width=0.3\linewidth]{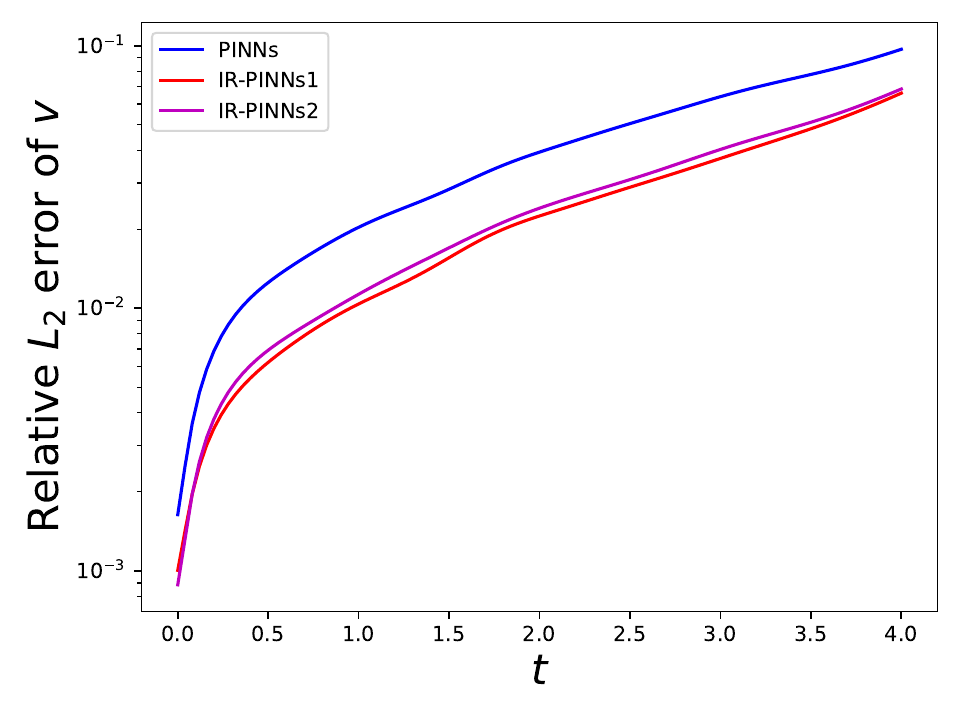}

    \caption{\textit{Boussinesq-Burgers equations:} Relative $L_2$ errors over time.}
    \label{fig:BBE_error}
\end{figure}

\begin{table}[H]
    \centering
    \begin{tabular}{cccc}
        \toprule
        Relative $L_2$ error & PINNs      & IR-PINNs1  & IR-PINNs2  \\
        \midrule
        $u$                  & 8.1022e-04 & 5.9748e-04 & 5.3934e-04 \\
        \midrule
        $v$                  & 5.0101e-02 & 3.0561e-02 & 3.2439e-02 \\
        \midrule
        Running time (hours) & 0.3288     & 0.3693     & 0.4173     \\
        \bottomrule
    \end{tabular}

    \caption{\textit{Boussinesq-Burgers equations:} Relative $L_2$ errors and running time of different methods.}
    \label{table:BBE}
\end{table}

In Figure \ref{fig:BBE_solution_u} and Figure \ref{fig:BBE_solution_v}, the numerical solutions and absolute errors of different methods are presented. We observe that IR-PINNs1 and IR-PINNs2 are still capable of enhancing the accuracy of the solutions when all methods employ the adaptive sampling strategy, as evidenced by Figure \ref{fig:BBE_error} and Table \ref{table:BBE}. However, it should be noted that if the evolution pattern of the solution in the temporal direction is relatively simple, the dominant source of error may stem from spatial discretization, such that the introduction of an integral-based regularization term results in modest improvement in accuracy. Additionally, the adaptive algorithm generates different spatial points at different times. Consequently, it is not feasible to discretize both residual terms simultaneously using the same set of points during the computation,  which leads to an increase in the overall computation time.

\begin{figure}[H]
    \centering
    \includegraphics[width=0.3\linewidth]{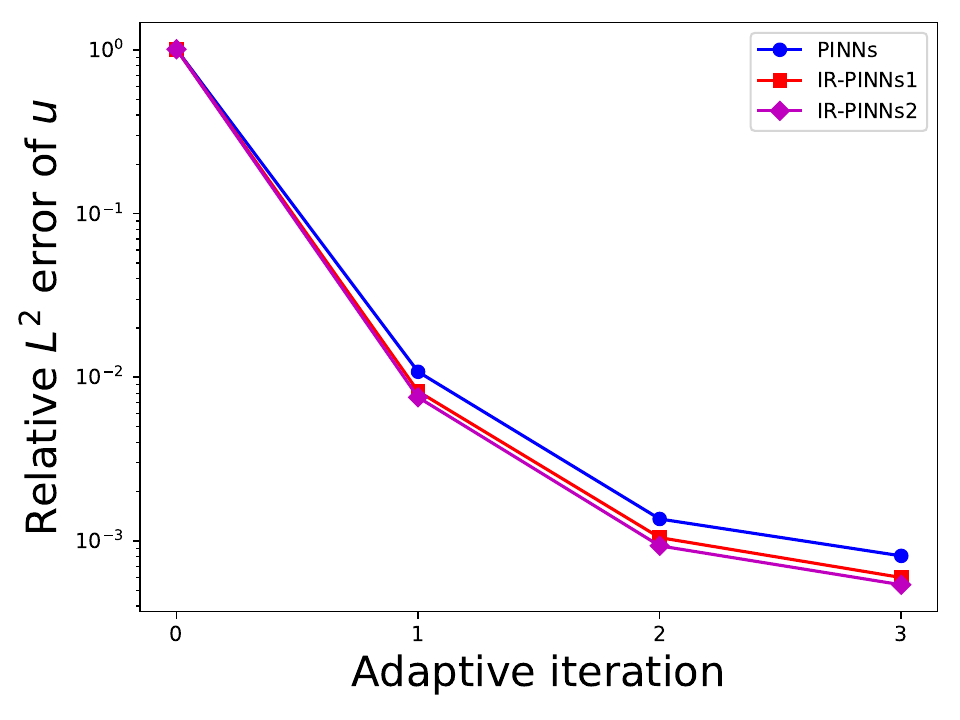}
    \includegraphics[width=0.3\linewidth]{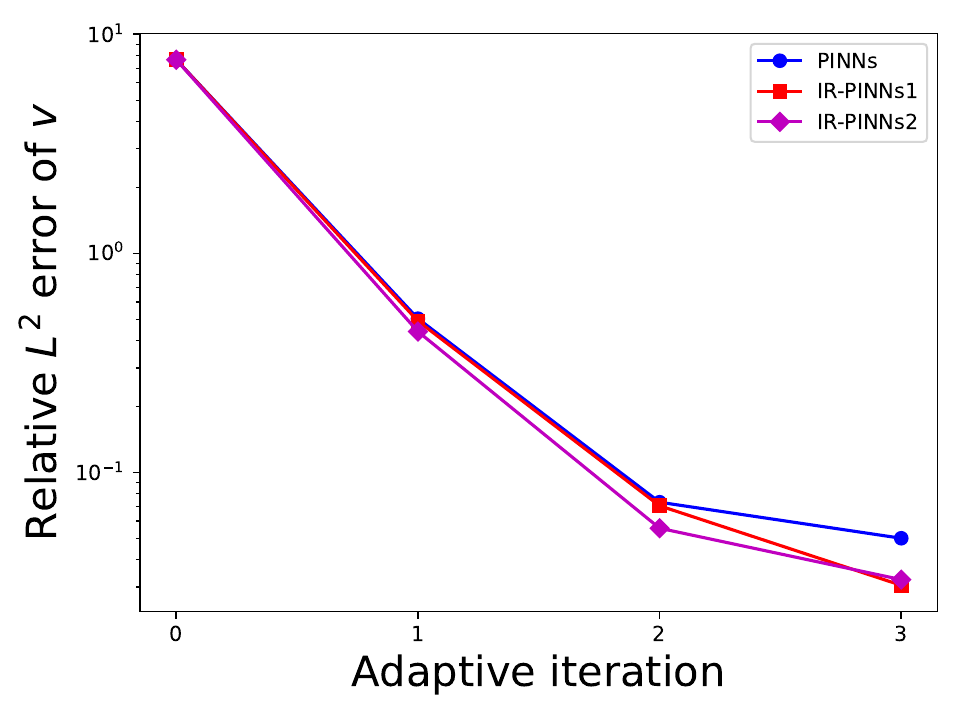}

    \caption{\textit{Boussinesq-Burgers equations:} Relative $L_2$ errors at different adaptive iterations.}
    \label{fig:BBE_adaptive_error}
\end{figure}

\begin{figure}[H]
    \centering
    \begin{subfigure}[t]{0.3\linewidth}
        \centering
        \includegraphics[width=\linewidth]{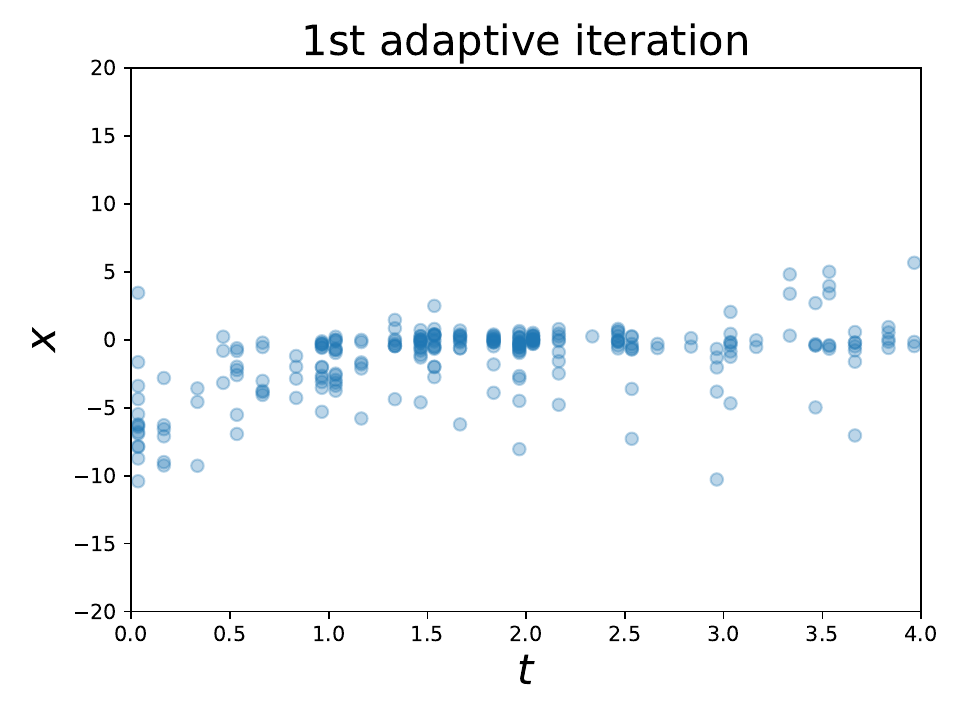}
        \includegraphics[width=\linewidth]{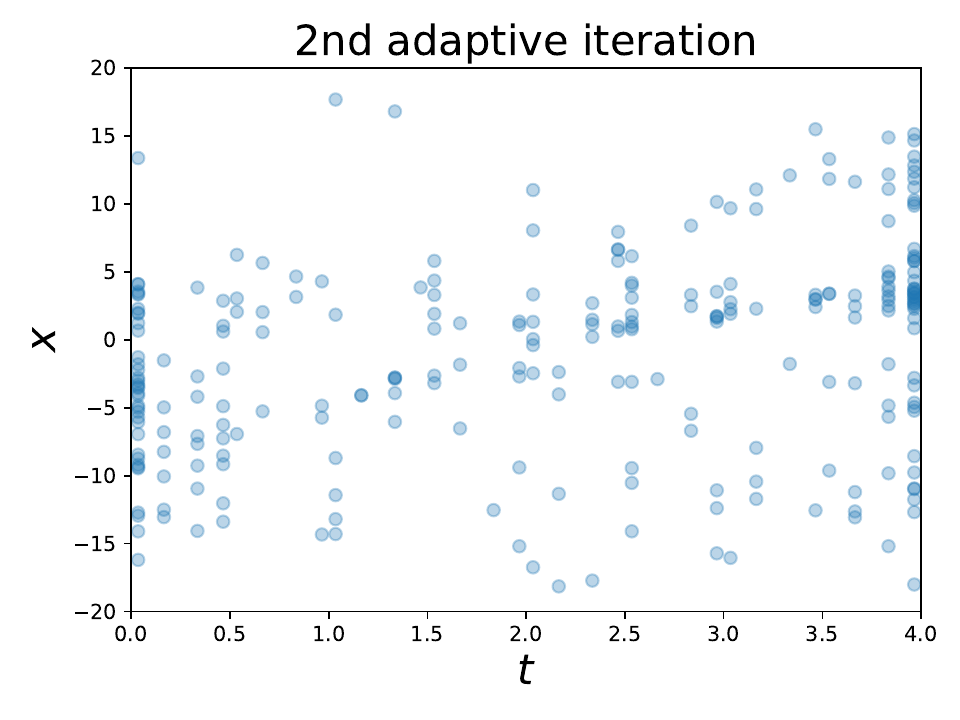}
        \caption{PINNs}
    \end{subfigure}
    \begin{subfigure}[t]{0.3\linewidth}
        \centering
        \includegraphics[width=\linewidth]{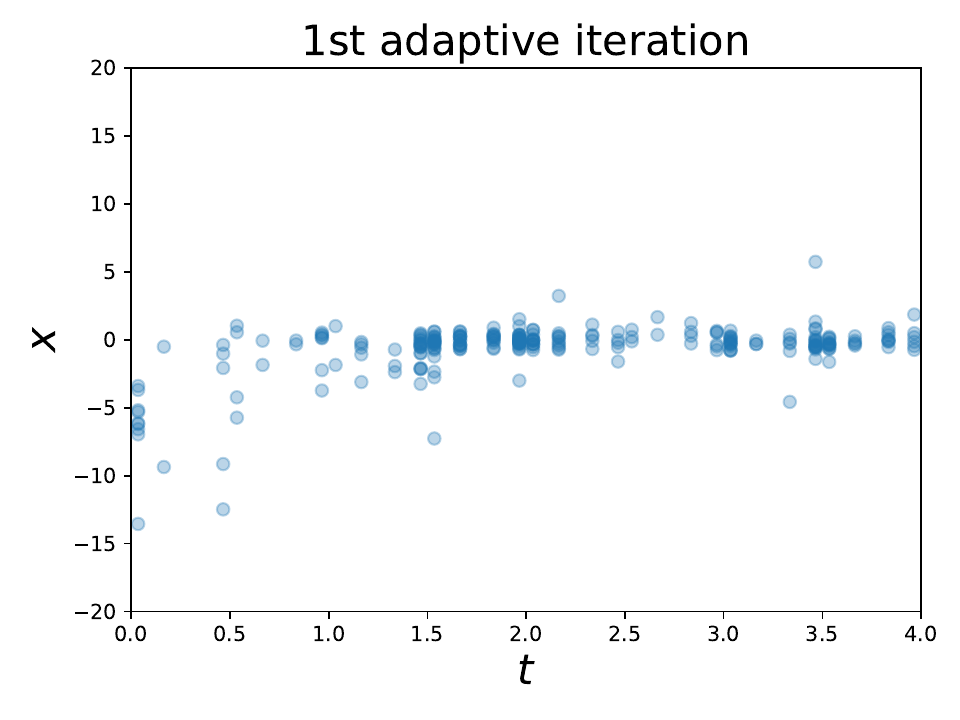}
        \includegraphics[width=\linewidth]{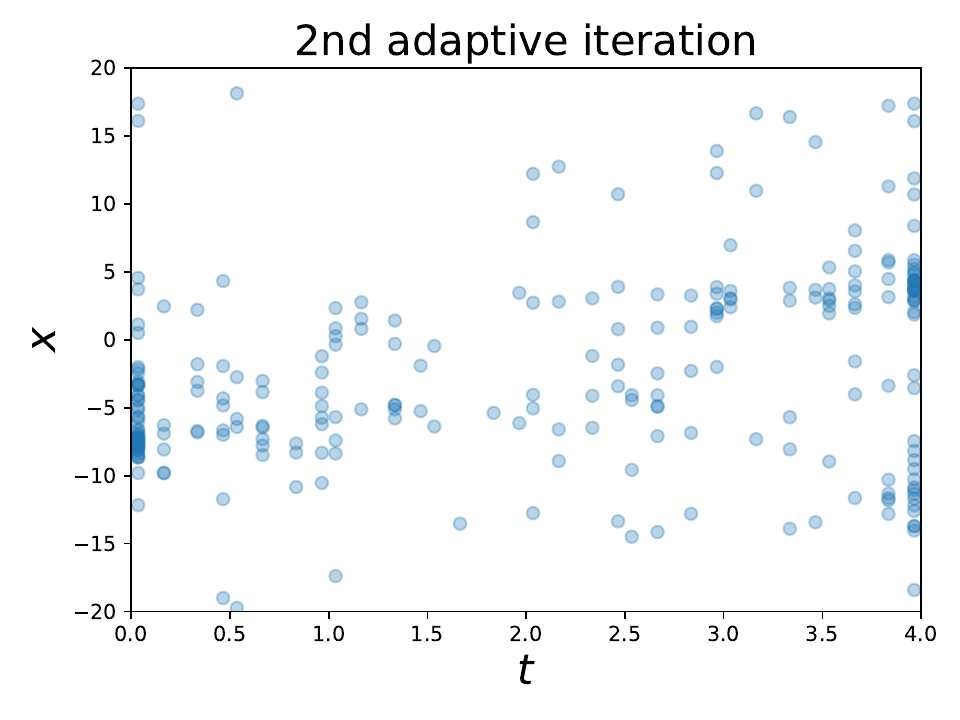}
        \caption{IR-PINNs1}
    \end{subfigure}
    \begin{subfigure}[t]{0.3\linewidth}
        \centering
        \includegraphics[width=\linewidth]{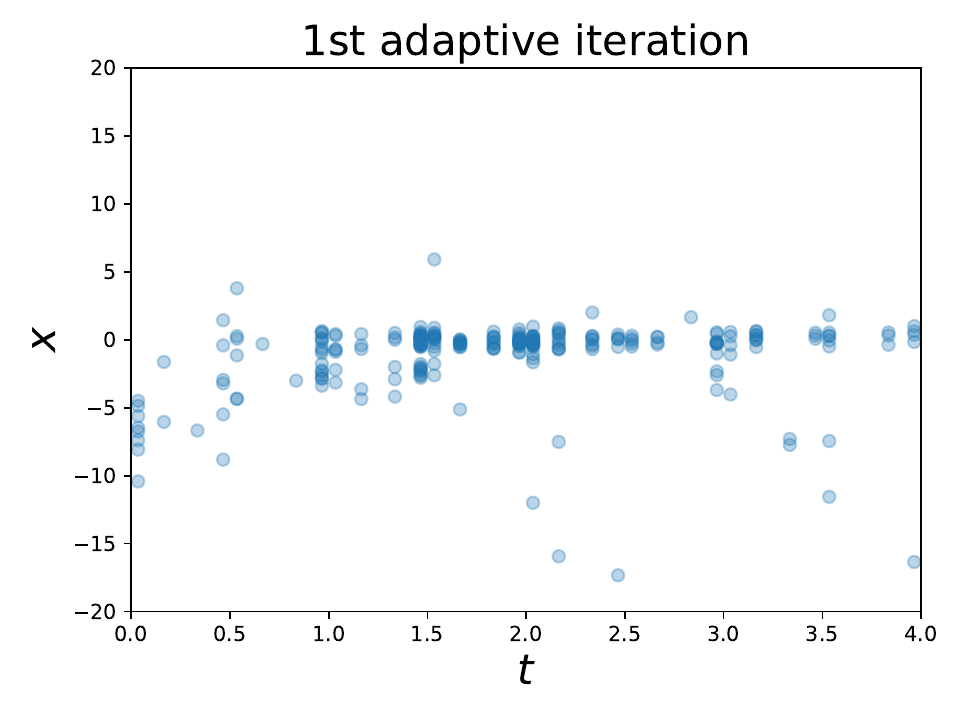}
        \includegraphics[width=\linewidth]{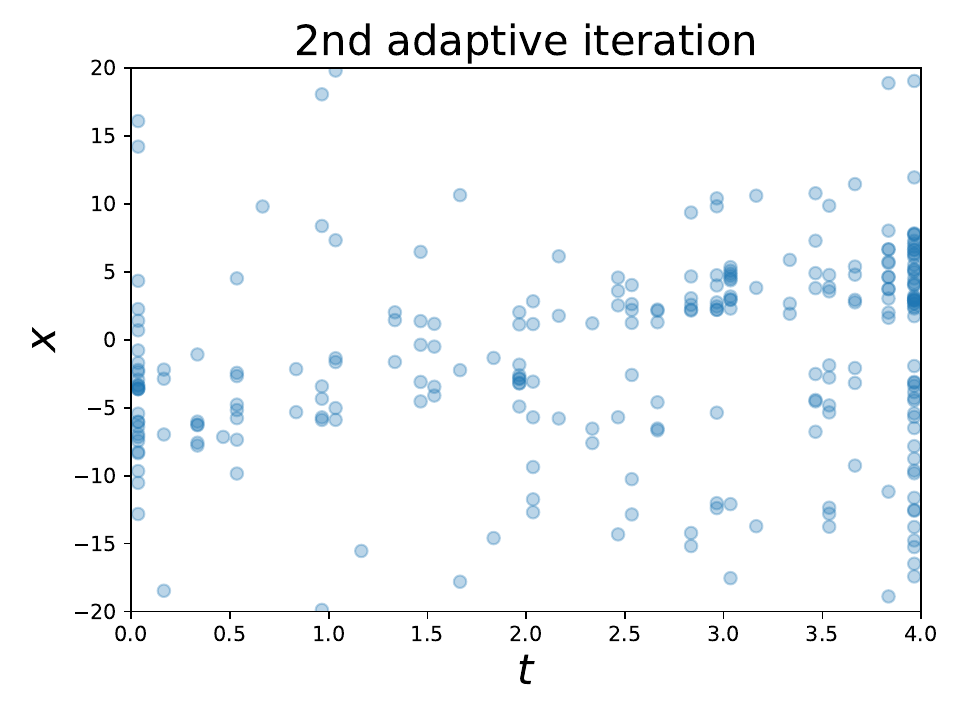}
        \caption{IR-PINNs2}
    \end{subfigure}

    \caption{\textit{Boussinesq-Burgers equations:} New samples generated by probability density model $p(t, \bm{x}; \theta_f)$ of two adaptive iterations.}
    \label{fig:BBE_adaptive_points}
\end{figure}

In Figure \ref{fig:BBE_adaptive_error} and Figure \ref{fig:BBE_adaptive_points}, we present the relative $L_2$ errors and newly generated points at different adaptive iterations. We can observe that the new training points are primarily concentrated in the middle region where the gradient is larger. As the training progresses, the residuals become more evenly distributed, resulting in a more uniform spatial distribution of the points.

\subsection{Time-dependent Fokker-Planck equation}
In this example, we consider an unbounded domain problem, time-dependent Fokker-Planck equation, which is a fundamental model in statistical physics and describes the evolution of the probability density function of a stochastic process. We start from the following stochastic differential equation
\begin{equation}
    \mathrm{d}\bm{X}_t = -\bm{\mu}(t, \bm{X}_t) \mathrm{d}t + \bm{\sigma}(t, \bm{X}_t) \mathrm{d}W_t,
    \label{Fokker-Planck_SDE}
\end{equation}
where $\bm{X}_t$ and $\bm{\mu}(t, \bm{X}_t)$ are $d$-dimensional random vectors, $\bm{\sigma}(t, \bm{X}_t)$ is a $d \times M$ matrix and $W_t$ is a $M$-dimensional standard Wiener process. The corresponding Fokker-Planck equation is given by
\begin{equation}
    \frac{\partial p(t, \bm{x})}{\partial t} = -\nabla \cdot \left(\bm{\mu}(t, \bm{x}) p(t, \bm{x})\right) + \frac{1}{2}\nabla \cdot \left(\nabla \cdot \left(\bm{\sigma}(t, \bm{x}) \bm{\sigma}(t, \bm{x})^T p(t, \bm{x})\right)\right),
    \label{Fokker-Planck}
\end{equation}
where $p(t, \bm{x})$ is the probability density function of $\bm{X}_t$, $\bm{\mu}$ is the drift term and $\bm{\sigma}$ is the diffusion term. Generally, the solution of the time-dependent Fokker-Planck equation is defined over an unbounded domain with the following boundary condition:
\begin{equation}
    p(t, \bm{x}) \rightarrow 0, \quad \text{as} \ \Vert\bm{x}\Vert \rightarrow \infty.
    \label{Fokker-Planck_BC}
\end{equation}
Moreover, since the solution is a probability density function, it must satisfy the non-negativity and normalization constraints:
\begin{equation}
    p(t, \bm{x}) \geq 0, \quad \int_{\mathbb{R}^d} p(t, \bm{x}) \mathrm{d}\bm{x} \equiv 1.
    \label{Fokker-Planck_constraints}
\end{equation}

We consider a two-dimensional nonlinear Fokker-Planck equation with the following drift and diffusion terms
\begin{equation}
    \bm{\mu} = (x_2, x_1-0.4x_2-0.1x_1^3), \quad \bm{\sigma} = \mathrm{diag}(0, \sqrt{0.8})
\end{equation}
with the initial condition
\begin{equation}
    p(0, \bm{x}) = \frac{1}{2\pi\vert \bm{\Sigma}_0 \vert^{1/2}} \exp\left(-\frac{1}{2}(\bm{x}-\bm{\mu}_0)^T \bm{\Sigma}_0^{-1} (\bm{x}-\bm{\mu}_0)\right),
\end{equation}
where $\bm{\mu}_0 = (0, 5)$ and $\bm{\Sigma}_0 = \rm{I}_2$. The temporal domain is $[0, 10]$ and we obtain the reference solution solution by the ADI scheme \cite{pichler2013numerical} in a truncated domain $[-10, 10]\times[-10, 10]$, with mesh size $\delta t = 0.005$ and $\delta h = 0.01$.

Here we directly approximate the solution of time-dependent Fokker-Planck equation by temporal normalizing flow \cite{feng2022solving}, a probability density model, which satisfies the properties in \eqref{Fokker-Planck_BC} and \eqref{Fokker-Planck_constraints} automatically. Since it is a long-time simulation, we employ the time-marching strategy and partition the temporal domain $[0, 10]$ into $10$ subdomains. For each subdomain, we use the actnorm layer and $8$ affine coupling layers with each layer parameterized by a $2$-layer fully-connected neural network with $64$ neurons. The initial learning rate is $10^{-3}$, with an exponential decay rate of $0.9$ applied every $500$ training epochs. The maximum number of training epochs is $20,000$. We divide the subdomain into $N = 10$ subintervals and in each subinterval, the number of Gaussian quadrature points is $M = 4$. The initial loss weight coefficient is set to $\lambda_2 = 100$.

To further improve the accuracy of the solution, we first generate $1,000$ sample paths by the Euler-Maruyama method \cite{kloeden1992stochastic} with $\delta t = 0.01$ and pre-train the neural network model for $2,000$ epochs with the data to obtain a good initialization. The details of the pre-training approach are provided in the Appendix \ref{sec:pre-training}. Subsequently, we generate $100$ spatial points for each temporal collocation point by the pre-trained model and $1,000$ initial points from the initial condition. Then we proceed with adaptive training. The number of adaptive iterations is $2$ for each subdomain and for each adaptive iteration, we newly generate $10\times4\times100$ adaptive sample points. Note that temporal normalizing flow is inherently a probability density model, so we can directly perform adaptive sampling using this model.

\begin{figure}[H]
    \centering
    \includegraphics[width=0.3\linewidth]{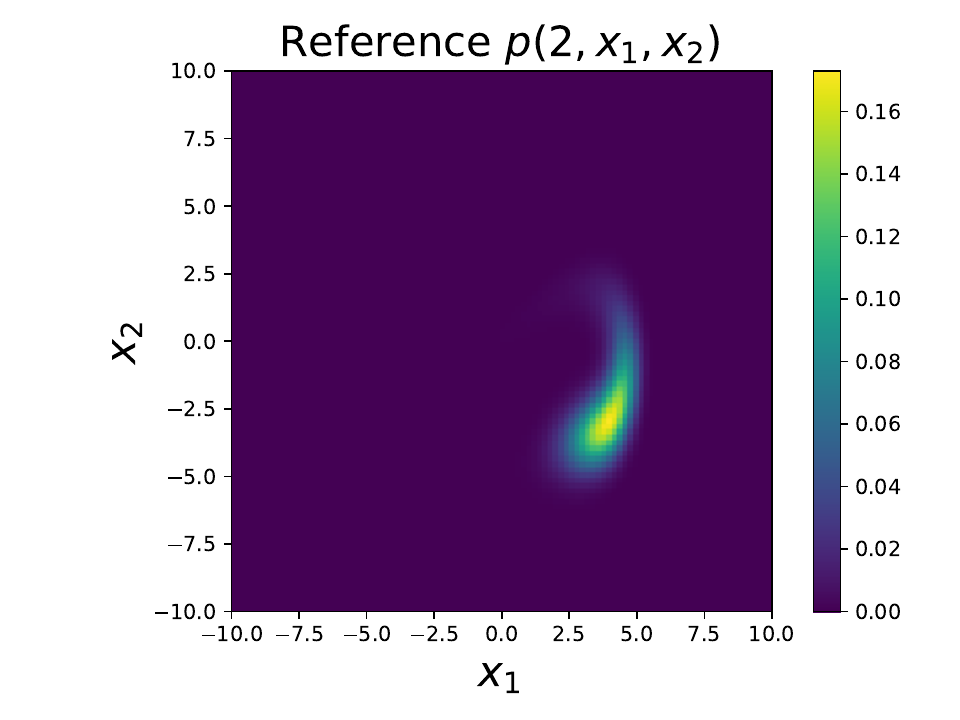}
    \includegraphics[width=0.3\linewidth]{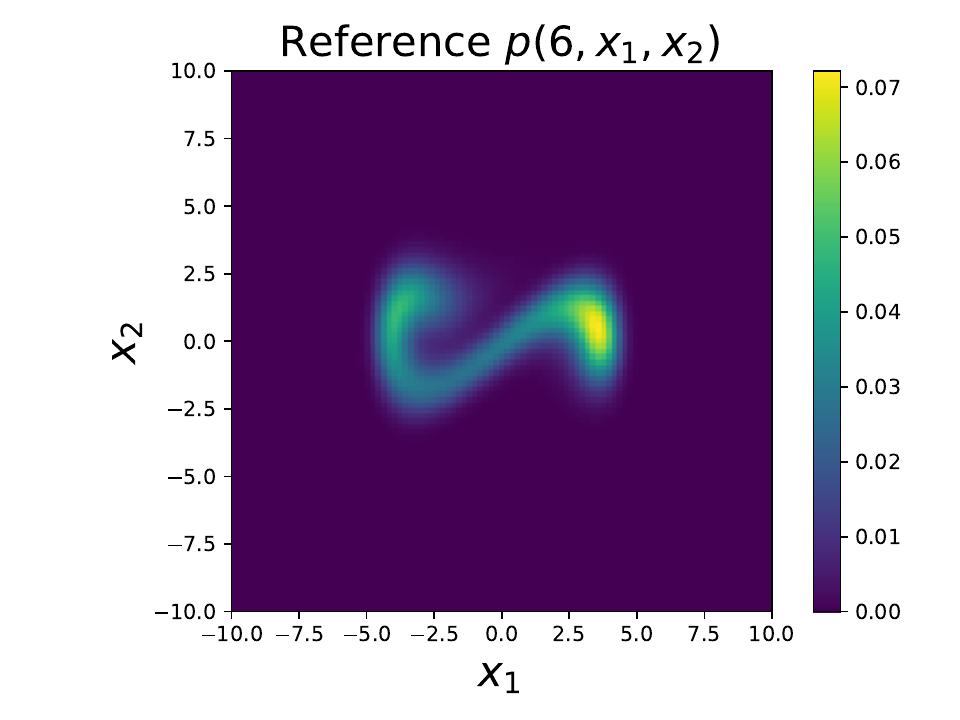}
    \includegraphics[width=0.3\linewidth]{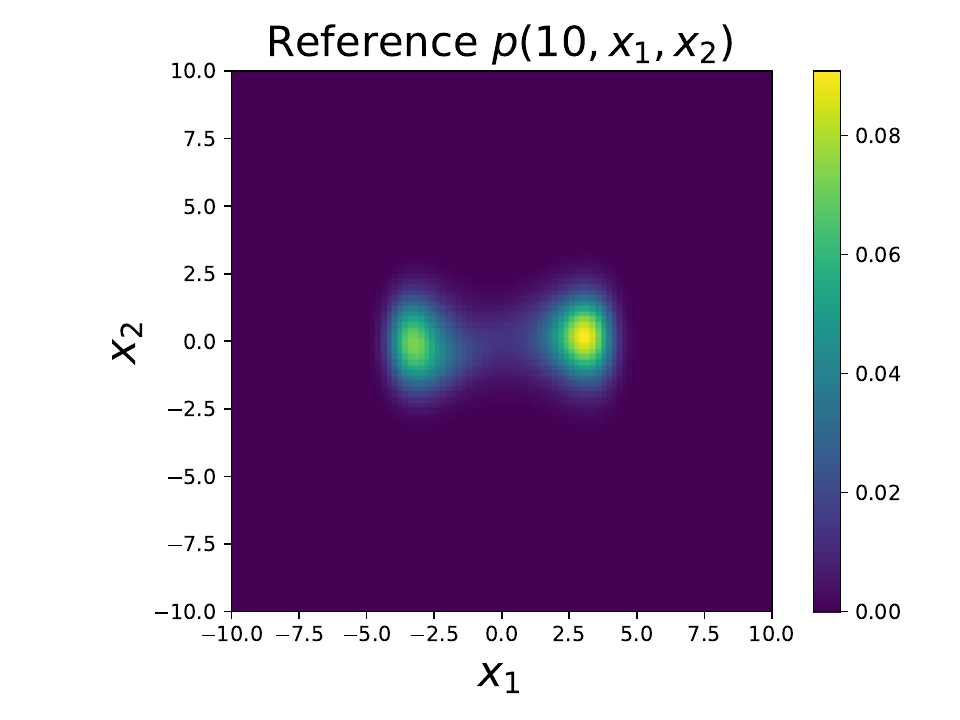}
    \caption{\textit{Time-dependent Fokker-Planck equation:} Reference solutions at $t = 2, 6, 10$.}
    \label{fig:FPE_reference_solution}
\end{figure}

\begin{figure}[H]
    \centering
    \includegraphics[width=0.3\linewidth]{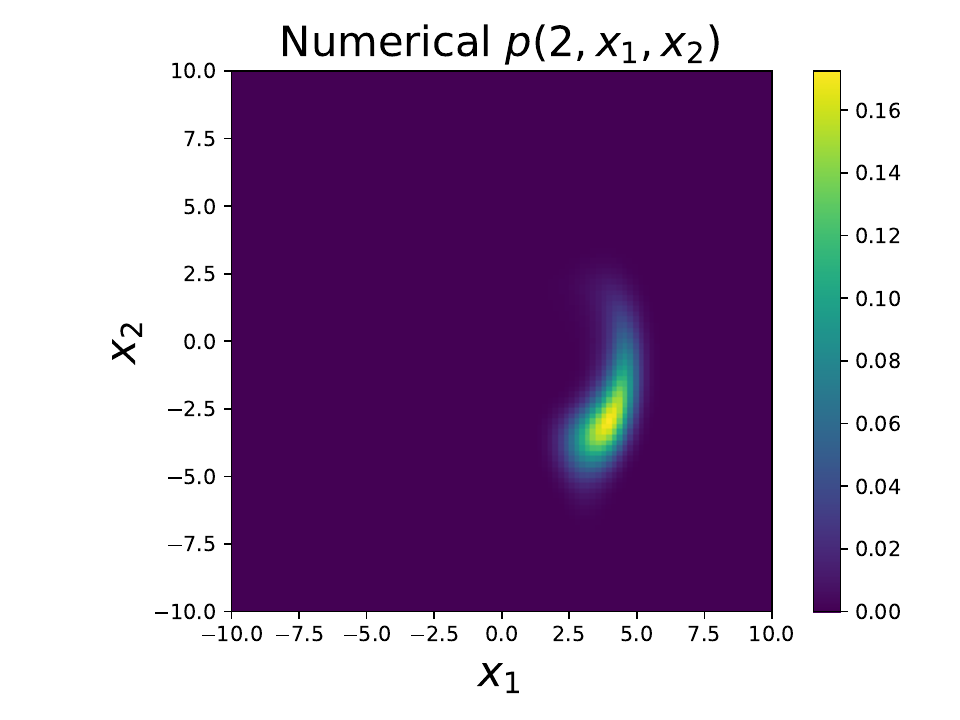}
    \includegraphics[width=0.3\linewidth]{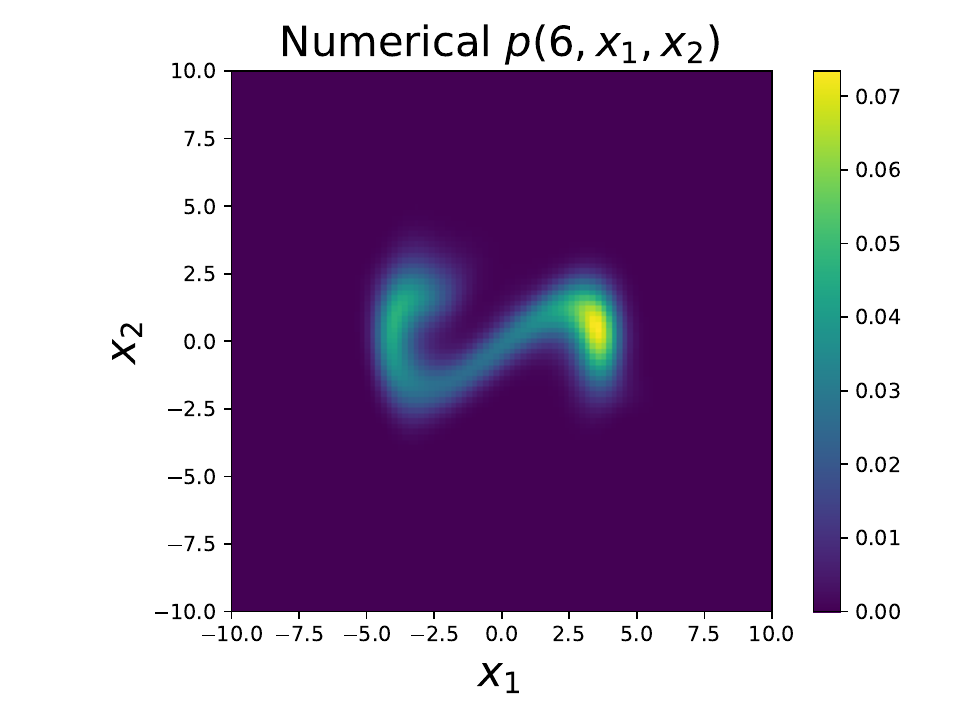}
    \includegraphics[width=0.3\linewidth]{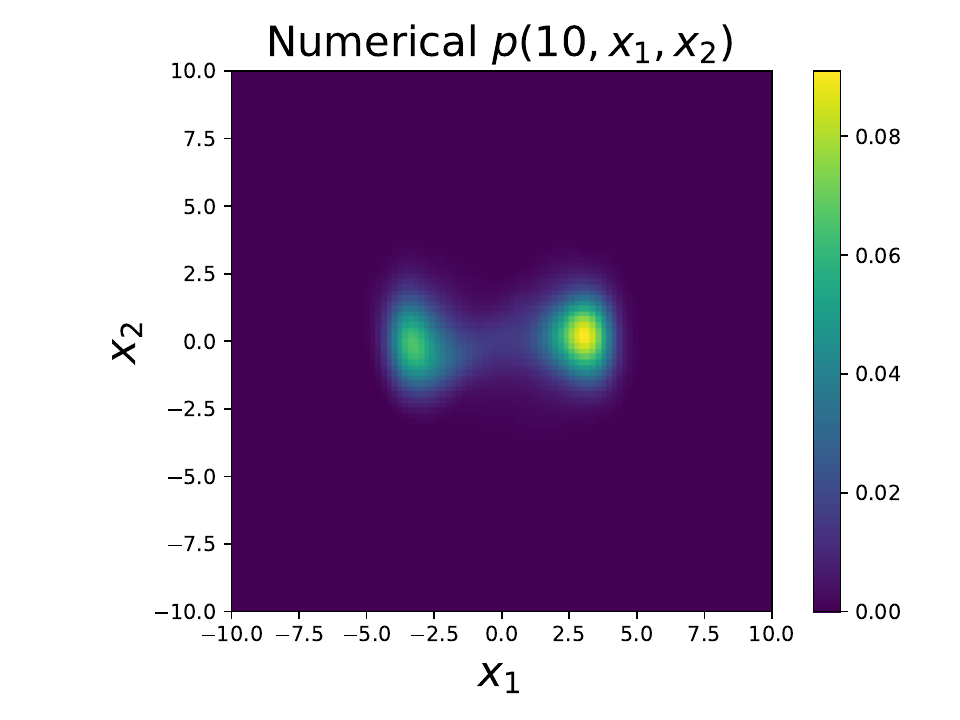}

    \includegraphics[width=0.3\linewidth]{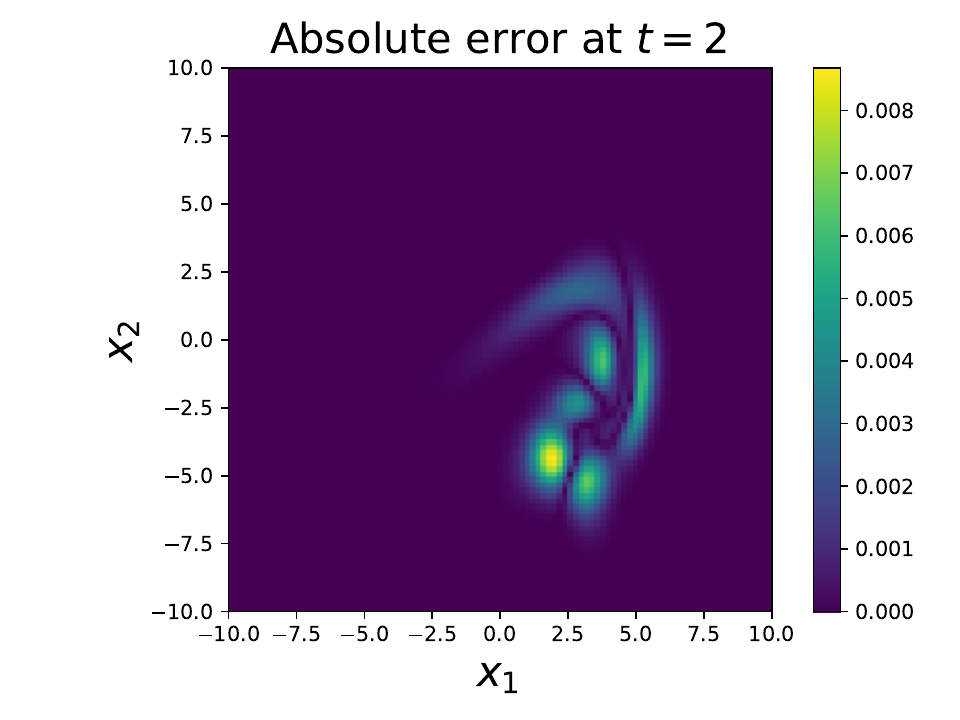}
    \includegraphics[width=0.3\linewidth]{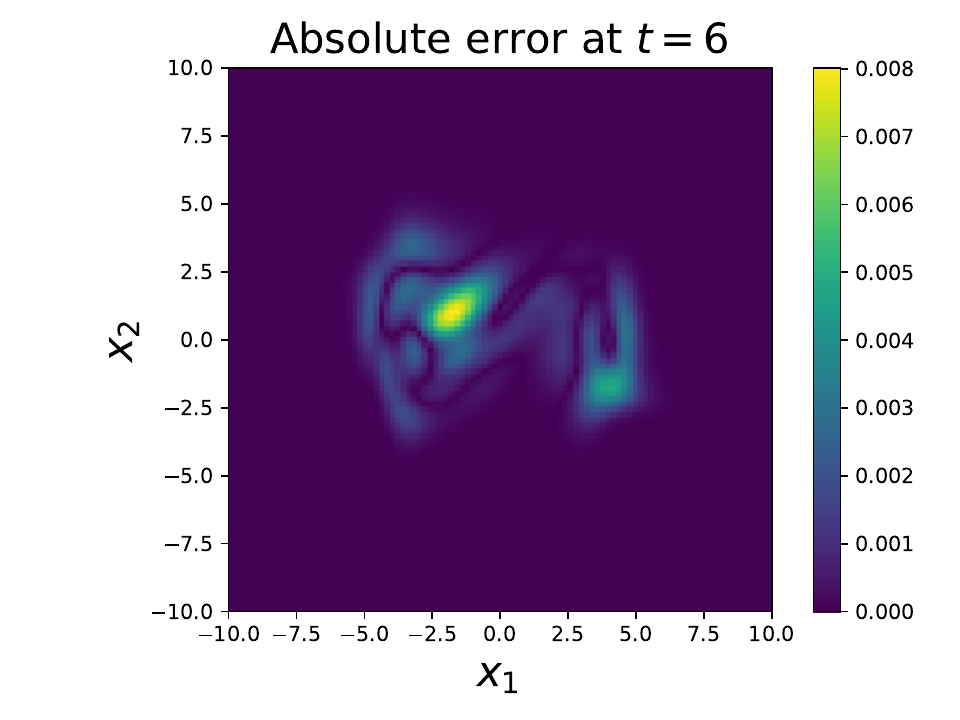}
    \includegraphics[width=0.3\linewidth]{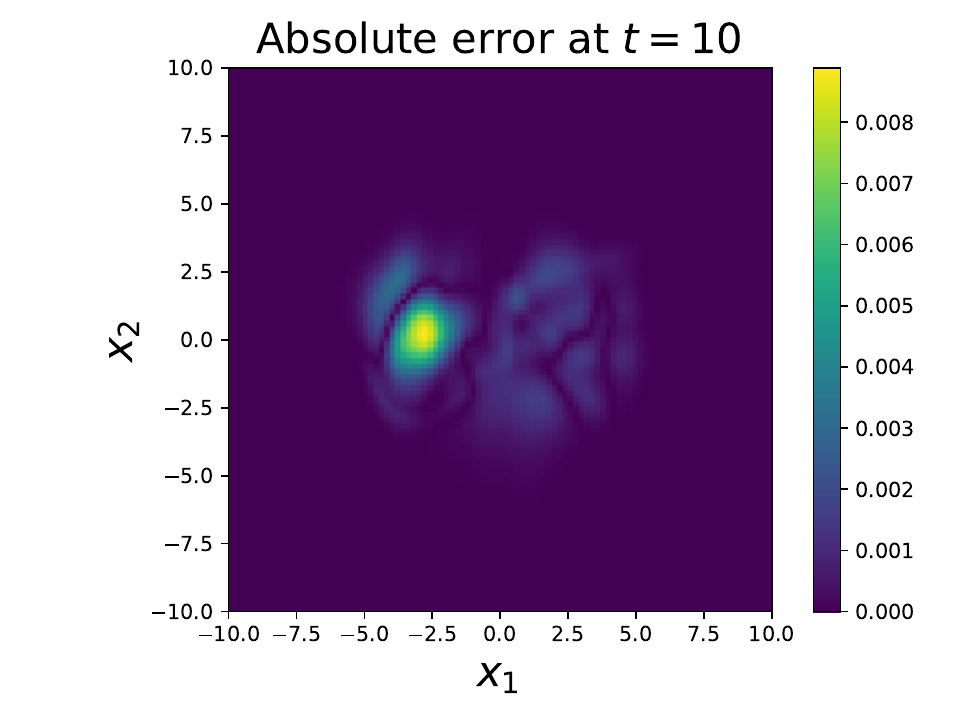}
    \caption{\textit{Time-dependent Fokker-Planck equation:} Numerical solutions and absolute errors of PINNs at $t = 2, 6, 10$.}
    \label{fig:FPE_solution_PINNs}
\end{figure}

\begin{figure}[H]
    \centering
    \includegraphics[width=0.3\linewidth]{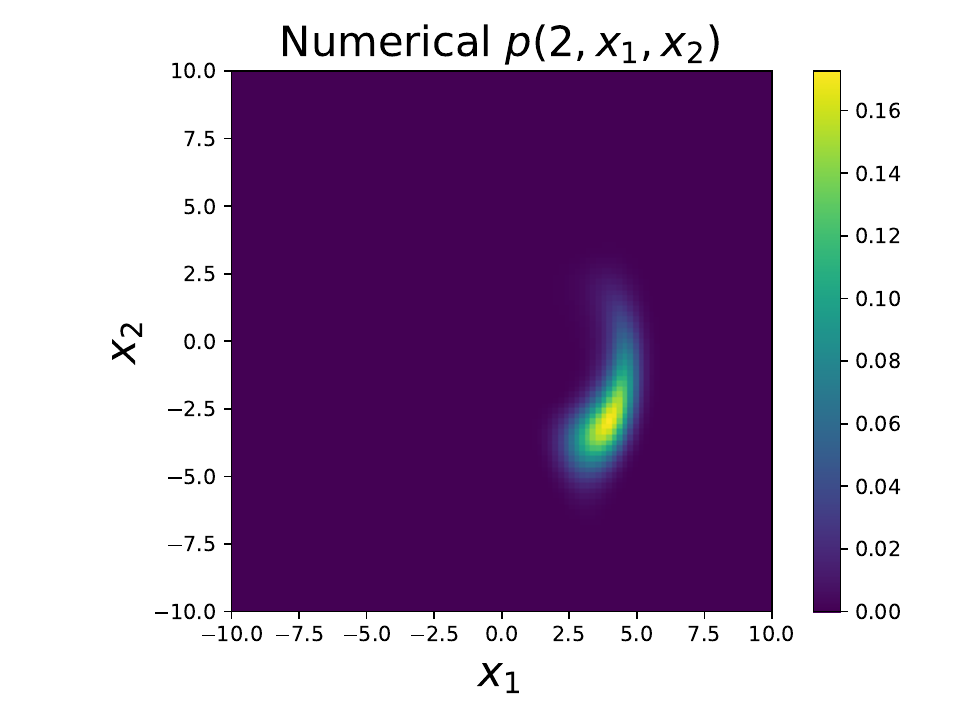}
    \includegraphics[width=0.3\linewidth]{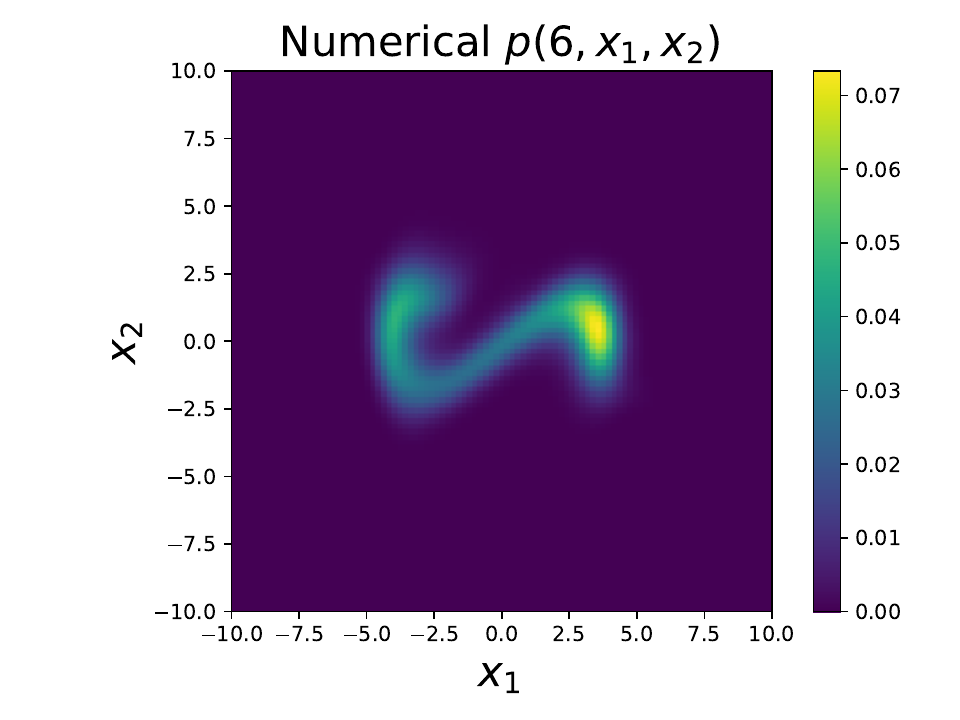}
    \includegraphics[width=0.3\linewidth]{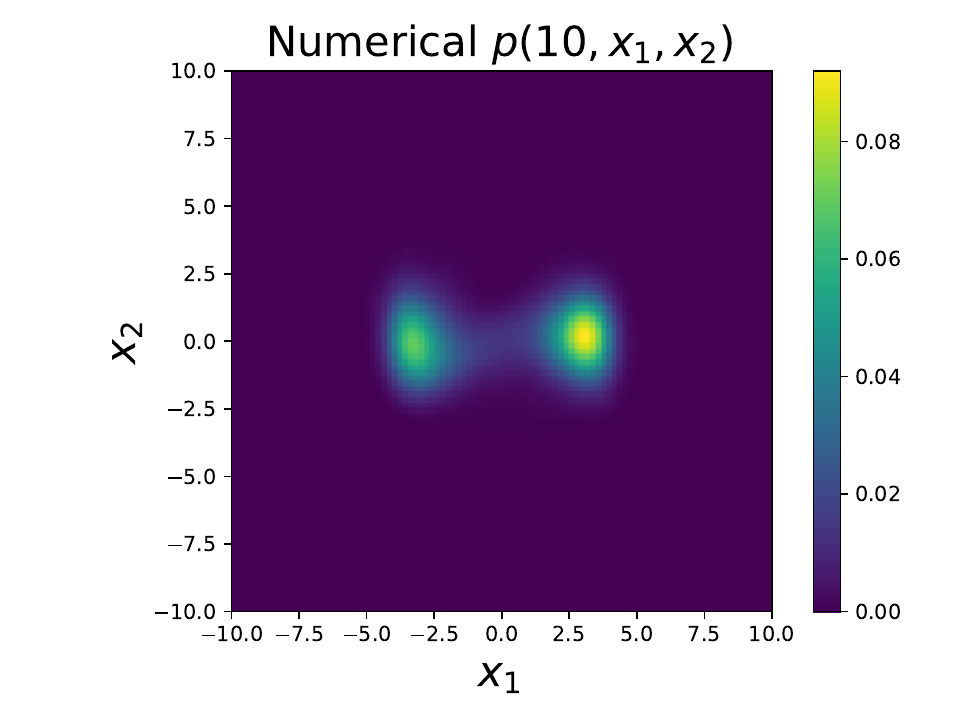}

    \includegraphics[width=0.3\linewidth]{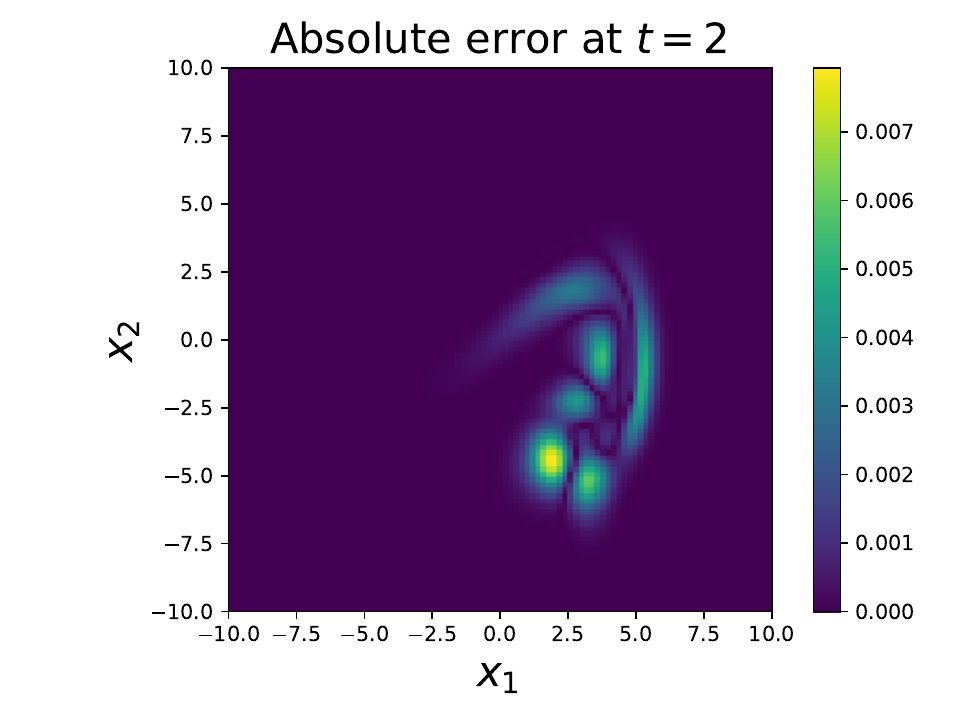}
    \includegraphics[width=0.3\linewidth]{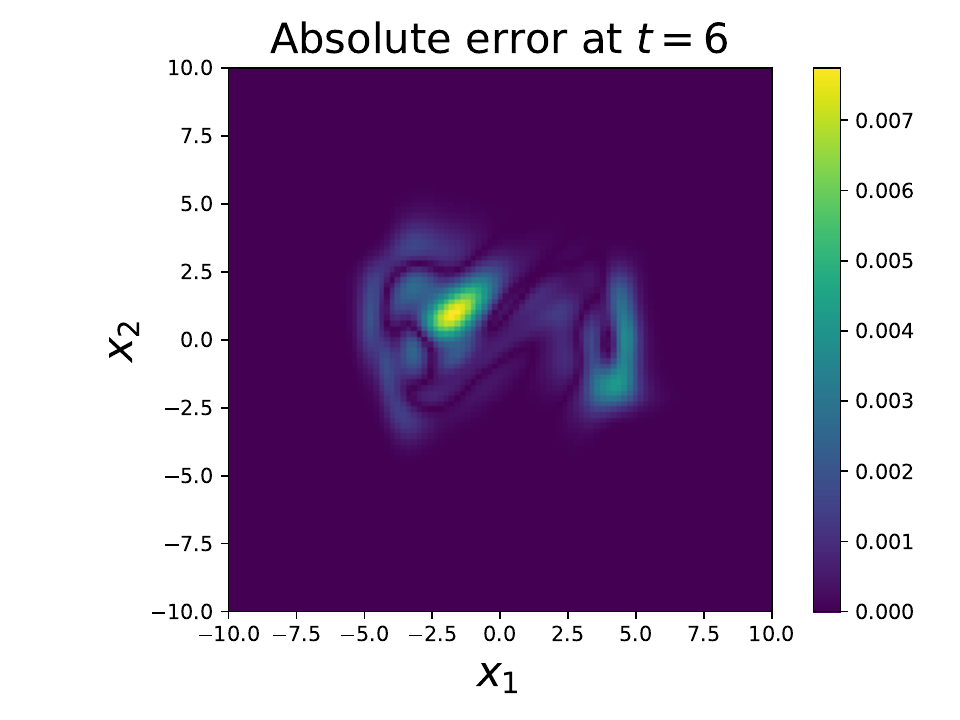}
    \includegraphics[width=0.3\linewidth]{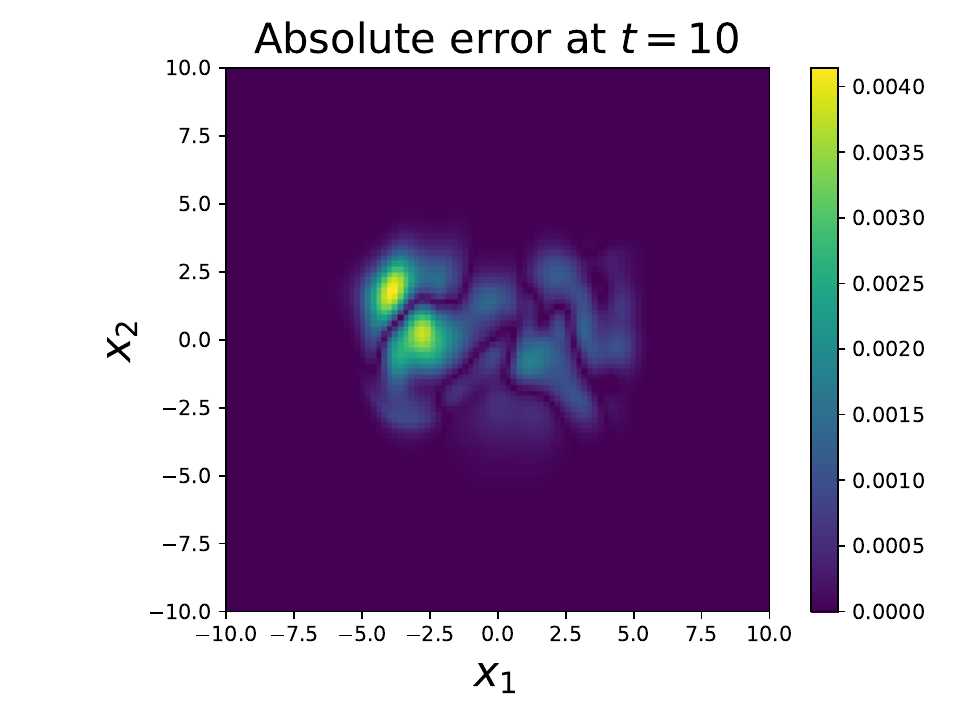}
    \caption{\textit{Time-dependent Fokker-Planck equation:} Numerical solutions and absolute errors of IR-PINNs1 at $t = 2, 6, 10$.}
    \label{fig:FPE_solution_IR-PINNs1}
\end{figure}

\begin{figure}[H]
    \centering
    \includegraphics[width=0.3\linewidth]{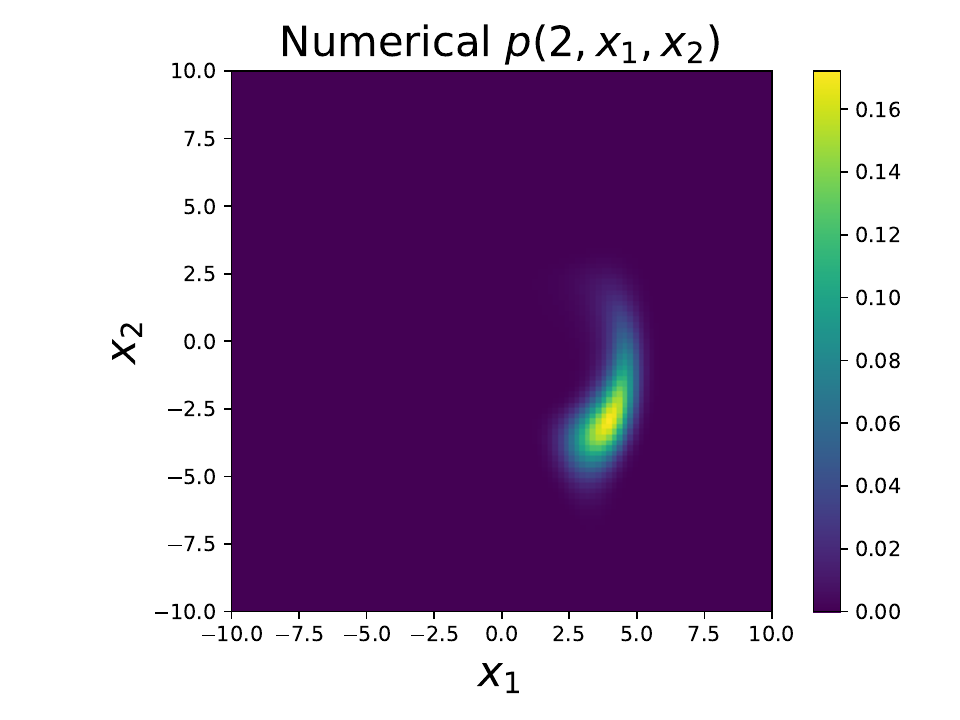}
    \includegraphics[width=0.3\linewidth]{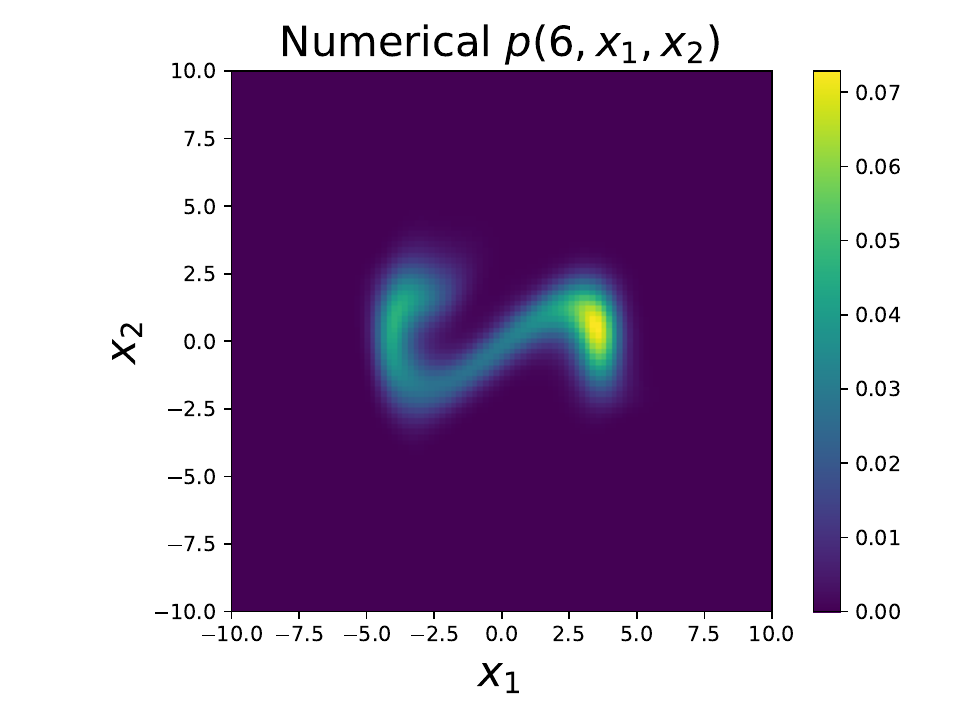}
    \includegraphics[width=0.3\linewidth]{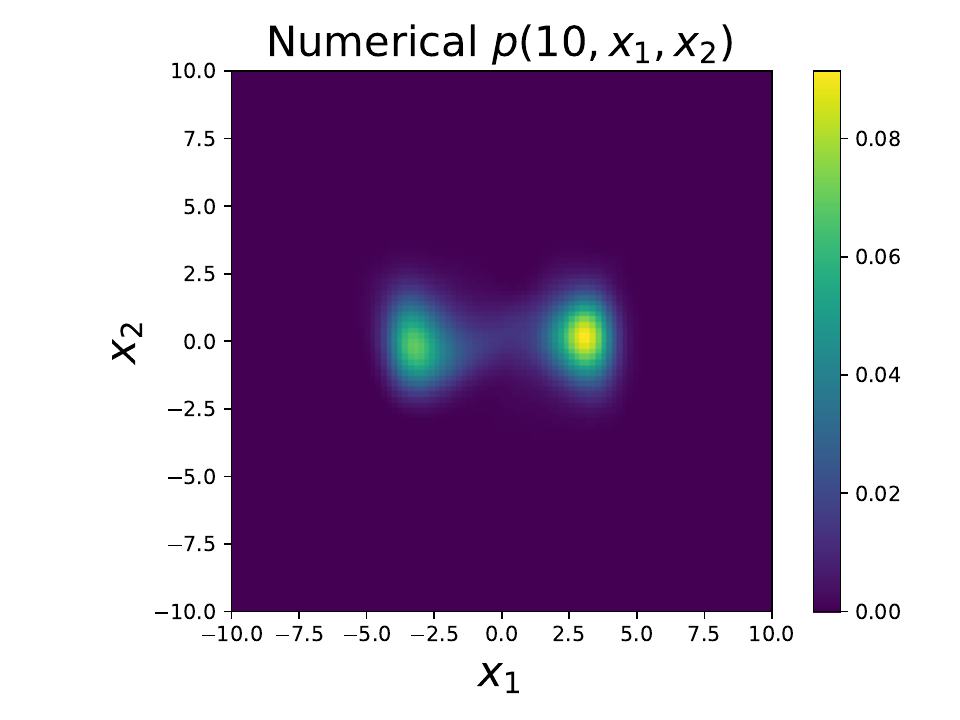}

    \includegraphics[width=0.3\linewidth]{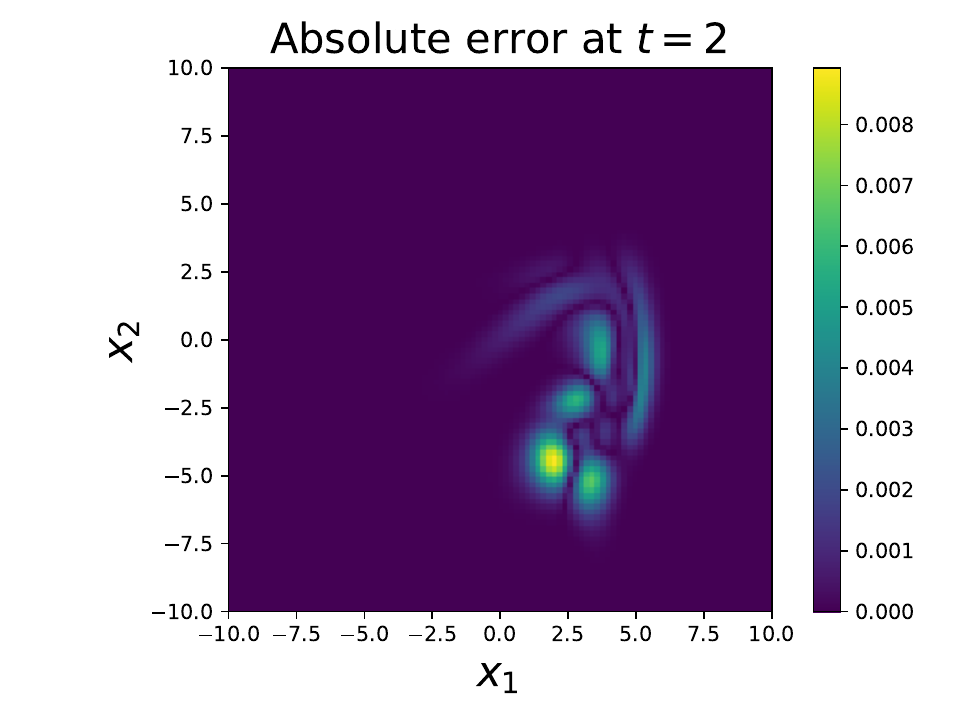}
    \includegraphics[width=0.3\linewidth]{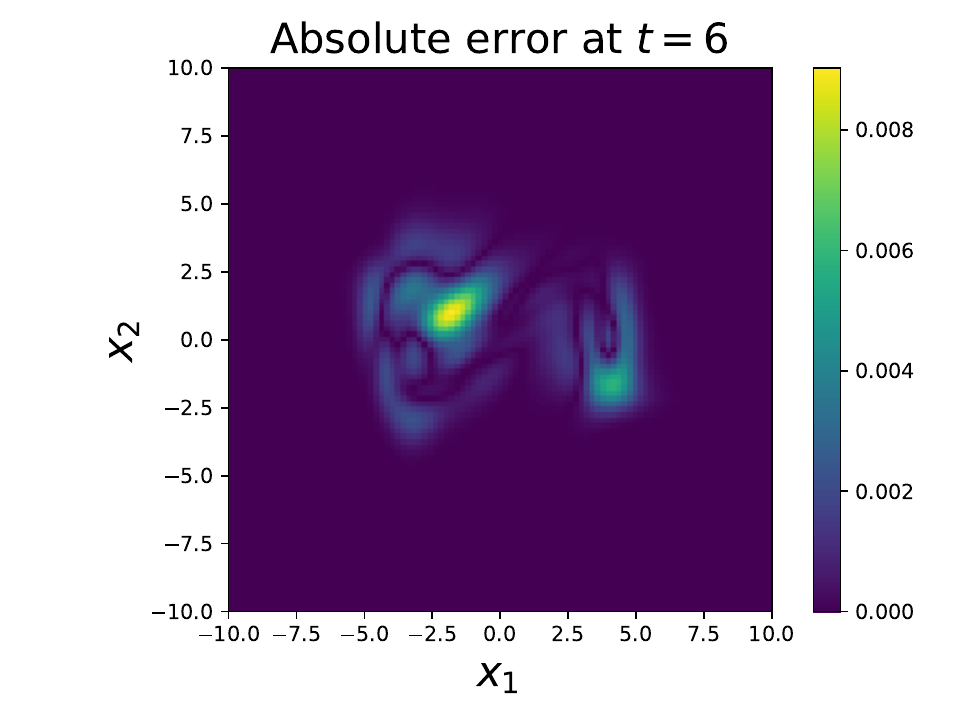}
    \includegraphics[width=0.3\linewidth]{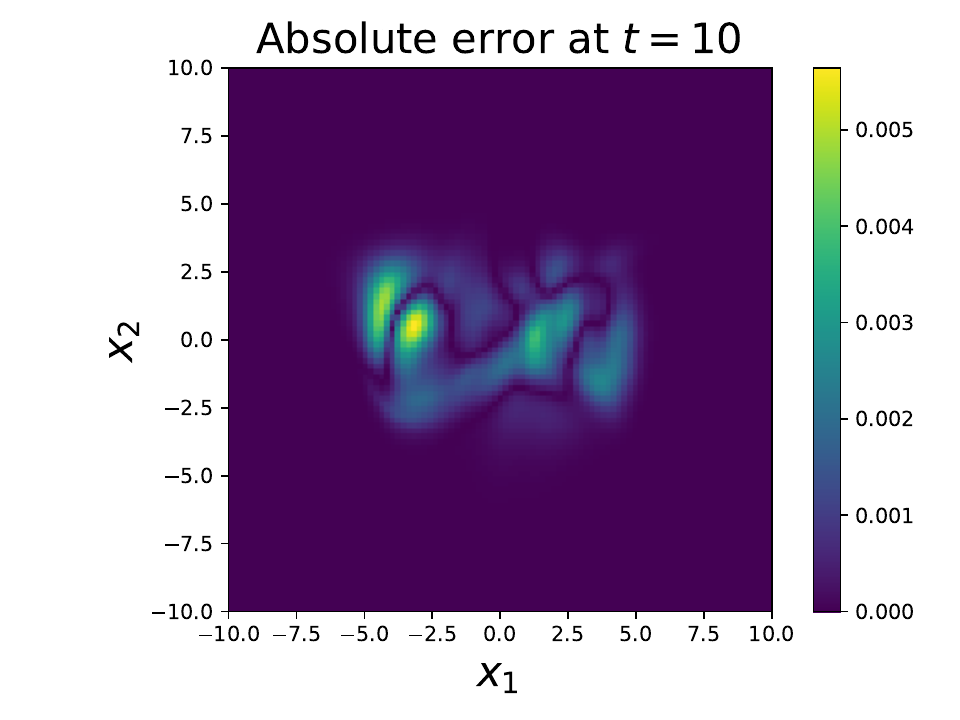}
    \caption{\textit{Time-dependent Fokker-Planck equation:} Numerical solutions and absolute errors of IR-PINNs2 at $t = 2, 6, 10$.}
    \label{fig:FPE_solution_IR-PINNs2}
\end{figure}

\begin{figure}[H]
    \centering
    \includegraphics[width=0.3\linewidth]{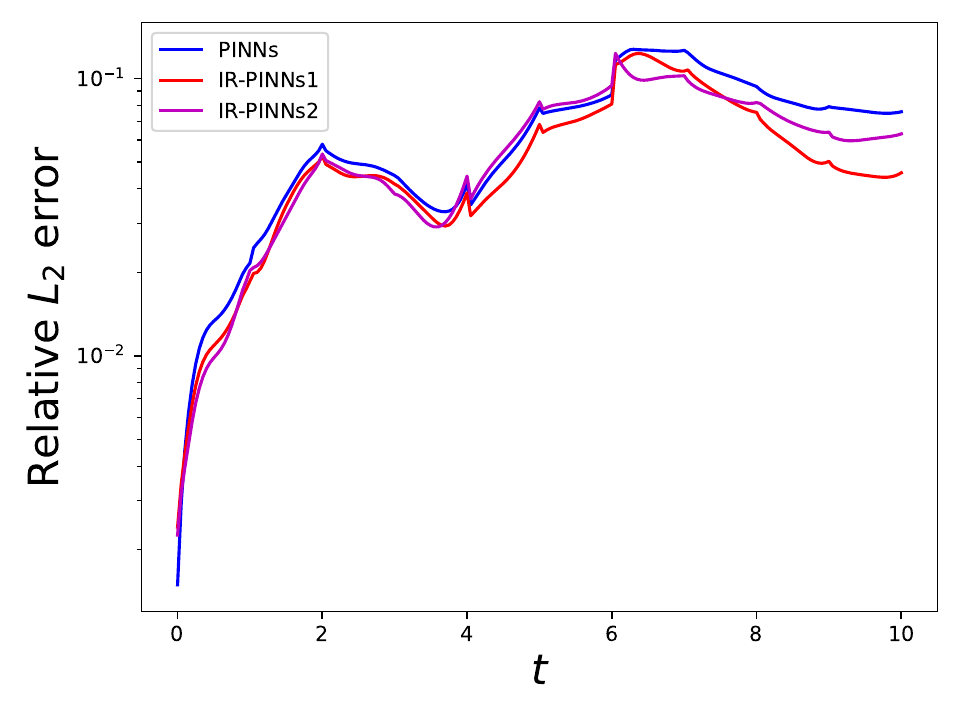}
    \includegraphics[width=0.3\linewidth]{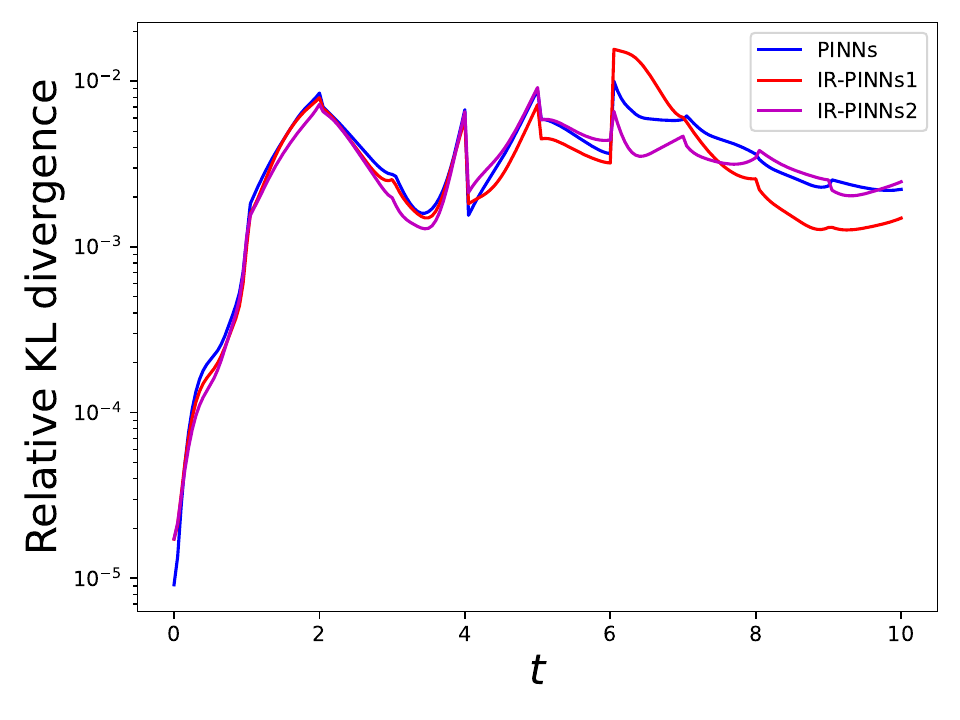}

    \caption{\textit{Time-dependent Fokker-Planck equation:} Relative $L_2$ errors and relative KL divergences over time.}
    \label{fig:FPE_error_KL}
\end{figure}

\begin{table}[H]
    \centering
    \begin{tabular}{cccc}
        \toprule
                               & PINNs      & IR-PINNs1  & IR-PINNs2  \\
        \midrule
        Relative $L_2$ error   & 6.3784e-02 & 5.3724e-02 & 5.6185e-02 \\
        \midrule
        Relative KL divergence & 2.0582e-03 & 2.0117e-03 & 1.8202e-03 \\
        \midrule
        Running time (hours)   & 9.612      & 21.06      & 24.47      \\
        \bottomrule
    \end{tabular}

    \caption{\textit{Time-dependent Fokker-Planck equation:} Relative $L_2$ errors and running time of different methods.}
    \label{table:FPE}
\end{table}

\begin{figure}[H]
    \centering
    \begin{subfigure}[t]{0.3\linewidth}
        \centering
        \includegraphics[width=\linewidth]{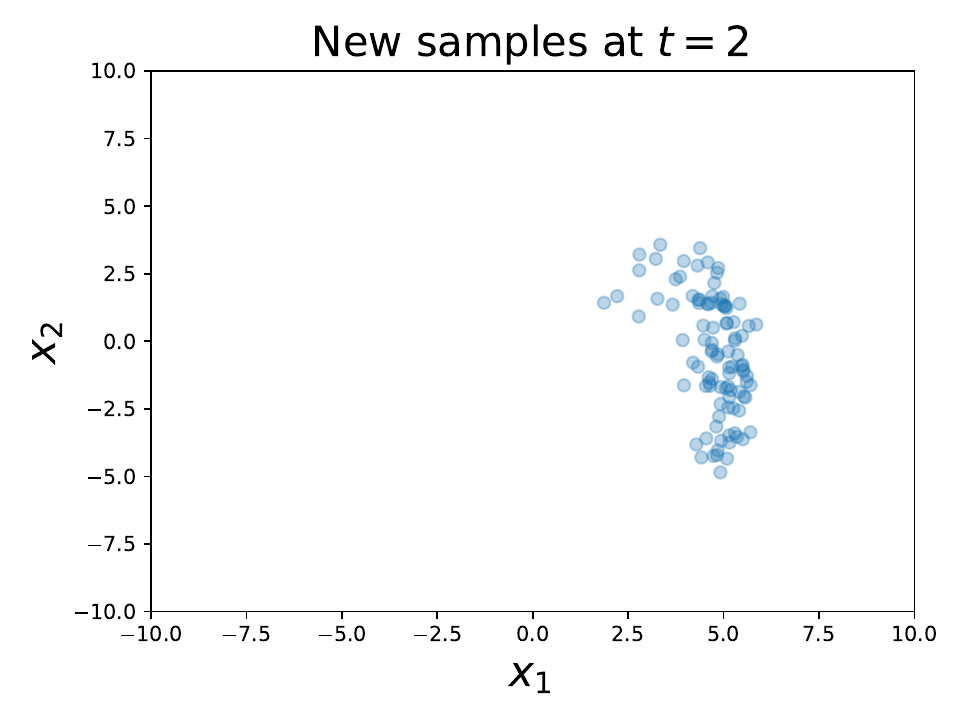}
        \includegraphics[width=\linewidth]{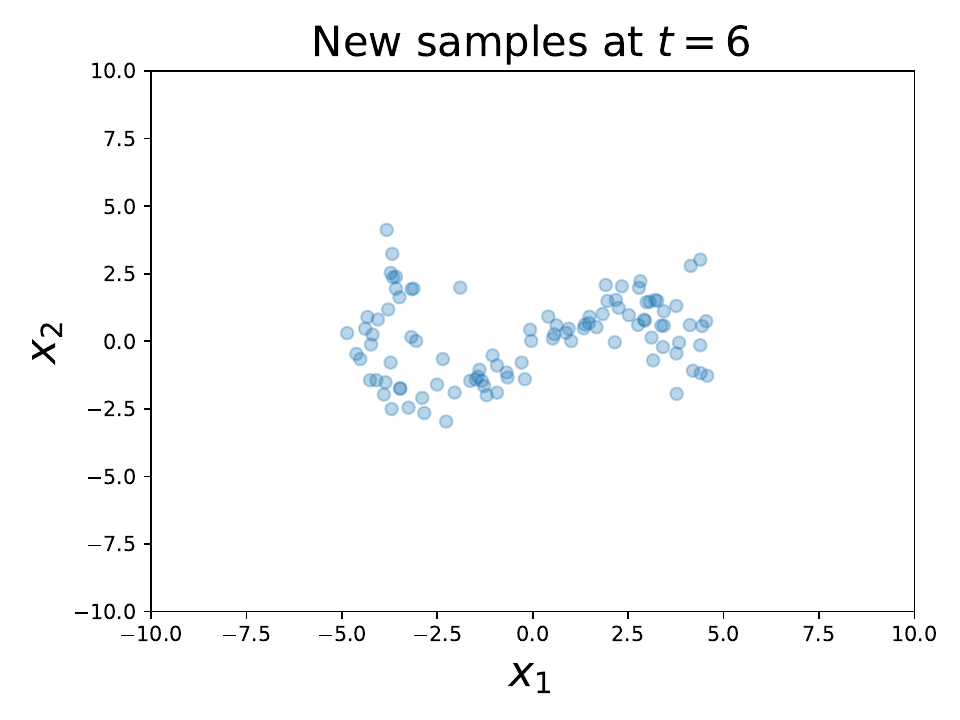}
        \includegraphics[width=\linewidth]{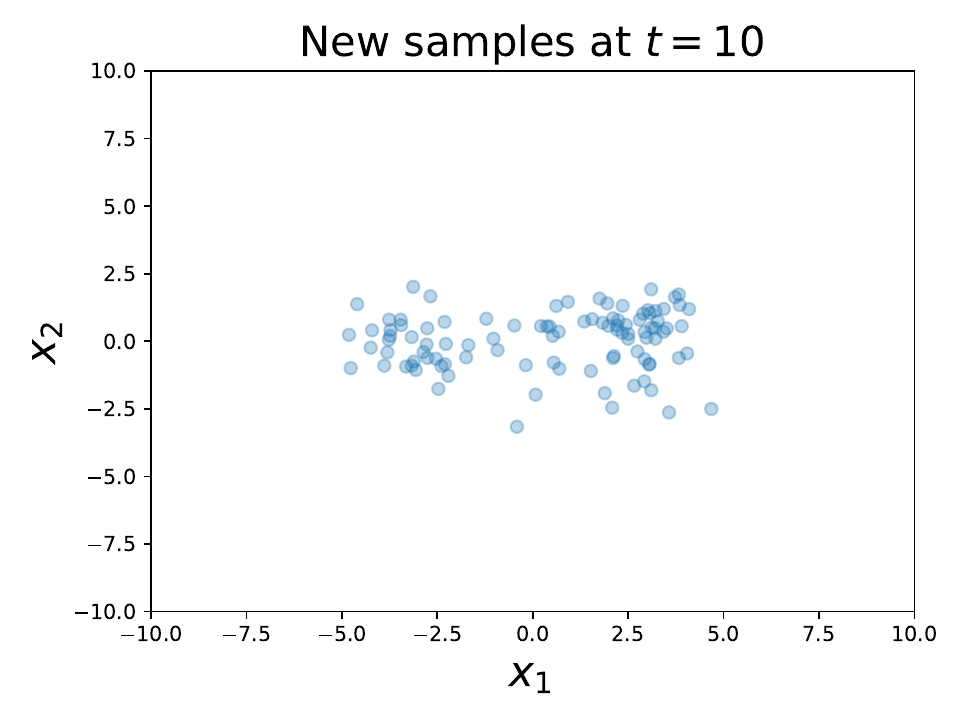}
        \caption{PINNs}
    \end{subfigure}
    \begin{subfigure}[t]{0.3\linewidth}
        \centering
        \includegraphics[width=\linewidth]{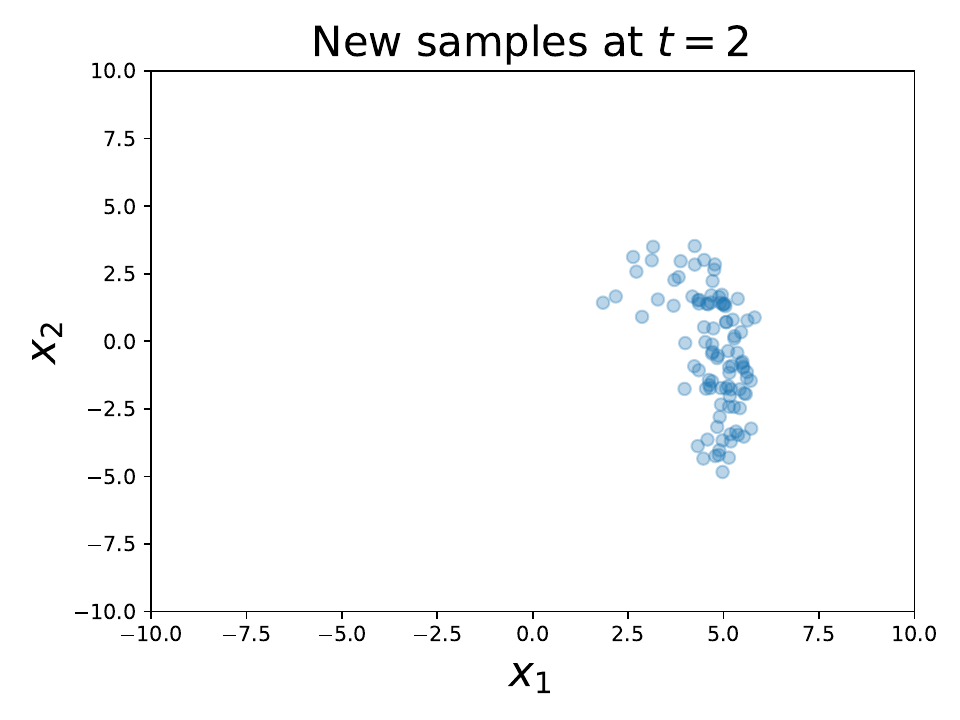}
        \includegraphics[width=\linewidth]{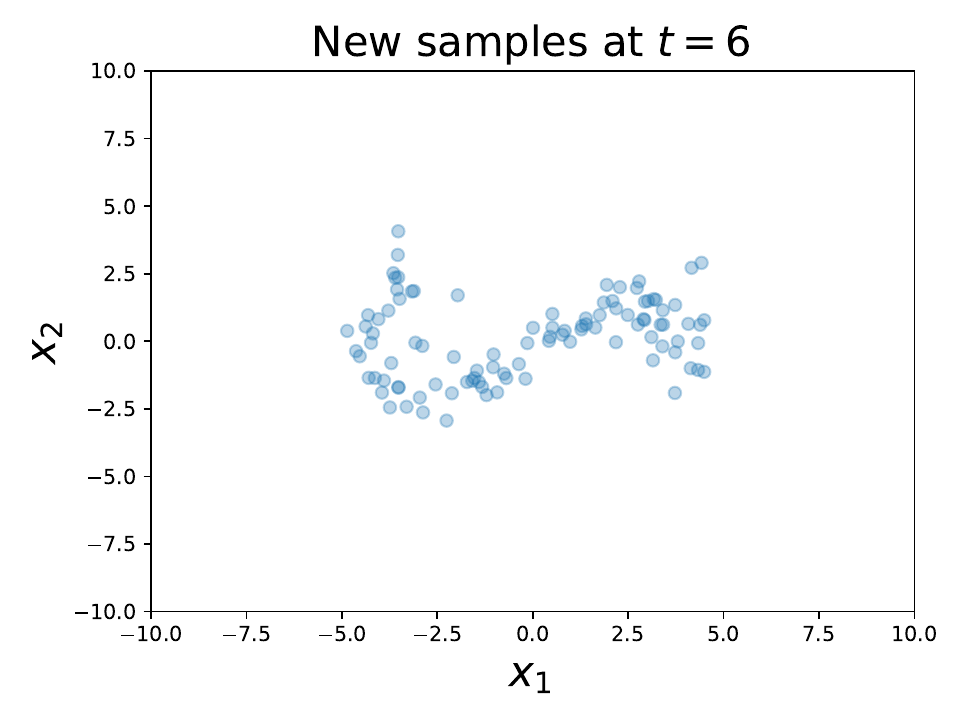}
        \includegraphics[width=\linewidth]{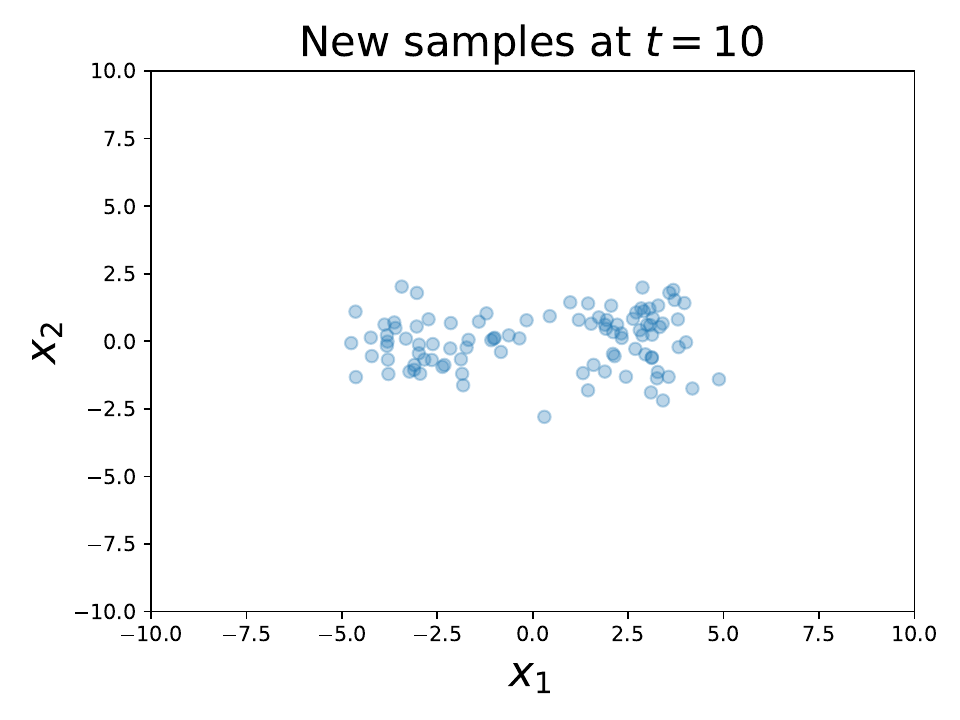}
        \caption{IR-PINNs1}
    \end{subfigure}
    \begin{subfigure}[t]{0.3\linewidth}
        \centering
        \includegraphics[width=\linewidth]{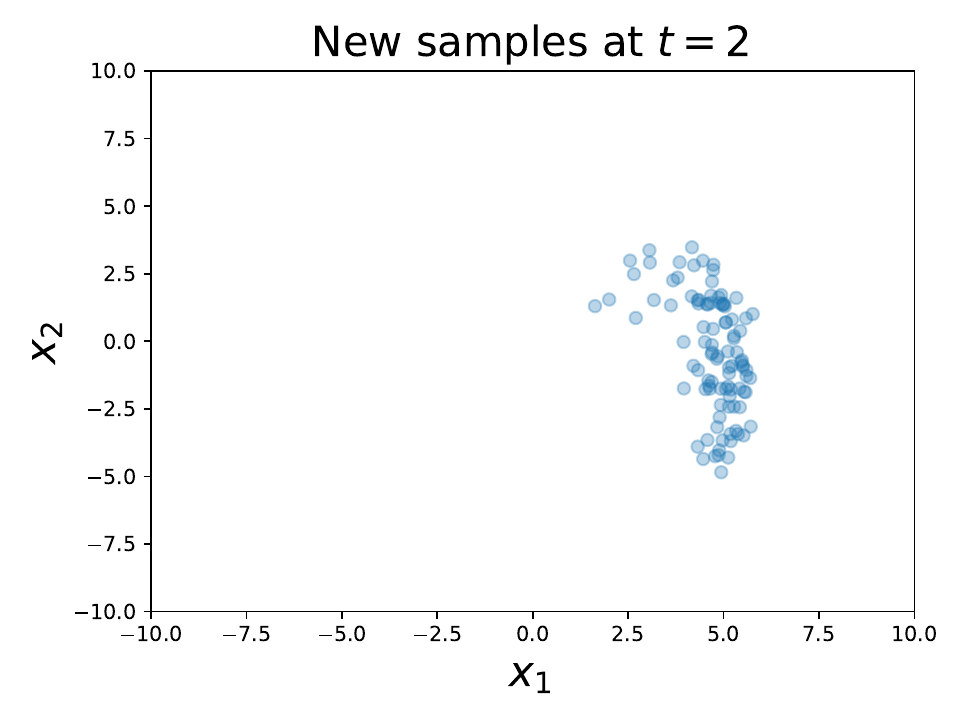}
        \includegraphics[width=\linewidth]{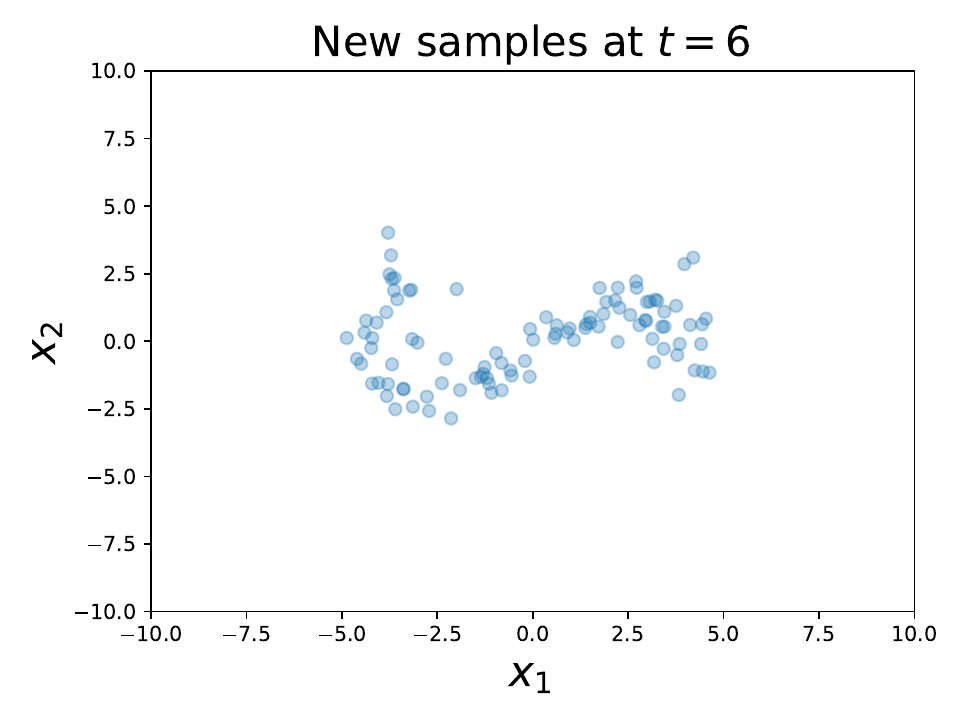}
        \includegraphics[width=\linewidth]{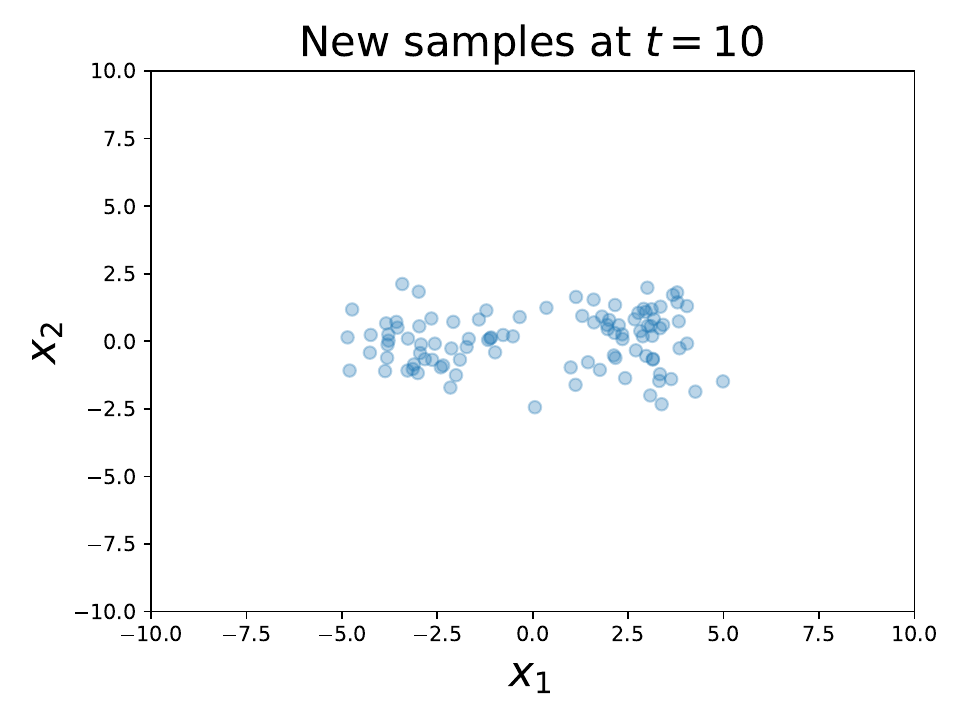}
        \caption{IR-PINNs2}
    \end{subfigure}

    \caption{\textit{Time-dependent Fokker-Planck equation:} New samples generated by probability density model $p(t, \bm{x}; \theta_f)$ of the last adaptive iteration at $t = 2, 6, 10$.}
    \label{fig:FPE_adaptive_points}
\end{figure}

We present the reference solutions, numerical solutions and absolute errors of different methods at $t = 2,6,10$ in Figure \ref{fig:FPE_reference_solution}-\ref{fig:FPE_solution_IR-PINNs2}. The results demonstrate that our method can be directly combined with an adaptive sampling strategy for solving time-dependent Fokker-Planck equation using temporal normalized flows. Similar to the previous example, due to the difficulties caused by the evolving physical localization, the improvements from temporal discretization are limited, as shown in Figure \ref{fig:FPE_error_KL} and Table \ref{table:FPE} in terms of the relative $L_2$ error and relative KL divergence. In Figure \ref{fig:FPE_adaptive_points}, we present the adaptively generated samples of the last adaptive iteration at $t = 2, 6, 10$.

\section{Conclusion}\label{sec:conclusion}
In this paper, we have proposed integral regularization physics-informed neural networks (IR-PINNs) to address the challenges of solving evolution equations, particularly in capturing long-time dynamics and reducing temporal error accumulation. By reformulating the equations into an integral form and incorporating an integral-based residual term, IR-PINNs enhance temporal accuracy and improve the resolution of temporal correlations. The addition of an adaptive sampling strategy further ensures higher accuracy in regions with sharp gradients or rapid variations. Numerical experiments on benchmark problems have demonstrated the effectiveness of IR-PINNs. Future research directions include investigating the underlying mechanisms of how the integral-based regularization term improves accuracy and determining its applicability to a broader range of equations. Additionally, exploring more effective adaptive sampling strategies could further enhance computational efficiency and solution precision. In summary, IR-PINNs offer a robust and accurate framework for solving evolution equations, advancing the field of physics-informed machine learning and its applications to dynamic systems.

\section*{Acknowledgments}
We would like to thank Professor Tao Zhou of the Chinese Academy of Sciences for valuable discussion. The third author is partially supported by Zhuhai Innovation and Entrepreneurship Team Project (2120004000498).

\bibliographystyle{plain}
\bibliography{bibliography.bib}

\begin{thebibliography}{10}

\bibitem{bradbury2018jax}
James Bradbury, Roy Frostig, Peter Hawkins, Matthew~James Johnson, Chris Leary,
  Dougal Maclaurin, George Necula, Adam Paszke, Jake VanderPlas, Skye
  Wanderman-Milne, et~al.
\newblock {JAX}: composable transformations of {P}ython+ {N}umpy programs.
\newblock 2018.

\bibitem{bruno2022fc}
Oscar~P Bruno, Jan~S Hesthaven, and Daniel~V Leibovici.
\newblock {FC}-based shock-dynamics solver with neural-network localized
  artificial-viscosity assignment.
\newblock {\em Journal of Computational Physics: X}, 15:100110, 2022.

\bibitem{chiu2022can}
Pao-Hsiung Chiu, Jian~Cheng Wong, Chinchun Ooi, My~Ha Dao, and Yew-Soon Ong.
\newblock {CAN-PINN}: A fast physics-informed neural network based on
  coupled-automatic-numerical differentiation method.
\newblock {\em Computer Methods in Applied Mechanics and Engineering},
  395:114909, 2022.

\bibitem{courant2008methods}
Richard Courant and David Hilbert.
\newblock {\em Methods of mathematical physics: partial differential
  equations}.
\newblock John Wiley \& Sons, 2008.

\bibitem{daw2022rethinking}
Arka Daw, Jie Bu, Sifan Wang, Paris Perdikaris, and Anuj Karpatne.
\newblock Rethinking the importance of sampling in physics-informed neural
  networks.
\newblock {\em arXiv preprint arXiv:2207.02338}, 2022.

\bibitem{dong2021method}
Suchuan Dong and Naxian Ni.
\newblock A method for representing periodic functions and enforcing exactly
  periodic boundary conditions with deep neural networks.
\newblock {\em Journal of Computational Physics}, 435:110242, 2021.

\bibitem{dresdner2022learning}
Gideon Dresdner, Dmitrii Kochkov, Peter Norgaard, Leonardo
  Zepeda-N{\'u}{\~n}ez, Jamie~A Smith, Michael~P Brenner, and Stephan Hoyer.
\newblock Learning to correct spectral methods for simulating turbulent flows.
\newblock {\em arXiv preprint arXiv:2207.00556}, 2022.

\bibitem{driscoll2014chebfun}
Tobin~A Driscoll, Nicholas Hale, and Lloyd~N Trefethen.
\newblock Chebfun guide, 2014.

\bibitem{feng2024hybrid}
Xiaodong Feng, Haojiong Shangguan, Tao Tang, Xiaoliang Wan, and Tao Zhou.
\newblock A hybrid {FEM-PINN} method for time-dependent partial differential
  equations.
\newblock {\em arXiv preprint arXiv:2409.02810}, 2024.

\bibitem{feng2022solving}
Xiaodong Feng, Li~Zeng, and Tao Zhou.
\newblock Solving time dependent {F}okker-{P}lanck equations via temporal
  normalizing flow.
\newblock {\em Communications in Computational Physics}, 32(2):401--423, 2022.

\bibitem{glorot2010understanding}
Xavier Glorot and Yoshua Bengio.
\newblock Understanding the difficulty of training deep feedforward neural
  networks.
\newblock In {\em Proceedings of the thirteenth international conference on
  artificial intelligence and statistics}, pages 249--256. JMLR Workshop and
  Conference Proceedings, 2010.

\bibitem{gu2022deep}
Yiqi Gu and Michael~K Ng.
\newblock Deep adaptive basis {G}alerkin method for high-dimensional evolution
  equations with oscillatory solutions.
\newblock {\em SIAM Journal on Scientific Computing}, 44(5):A3130--A3157, 2022.

\bibitem{jung2024ceens}
Jeahan Jung, Heechang Kim, Hyomin Shin, and Minseok Choi.
\newblock {CEEN}s: Causality-enforced evolutional networks for solving
  time-dependent partial differential equations.
\newblock {\em Computer Methods in Applied Mechanics and Engineering},
  427:117036, 2024.

\bibitem{kingma2014adam}
Diederik~P Kingma and Jimmy Ba.
\newblock Adam: A method for stochastic optimization.
\newblock {\em arXiv preprint arXiv:1412.6980}, 2014.

\bibitem{kloeden1992stochastic}
Peter~E Kloeden, Eckhard Platen, Peter~E Kloeden, and Eckhard Platen.
\newblock {\em Stochastic differential equations}.
\newblock Springer, 1992.

\bibitem{krishnapriyan2021characterizing}
Aditi Krishnapriyan, Amir Gholami, Shandian Zhe, Robert Kirby, and Michael~W
  Mahoney.
\newblock Characterizing possible failure modes in physics-informed neural
  networks.
\newblock {\em Advances in Neural Information Processing Systems},
  34:26548--26560, 2021.

\bibitem{lippe2023pde}
Phillip Lippe, Bas Veeling, Paris Perdikaris, Richard Turner, and Johannes
  Brandstetter.
\newblock Pde-refiner: Achieving accurate long rollouts with neural pde
  solvers.
\newblock {\em Advances in Neural Information Processing Systems},
  36:67398--67433, 2023.

\bibitem{lu2022learning}
Yubin Lu, Romit Maulik, Ting Gao, Felix Dietrich, Ioannis~G Kevrekidis, and
  Jinqiao Duan.
\newblock Learning the temporal evolution of multivariate densities via
  normalizing flows.
\newblock {\em Chaos: An Interdisciplinary Journal of Nonlinear Science},
  32(3), 2022.

\bibitem{mattey2022novel}
Revanth Mattey and Susanta Ghosh.
\newblock A novel sequential method to train physics informed neural networks
  for allen cahn and cahn hilliard equations.
\newblock {\em Computer Methods in Applied Mechanics and Engineering},
  390:114474, 2022.

\bibitem{muller2019neural}
Thomas M{\"u}ller, Brian McWilliams, Fabrice Rousselle, Markus Gross, and Jan
  Nov{\'a}k.
\newblock Neural importance sampling.
\newblock {\em ACM Transactions on Graphics (ToG)}, 38(5):1--19, 2019.

\bibitem{penwarden2023unified}
Michael Penwarden, Ameya~D Jagtap, Shandian Zhe, George~Em Karniadakis, and
  Robert~M Kirby.
\newblock A unified scalable framework for causal sweeping strategies for
  physics-informed neural networks ({PINNs}) and their temporal decompositions.
\newblock {\em Journal of Computational Physics}, 493:112464, 2023.

\bibitem{pichler2013numerical}
Lukas Pichler, Arif Masud, and Lawrence~A Bergman.
\newblock Numerical solution of the {Fokker--Planck} equation by finite
  difference and finite element methods—a comparative study.
\newblock {\em Computational Methods in Stochastic Dynamics: Volume 2}, pages
  69--85, 2013.

\bibitem{raissi2019physics}
Maziar Raissi, Paris Perdikaris, and George~E Karniadakis.
\newblock Physics-informed neural networks: A deep learning framework for
  solving forward and inverse problems involving nonlinear partial differential
  equations.
\newblock {\em Journal of Computational Physics}, 378:686--707, 2019.

\bibitem{rathore2024challenges}
Pratik Rathore, Weimu Lei, Zachary Frangella, Lu~Lu, and Madeleine Udell.
\newblock Challenges in training {PINNs}: A loss landscape perspective.
\newblock {\em arXiv preprint arXiv:2402.01868}, 2024.

\bibitem{wang2023long}
Sifan Wang and Paris Perdikaris.
\newblock Long-time integration of parametric evolution equations with
  physics-informed deeponets.
\newblock {\em Journal of Computational Physics}, 475:111855, 2023.

\bibitem{wang2024respecting}
Sifan Wang, Shyam Sankaran, and Paris Perdikaris.
\newblock Respecting causality for training physics-informed neural networks.
\newblock {\em Computer Methods in Applied Mechanics and Engineering},
  421:116813, 2024.

\bibitem{wang2021understanding}
Sifan Wang, Yujun Teng, and Paris Perdikaris.
\newblock Understanding and mitigating gradient flow pathologies in
  physics-informed neural networks.
\newblock {\em SIAM Journal on Scientific Computing}, 43(5):A3055--A3081, 2021.

\bibitem{wang2021learning}
Sifan Wang, Hanwen Wang, and Paris Perdikaris.
\newblock Learning the solution operator of parametric partial differential
  equations with physics-informed {DeepONets}.
\newblock {\em Science Advances}, 7(40):eabi8605, 2021.

\bibitem{wang2022and}
Sifan Wang, Xinling Yu, and Paris Perdikaris.
\newblock When and why {PINNs} fail to train: A neural tangent kernel
  perspective.
\newblock {\em Journal of Computational Physics}, 449:110768, 2022.

\bibitem{wight2020solving}
Colby~L Wight and Jia Zhao.
\newblock Solving {A}llen-{C}ahn and {C}ahn-{H}illiard equations using the
  adaptive physics informed neural networks.
\newblock {\em Communications in Computational Physics}, 29(3):930--954, 2021.

\bibitem{wu2025coast}
Zhikai Wu, Shiyang Zhang, Sizhuang He, Sifan Wang, Min Zhu, Anran Jiao, Lu~Lu,
  and David van Dijk.
\newblock Coast: Intelligent time-adaptive neural operators.
\newblock {\em arXiv preprint arXiv:2502.08574}, 2025.

\bibitem{zeng2023bounded}
Li~Zeng, Xiaoliang Wan, and Tao Zhou.
\newblock Bounded {KR}net and its applications to density estimation and
  approximation.
\newblock {\em arXiv preprint arXiv:2305.09063}, 2023.

\end{thebibliography}

\appendix
\section{Time-marching strategy}\label{sec:time_marching}
Following the methodology proposed in \cite{wight2020solving}, we implement a time-marching strategy to enhance the convergence of IR-PINNs when applied to long-time integration problems. Specifically, we partition the temporal domain $[0, T]$ into subdomains
\begin{equation}
    [0, \Delta t], [\Delta t, 2\Delta t], \cdots, [(N-1)\Delta t, N\Delta t], \quad \Delta t = T/N.
\end{equation}
We build a single neural network for each subdomain, training consecutively by using the final prediction of the current subdomain as the initial condition for the next, until the entire temporal domain is fully trained.

\section{Exact periodic boundary conditions}\label{sec:periodic_BC}
Following the work from \cite{dong2021method, wang2024respecting}, we can enforce exact $C^{\infty}$ periodic boundary conditions by constructing a Fourier feature embedding of the form
\begin{equation}
    v(x) = \left\{ 1,\cos(\omega x),\sin(\omega x) ,\cdots,\cos(M\omega x),\sin(M\omega x)\right\}
\end{equation}
as the spatial input to the neural network, where $\omega =\frac{2\pi}{L}, L = x_{\mathrm{max}} - x_{\mathrm{min}}$, and $M$ is a non-negative integer representing the frequency of the input. In this work, we take $M = 5$.

\section{Pre-training approach for the time-dependent Fokker-Planck equation}\label{sec:pre-training}
Assuming we can obtain some sample paths directly from the SDE \eqref{Fokker-Planck_SDE} using Euler-Maruyama method \cite{kloeden1992stochastic}, the accuracy of the numerical solution will be further improved. Specifically, for a given set of $N_t$ temporal points $\{t_j\}_{j=1}^{N_t}$, we generate $N_{p}$ paths $\{X_t^{(i)}\}_{i=1}^{N_p}$:
\begin{equation}
    X_t^{(i)} = \left((t_1, \bm{x}_1^{(i)}), (t_2, \bm{x}_2^{(i)}), \cdots, (t_{N_t}, \bm{x}_{N_t}^{(i)})\right), \quad i = 1, 2, \cdots, N_p.
\end{equation}
Subsequently, the neural network can be pre-trained by minimizing the negative log-likelihood loss \cite{lu2022learning}, formulated as:
\begin{equation}
    \mathcal{L}_{\mathrm{data}} = -\frac{1}{N_p N_t} \sum_{i=1}^{N_p} \sum_{j=1}^{N_t} \log p(t_j, \bm{x}_j^{(i)}).
\end{equation}

\typeout{get arXiv to do 4 passes: Label(s) may have changed. Rerun}

\end{document}